\documentclass[12pt]{article}
\usepackage{amsfonts}
\usepackage{full page,
amssymb, amscd, amsmath, graphicx, %showkeys
}

\makeatletter
\DeclareFontFamily{U}{tipa}{}
\DeclareFontShape{U}{tipa}{m}{n}{<->tipa10}{}
\newcommand{\arc@char}{{\usefont{U}{tipa}{m}{n}\symbol{62}}}%

\newcommand{\arc}[1]{\mathpalette\arc@arc{#1}}

\newcommand{\arc@arc}[2]{%
  \sbox0{$\m@th#1#2$}%
  \vbox{
    \hbox{\resizebox{\wd0}{\height}{\arc@char}}
    \nointerlineskip
    \box0
  }%
}
\makeatother

\makeatletter
\def\widebreve{\mathpalette\wide@breve}
\def\wide@breve#1#2{\sbox\z@{$#1#2$}%
     \mathop{\vbox{\m@th\ialign{##\crcr
\kern0.08em\brevefill#1{0.8\wd\z@}\crcr\noalign{\nointerlineskip}%
                    $\hss#1#2\hss$\crcr}}}\limits}
\def\brevefill#1#2{$\m@th\sbox\tw@{$#1($}%
  \hss\resizebox{#2}{\wd\tw@}{\rotatebox[origin=c]{90}{\upshape(}}\hss$}
\makeatletter

\newtheorem{thm}{Theorem}[section]
\newtheorem{defn}[thm]{Definition}
\newtheorem{prop}[thm]{Proposition}
\newtheorem{cor}[thm]{Corollary}
\newtheorem{lemma}[thm]{Lemma}
\newtheorem{rema}[thm]{Remark}

\newcommand{\halmos}{\rule{1ex}{1.4ex}}

\newcommand{\nn}{\nonumber \\}
\newcommand{\overarc}[1]{\,\arc{\!{#1}}}
 
 \newcommand{\res}{\mbox{\rm Res}}

\renewcommand{\hom}{\mbox{\rm Hom}}

 \newcommand{\pf}{{\it Proof.}\hspace{2ex}}
 
 \newcommand{\epfv}{\hspace*{\fill}\mbox{$\halmos$}\vspace{1em}}

\newcommand{\wt}{\mbox{\rm wt}\,}

\newcommand{\C}{\mathbb{C}}
\newcommand{\Z}{\mathbb{Z}}
\newcommand{\R}{\mathbb{R}}

\newcommand{\N}{\mathbb{N}}

\newcommand{\I}{\mathbb{I}}

\newcommand{\one}{\mathbf{1}}

\renewcommand{\mod}{\;\;\mbox{\rm mod}\;}
\newcommand{\g}{\mathfrak{g}}

\title{ {\bf Lower-bounded and grading-restricted twisted modules
for affine vertex (operator) algebras} }
\date{}
\author{Yi-Zhi Huang}
\begin{document}

\bibliographystyle{alpha}
\maketitle
\begin{abstract}
We apply the construction of the universal lower-bounded generalized twisted modules 
by the author to construct  universal 
lower-bounded and grading-restricted generalized twisted modules 
for affine vertex (operator) algebras.
We prove that these universal twisted modules for affine vertex (operator) algebras
are equivalent to suitable induced modules of the corresponding twisted affine Lie algebra
or quotients of such induced modules by explicitly given submodules. 
\end{abstract}

\renewcommand{\theequation}{\thesection.\arabic{equation}}
\renewcommand{\thethm}{\thesection.\arabic{thm}}
\setcounter{equation}{0}
\setcounter{thm}{0}
\section{Introduction}

In \cite{H-const-twisted-mod}, the author constructed universal 
lower-bounded generalized twisted modules for a grading-restricted vertex algebra. 
In the present paper, we apply this construction to construct and identify explicitly universal 
lower-bounded and grading-restricted generalized twisted modules 
for affine vertex (operator) algebras. In particular,  general classes of 
lower-bounded and grading-restricted generalized twisted modules can be studied using these universal ones. 

Let $\g$ be a finite-dimensional Lie algebra with a nondegenerate invariant symmetric bilinear form $(\cdot, \cdot)$
and $g$ an automorphism of $\g$. Then an induced module $M(\ell, 0)$ of level $\ell\in \C$ 
for the affine Lie algebra $\hat{\g}$ generated by the trivial module $\C$ for $\g$ 
has a structure of grading-restricted vertex algebra. In the case that $\g$ is simple and 
$\ell\ne -h^{\vee}$, where $h^{\vee}$ is the dual Coxeter number of $\g$, $M(\ell, 0)$ has 
a conformal vector and is thus a vertex operator algebra. Let $L(\ell, 0)$ be the 
irreducible quotient of $M(\ell, 0)$. Then $L(\ell, 0)$ is also a graidng-restricted vertex algebra
and, when $\ell+h^{\vee}\ne 0$, is a vertex operator algebra. 
An automorphism $g$ of $\g$ induces  automorphisms, still denoted by $g$, of 
$\hat{\g}$, $M(\ell, 0)$ and $L(\ell, 0)$. There is also a twisted affine Lie algebra $\hat{\g}^{[g]}$
constructed using $\g$, $(\cdot, \cdot)$ and $g$. Note that $\g$ has a rich automorphism group
containing the Lie group corresponding to $\g$. Automorphisms
of $\g$, $M(\ell, 0)$ and $L(\ell, 0)$
are mostly of infinite orders and many of them do not act on $\g$, $M(\ell, 0)$ and $L(\ell, 0)$ semisimply. 

Twisted modules associated to automorphisms of finite orders of a vertex operator algebra
were introduced and studied first 
by Frenkel, Lepowsky and Meurman in \cite{FLM1}, 
\cite{FLM2} and
\cite{FLM3} and by Lepowsky in \cite{Le1} and \cite{Le2}.  
In \cite{H-log-twisted-mod}, the author introduced twisted modules associated to
general automorphisms of a vertex operator algebra, including in particular, automorphisms
which do not act on the vertex operator algebra semisimply. A particular class of 
examples associated to such general automorphisms were also given in \cite{H-log-twisted-mod}.
In \cite{H-const-twisted-mod}, the author gave a construction of universal 
lower-bounded generalized twisted modules associated to such general automorphisms of
 a grading-restricted vertex algebra. 

Applying the construction in \cite{H-const-twisted-mod} to $M(\ell, 0)$, we construct universal 
 lower-bounded (grading-restricted) generalized $g$-twisted $M(\ell, 0)$-modules
generated by a vector space (a finite-dimensional module for a  suitable subalgebra $\hat{\g}^{[g]}_{\I}$ 
of $\hat{\g}^{[g]}$)  with actions of $g$, its semisimple and unipotent parts,
and some other operators and annihilated by the positive part of $\hat{\g}^{[g]}$
when $M(\ell, 0)$ is viewed as a grading-restricted vertex algebra. When 
$\ell+h^{\vee}\ne 0$ and $M(\ell, 0)$ is viewed as a vertex operator algebra, we also construct universal 
 lower-bounded (grading-restricted) generalized $g$-twisted $M(\ell, 0)$-modules
generated by a vector space (a finite-dimensional module for $\hat{\g}^{[g]}_{\I}$) 
with additional structures as above. These universal  lower-bounded and grading-restricted generalized 
$g$-twisted $M(\ell, 0)$-modules as  $\hat{\g}^{[g]}$-modules are then proved to be equivalent 
to suitable induced modules for $\hat{\g}^{[g]}$. The proofs of these equivalences
use the results in Section 2 of \cite{H-exist-twisted-mod} on a linearly independent 
set of generators of the universal lower-bounded generalized twisted modules 
constructed in Section 5 of \cite{H-const-twisted-mod}. In the case that $M(\ell, 0)$ is viewed as 
a vertex operator algebra, we also give explicit formulas for the Virasoro operators on 
the universal lower-bounded generalized twisted $M(\ell, 0)$-modules. These formulas are 
needed in the proof of their equivalences to  suitable induced module for $\hat{\g}^{[g]}$. 

When $\g$ is simple and $\ell\in \Z_{+}$ and $L(\ell, 0)$ is viewed as a vertex operator algebra, 
we construct universal  lower-bounded (grading-restricted) generalized $g$-twisted 
$M(\ell, 0)$-modules generated by a vector space (a finite-dimensional $\hat{\g}^{[g]}_{\I}$-module) with additional structures as above. 
We also prove that these universal  lower-bounded and grading-restricted 
generalized $g$-twisted $L(\ell, 0)$-modules as  $\hat{\g}^{[g]}$-modules 
are equivalent to quotients by  certain explicitly given submodules of the 
induced modules for $\hat{\g}^{[g]}$ equivalent 
to the universal lower-bounded and grading-restricted generalized $g$-twisted $M(\ell, 0)$-modules. 
To prove these equivalences, 
we also generalize a result of Kac (see Proposition 8.1 in \cite{K})
on automorphisms of finite orders of a finite-dimensional simple Lie algebra to 
semisimple automorphisms of arbitrary orders. 

Immediate consequences of the universal properties satisfied by those universal twisted 
modules are that 
 lower-bounded and grading-restricted generalized $g$-twisted 
$M(\ell, 0)$- and $L(\ell, 0)$-modules generated by subspaces and finite-dimensional 
$\hat{\g}^{[g]}_{\I}$-submodules with additional structures as above are quotients 
of these universal ones. Thus we can study these types of twisted modules, including untwisted ones,  
using our results on the universal ones in the present paper. 

In the case that $g$ is of finite order, 
Li gave the relationship between weak twisted modules for an affine vertex operator algebra
and  restricted modules for the corresponding twisted affine Lie algebra  in \cite{Li2}. 
In \cite{B}, Bakalov introduced twisted affine Lie algebras in the case that 
$g$ is a general automorphsim of $\g$ and gave the relationship 
between weak twisted modules for an affine vertex operator algebra
and  restricted modules for the corresponding twisted affine Lie algebra.
In the present paper, we do not study these most general weak twisted modules and restricted modules. 
We study only lower-bounded and grading-restricted generalized twisted modules for affine vertex (operator) 
algebras and lower-bounded and grading-restricted modules for the twisted affine Lie algebras. 
We would like to emphasize that in the representation theory of vertex (operator) algebras,
to obtain substantial results, we have to restrict ourselves to grading-restricted 
generalized (twisted) modules and we often have to further restrict ourselves to 
such modules of finite lengths. On the other hand, lower-bounded 
generalized (twisted) modules always appear in various constructions and proofs.
One of the difficult problems is to prove that these  lower-bounded 
generalized (twisted) modules appearing in our constructions and proofs are actually 
grading-restricted 
generalized (twisted) modules of finite lengths. So these two types of twisted modules 
are what we are mainly interested. Moreover, for such modules of finite lengths, we can reduce 
their study to those modules generated by subspaces annihilated by the positive part of 
the twisted affine Lie algebra. This is the reason why we choose to construct, identify 
and study  these types of twisted modules in this paper. Though weak twisted modules are more general, 
usually we need them only in the formulations of certain notions in the representation 
theory of vertex (operator) algebras.

As is mentioned in the preceding paragraph, one of the difficult problems in the representation 
theory of vertex operator algebras is to prove that suitable  lower-bounded 
generalized twisted modules are actually grading restricted. In fact,  
universal lower-bounded generalized twisted modules are in a certain sense analogous to 
Verma modules in the representation theory of finite-dimensional Lie algebras. 
In the case of finite-dimensional Lie algebras, we know that a Verma module generated from 
a highest weight vector has a finite-dimensional quotient module if and only if the highest weight is 
dominant integral. For vertex operator algebra, we can ask an analogous question: 
Under what conditions, a universal lower-bounded generalized twisted module has a grading-restricted 
quotient. In this paper, our construction and identification of 
grading-restricted generalized twisted modules for $M(\ell, 0)$ and $L(\ell, 0)$
gives an answer to this question for affine vertex (operator) algebras. 

One main goal of studying these twisted modules is to use their properties and structures 
to study twisted intertwining operators among them (see \cite{H-twisted-int}). 
We expect that the constructions and results in
the present paper will play an important role in the study of twisted intertwining operators 
for affine vertex operator algebras. 

The present paper is organized as follows: In Section 2, we recall some basic material 
on  the affine Lie algebra $\hat{\g}$ of a finite-dimensional Lie algebra $\g$,
an automorphism $g$  of $\g$ and the twisted affine Lie algebra $\hat{\g}^{[g]}$. In Section 3, we recall 
vertex (operator) algebras $M(\ell, 0)$ and $L(\ell, 0)$ associated to affine Lie algebras and their automorphisms
induced from those of $\g$. The construction, 
identification and basic properties of lower-bounded and grading-restricted generalized twisted modules for 
$M(\ell, 0)$ are given in Section 4. In Subsection \ref{4.1}, we construct and identify explicitly 
 lower-bounded and grading-restricted generalized twisted modules for 
$M(\ell, 0)$ viewed as a grading-restricted vertex algebra. 
In Subsection \ref{4.2}, we construct and identify explicitly such twisted modules for 
$M(\ell, 0)$ viewed as a vertex operator algebra. In Subsection \ref{4.3}, basic properties of 
these twisted modules for $M(\ell, 0)$, including their universal properties and their quotients, are given. 
The construction, 
identification and  basic properties of lower-bounded and grading-restricted generalized twisted modules for 
$L(\ell, 0)$ are given in Section 5. 

\paragraph{Acknowledgments}
The author is grateful to Sven M\"{o}ller for the argument (Lemma 8.3 in \cite{EMS})
needed in the last step of the proof of Proposition \ref{autom-decomp}. 

\renewcommand{\theequation}{\thesection.\arabic{equation}}
\renewcommand{\thethm}{\thesection.\arabic{thm}}
\setcounter{equation}{0}
\setcounter{thm}{0}
\section{Twisted affine Lie algebras}

Let $\mathfrak{g}$ be  a finite-dimensional  Lie algebra
and  $(\cdot, \cdot)$ a nondegenerate invariant symmetric bilinear form on $\mathfrak{g}$.
Recall that the affine Lie algebra $\hat{\mathfrak{g}}$ is
the vector space $\mathfrak{g} \otimes \C[t, t^{-1}] \oplus \C{\bf k}$ equipped with the bracket operation
\begin{align*}
[a \otimes t^m, b \otimes t^n] &= [a, b]\otimes t^{m+n} + (a, b)m\delta_{m+n,0}{\bf k},\\
[a \otimes t^m, {\bf k}] &= 0,
\end{align*}
for $a, b \in \mathfrak{g}$ and $m, n \in \Z$. 
Let $\hat{\mathfrak{g}}_{\pm} = \mathfrak{g}\otimes t^{\pm 1}\C[t^{\pm 1}]$. Then
$$\hat{\mathfrak{g}} = \hat{\mathfrak{g}}_{-} \oplus \mathfrak{g} \oplus \C{\bf k} \oplus \hat{\mathfrak{g}}_{+}.$$

Let $g$ be an automorphism of $\mathfrak{g}$. Assume also that 
$(\cdot, \cdot)$ is invariant under $g$. This is true in the case that $\g$ is semisimple and 
$(\cdot, \cdot)$ is proportional to the Killing form. 
Since $\mathfrak{g}$ is finite dimensional,
there exist operators $\mathcal{L}_{g}$, $\mathcal{S}_{g}$ and
 $\mathcal{N}_{g}$ on $\mathfrak{g}$
such that $g=e^{2\pi i\mathcal{L}_{g}}$ and $\mathcal{S}_{g}$ and
 $\mathcal{N}_{g}$ are the semisimple and nilpotent parts of $\mathcal{L}_{g}$, 
respectively. Then $g$, $\mathcal{L}_{g}$, $\mathcal{S}_{g}$ and
 $\mathcal{N}_{g}$ induce operators,
still denoted by $g$, $\mathcal{L}_{g}$, $\mathcal{S}_{g}$ and
 $\mathcal{N}_{g}$, on
the affine Lie algebra $\hat{\mathfrak{g}}$ such that $g$ is also an automorphism
of $\hat{\mathfrak{g}}$. 

Let 
$$P_{\g}=\{\alpha\in \C\;|\; \Re(\alpha)\in [0, 1), 
e^{2 \pi i\alpha} \;\text{is an eigenvalue of}\;g\}.$$
Then 
$$\mathfrak{g}=\coprod_{\alpha\in P_{\g}}
\mathfrak{g}^{[\alpha]},$$
where for $\alpha\in P_{\g}$,
$\mathfrak{g}^{[\alpha]}$ is the generalized eigenspace of $g$ (or the eigenspace of 
$e^{2\pi i\mathcal{S}_{g}}$) with the eigenvalue $e^{2\pi i\alpha}$.

%Let $P_{K}(t)$ 
%be the subspace of all polynomials in $t$ of degree less thsan or equal to $K$. 

For $\alpha, \beta\in [0, 1)+i\R$,  let 
$$s(\alpha, \beta)=\left\{\begin{array}{ll}\alpha+\beta&\Re(\alpha+\beta)<1\\
\alpha+\beta-1&\Re(\alpha+\beta)\ge 1.\end{array}\right.$$
Then $\Re(s(\alpha, \beta))\in [0, 1)$ and for $\alpha, \beta\in P_{\g}$,
 $s(\alpha, \beta)\in (P_{g}+P_{g})\cup (P_{g}+P_{g}-1)$. 

\begin{lemma}\label{brac-eigenvalue}
For $\alpha, \beta\in P_{\g}$, 
$[\mathfrak{g}^{[\alpha]}, \mathfrak{g}^{[\beta]}]\subset 
\mathfrak{g}^{[s(\alpha, \beta)]}$.
In particular,  in the case that  
$[\mathfrak{g}^{[\alpha]}, \mathfrak{g}^{[\beta]}]\ne 0$, $s(\alpha, \beta)\in P_{g}\cap ((P_{g}+P_{g})\cup (P_{g}+P_{g}-1))\subset
P_{g}$ and $e^{2\pi i(\alpha+\beta)}$ is an eigenvalue of $g$. 
\end{lemma}
\pf
For $a\in \mathfrak{g}^{[\alpha]}$ 
and $b\in \mathfrak{g}^{[\beta]}$,
we have 
\begin{align*}
(g-&e^{2\pi i(\alpha+\beta)})[a, b]\nn
&=
[ga, gb]-e^{2\pi i(\alpha+\beta)}[a, b]\nn
&=e^{2\pi i(\alpha+\beta)}[e^{2\pi i\mathcal{N}_{g}}a, e^{2\pi i\mathcal{N}_{g}}b]
-e^{2\pi i(\alpha+\beta)}[a, b]\nn
&=e^{2\pi i(\alpha+\beta)}\left([(1_{\mathfrak{g}}
+(e^{2\pi i\mathcal{N}_{g}}-1_{\mathfrak{g}}))a, (1_{\mathfrak{g}}
+(e^{2\pi i\mathcal{N}_{g}}-1_{\mathfrak{g}}))b]
-[a, b]\right)\nn
&=e^{2\pi i(\alpha+\beta)}\left([(e^{2\pi i\mathcal{N}_{g}}-1_{\mathfrak{g}})a, b]
+[a, (e^{2\pi i\mathcal{N}_{g}}-1_{\mathfrak{g}})b]+[(e^{2\pi i\mathcal{N}_{g}}-1_{\mathfrak{g}})a, (e^{2\pi i\mathcal{N}_{g}}-1_{\mathfrak{g}})b]\right).
\end{align*}
Then there exists $\widetilde{K}\in \Z_{+}$ such that 
$$(g-e^{2\pi i(\alpha+\beta)})^{\widetilde{K}}[a, b]=0.$$ 
(Note that we can always take $\widetilde{K}=\dim \g$.)
 If $[a, b]=0$, we have $[a, b]\subset \mathfrak{g}^{[s(\alpha, \beta)]}$.
If $[a, b]\ne 0$, it is a generalized eigenvector of $g$ with 
eigenvalue $e^{2\pi i(\alpha+\beta)}$ and thus is in $\mathfrak{g}^{[s(\alpha, \beta)]}$.
In this case, $\mathfrak{g}^{[s(\alpha, \beta)]}\ne 0$. So $s(\alpha, \beta)\in P_{g}$.
We also have either $s(\alpha, \beta)=\alpha+\beta\in P_{g}+P_{g}$ or 
$s(\alpha, \beta)=\alpha+\beta-1\in P_{g}+P_{g}-1$. Thus 
$s(\alpha, \beta)\in P_{g}\cap ((P_{g}+P_{g})\cup (P_{g}+P_{g}-1))$. 
\epfv

\begin{cor}\label{N-g-der}
The operators $e^{2\pi i\mathcal{S}_{g}}$ and  $e^{2\pi i\mathcal{N}_{g}}$ 
are also automorphisms of $\g$. 
The operator $\mathcal{N}_{g}$ is a derivation of the Lie algebra $\g$. 
\end{cor}
\pf
Let $a\in \g^{[\alpha]}$ and $b\in \g^{[\beta]}$. By Lemma \ref{brac-eigenvalue},
$$e^{2\pi i\mathcal{S}_{g}}[a, b]=e^{2\pi i(\alpha+\beta)}[a, b]
=[e^{2\pi i\alpha}a, e^{2\pi i\alpha}b]
=[e^{2\pi i\mathcal{S}_{g}}a, e^{2\pi i\mathcal{S}_{g}}b].$$
So $e^{2\pi i\mathcal{S}_{g}}$ is an automorphism of $\g$. Therefore 
$e^{-2\pi i\mathcal{S}_{g}}$ is also an automorphism of $\g$. Thus
$e^{2\pi i\mathcal{N}_{g}}=e^{-2\pi i\mathcal{S}_{g}}g$ is an automorphism 
of $\g$. 

For $a, b\in \g$, we have 
\begin{align*}
(\text{ad}\;e^{2\pi i\mathcal{N}_{g}}a)b
&=[e^{2\pi i\mathcal{N}_{g}}a, b]\nn
&=e^{2\pi i\mathcal{N}_{g}} [a, e^{-2\pi i\mathcal{N}_{g}}b]\nn
&=e^{2\pi i\mathcal{N}_{g}}(\text{ad}\;a) e^{-2\pi i\mathcal{N}_{g}}b\nn
&=((\text{Ad}\;2\pi i\mathcal{N}_{g})(\text{ad}\;a))b.
\end{align*}
Thus 
\begin{align*}
[\mathcal{N}_{g}a, b]&=\frac{1}{2\pi i}
(\text{ad}\; \log e^{2\pi i\mathcal{N}_{g}}a)b\nn
&=\frac{1}{2\pi i} (\log (\text{Ad}\;2\pi i\mathcal{N}_{g})(\text{ad}\;a))b\nn
&=((\text{ad}\; \mathcal{N}_{g})(\text{ad}\;a))b\nn
&=\mathcal{N}_{g}[a, b]-[a, \mathcal{N}_{g}b],
\end{align*}
proving that $\mathcal{N}_{g}$ is a derivation of $\g$.
\epfv

\begin{lemma}\label{S-g-ortho}
If $\alpha+ \beta\not \in \{ 0, 1\}$, then $\g^{[\alpha]}$ and $\g^{[\beta]}$ are orthogonal.
If $\alpha+ \beta\in \{0, 1\}$, then $(\cdot, \cdot)$ restricted to $\g^{[\alpha]}\times \g^{[\beta]}$
is nondegenerate. 
\end{lemma}
\pf
For $a\in \g^{[\alpha]}$, there exists $p\in \Z_{+}$ such that $(g-e^{2\pi i\alpha})^{p}a=0$. 
On the other hand, since $\alpha+ \beta\not \in \{0,  1\}$, 
the restriction of $(g^{-1}-e^{2\pi i\alpha})^{p}$ to $\g^{[\beta]}$ 
is a linear isomorphism from $\g^{[\beta]}$ to itself. If there exist 
$a\in \g^{[\alpha]}$ and $b\in \g^{[\beta]}$ such that $(a, b)\ne 0$. 
then there exists $c\in \g^{[\beta]}$ such that $b=(g^{-1}-e^{2\pi i\alpha})^{p}c$.
Then $(a, (g^{-1}-e^{2\pi i\alpha})^{p}c)=(a, b)\ne 0$. But 
since $(\cdot, \cdot)$ is invariant under $g$, we have 
$(a, (g^{-1}-e^{2\pi i\alpha})^{p}c)=((g-e^{2\pi i\alpha})^{p}a, c)=0$. Contradiction. 
So we must have $(a, b)=0$ for $a\in \g^{[\alpha]}$ and $b\in \g^{[\beta]}$. 

In the case that $\alpha+\beta\in \{0, 1\}$, if $(\cdot, \cdot)$ restricted to $\g^{[\alpha]}\times \g^{[\beta]}$
is degenerate, then there exists $a\in \g^{[\alpha]}\setminus \{0\}$ such that $(a, b)=0$ for $b\in \g^{[\beta]}$.
But for $\beta\in P_{g}$ such that $\alpha+\beta\not\in \{0, 1\}$,
we just proved that $(a, b)=0$ for $b\in \g^{[\beta]}$. 
Thus $(a, b)=0$ for $b\in \g$. Contradiction to the nondegeneracy of $(\cdot, \cdot)$. 
\epfv

\begin{prop}
The nondegenerate invariant symmetric bilinear form $(\cdot, \cdot)$ is also invariant under 
$e^{2\pi i\mathcal{S}_{g}}$ and $e^{2\pi i\mathcal{N}_{g}}$. 
\end{prop}
\pf
Let $a\in \g^{[\alpha]}$ and $b\in \g^{[\beta]}$. If $\alpha+\beta\in \{0, 1\}$, then 
$(e^{2\pi i\mathcal{S}_{g}}a, e^{2\pi i\mathcal{S}_{g}}b)
=e^{2\pi i(\alpha+\beta)}(a, b)=(a, b)$. If $\alpha+\beta\not\in \{ 0, 1\}$, then 
by Lemma \ref{S-g-ortho},
$(e^{2\pi i\mathcal{S}_{g}}a, e^{2\pi i\mathcal{S}_{g}}b)=0=(a, b)$. 
So $(\cdot, \cdot)$ is  invariant under 
$e^{2\pi i\mathcal{S}_{g}}$. 

Since $e^{2\pi i\mathcal{N}_{g}}=e^{-2\pi i\mathcal{S}_{g}} g$ and certainly 
$(\cdot, \cdot)$ is also invariant under $e^{-2\pi i\mathcal{S}_{g}}$,
$(\cdot, \cdot)$ is invariant under 
$e^{2\pi i\mathcal{N}_{g}}$.
\epfv

\begin{cor}\label{N-g-inv}
For $a, b\in \g$, we have $(\mathcal{N}_{g}a, b)+(a, \mathcal{N}_{g}b)=0$. 
\end{cor}
\pf
For $a, b\in \g$, we have 
\begin{align*}
(\mathcal{N}_{g}a, b)&=\left(\frac{1}{2\pi i}(1_{\g}+(\log e^{2\pi i\mathcal{N}_{g}}-1_{\g}))a, b\right)\nn
&=\left(a, \frac{1}{2\pi i}(1_{\g}+(\log e^{-2\pi i\mathcal{N}_{g}}-1_{\g})b\right)\nn
&=-(a, \mathcal{N}_{g}b).
\end{align*}
\epfv

\begin{rema}
{\rm Note that if $\Re\{\alpha\}=\Re\{\beta\}=0$, then
$\Re\{\alpha+\beta\}=\Re\{s(\alpha, \beta)\}=0$. In particular, 
$$\coprod_{\Re\{\alpha\}=0}\mathfrak{g}^{[\alpha]}$$
is a subalgebra of $\mathfrak{g}$. The fixed-point subalgebra 
$\g^{[0]}$ is a subalgebra of this subalgebra.}
\end{rema}

The decomposition 
$$\mathfrak{g}=\coprod_{\alpha\in P_{\g}}
\mathfrak{g}^{[\alpha]}$$
induces decompositions
$$\hat{\mathfrak{g}}=\coprod_{\alpha\in P_{\g}}
\hat{\mathfrak{g}}^{[\alpha]}$$
where $\hat{\mathfrak{g}}^{[\alpha]}$ for $\alpha\in P_{\g}$ 
are the generalized eigenspaces of $g$  (or the eigenspaces of 
$e^{2\pi i\mathcal{S}_{g}}$)  on $\hat{\g}$ 
with the eigenvalue $e^{2\pi i\alpha}$. 

We now define the twisted affine Lie algebra associated to 
$\g$, $(\cdot, \cdot)$ and $g$ (see, for example, \cite{K} and \cite{B}). 
Let 
$$\hat{\mathfrak{g}}^{[g]}=\coprod_{\alpha\in P_{\g}}
\mathfrak{g}^{[\alpha]}\otimes 
t^{\alpha}\C[t, t^{-1}] \oplus \C\mathbf{k}.$$
We define a bracket operation 
on $\hat{\mathfrak{g}}^{[g]}$ by
\begin{align}
[a\otimes t^{m}, b\otimes t^{n}]
&=[a, b]\otimes t^{m+n}
+m(a, b)\delta_{m+n, 0}\mathbf{k}+(\mathcal{N}_{g}a, b)\delta_{m+n, 0}
\mathbf{k},\label{bracket-1}\\
[\mathbf{k}_{1}, a\otimes t^{m}]&=0,\label{bracket-2}
\end{align}
for $a\in \mathfrak{g}^{[\alpha]}$, $b\in 
\mathfrak{g}^{[\beta]}$, $m\in \alpha+\Z$, $n\in 
\beta+\Z$, $\alpha, \beta\in P_{\g}$.
Then it is straightforward to verify that 
the vector space $\hat{\mathfrak{g}}^{[g]}$ equipped with the bracket operation 
defined above is a Lie algebra.

Let 
\begin{align*}
\hat{\mathfrak{g}}^{[g]}_{+}&=\left(\bigoplus_{\alpha\in P_{\g}, 
\Re\{\alpha\}>0}\mathfrak{g}^{[\alpha]}
\otimes t^{\alpha}\C[t]\right)\oplus \left( \bigoplus_{\alpha\in P_{\g}, 
\Re\{\alpha\}=0}
\mathfrak{g}^{[\alpha]}\otimes t^{\alpha+1}[t]\right),\\
\hat{\mathfrak{g}}^{[g]}_{-}&=\bigoplus_{\alpha\in P_{\g}}\mathfrak{g}^{[\alpha]}
\otimes t^{\alpha-1}\C[t^{-1}],\\
\hat{\mathfrak{g}}^{[g]}_{\I}&=\left(\bigoplus_{\alpha\in P_{\g}, 
\Re\{\alpha\}=0}\mathfrak{g}^{[\alpha]}
\otimes \C t^{\alpha}\right),\\
\hat{\mathfrak{g}}^{[g]}_{0}&=\hat{\mathfrak{g}}^{[g]}_{\I}\oplus \C\mathbf{k}.
\end{align*}
Then $\hat{\mathfrak{g}}^{[g]}_{+}$, $\hat{\mathfrak{g}}^{[g]}_{-}$, $\hat{\mathfrak{g}}^{[g]}_{\I}$ and 
$\hat{\mathfrak{g}}^{[g]}_{0}$ are subalgebras of $\hat{\mathfrak{g}}^{[g]}$ and 
$\hat{\mathfrak{g}}^{[g]}_{\I}$ is a subalgebra of $\hat{\mathfrak{g}}^{[g]}_{0}$. 
Moreover, we have a triangular decomposition 
$$\hat{\mathfrak{g}}^{[g]} = \hat{\mathfrak{g}}^{[g]}_{-} 
\oplus  \hat{\mathfrak{g}}^{[g]}_{0}  \oplus \hat{\mathfrak{g}}^{[g]}_{+}.$$

In this paper, we are interested in only those $\hat{\mathfrak{g}}^{[g]}$-modules with
 lower-bounded $\C$-gradings compatible with the grading of $\hat{\mathfrak{g}}^{[g]}$ and 
with actions of $g$. To be precise, we give the following definition:

\begin{defn}
{\rm A {\it graded $\hat{\mathfrak{g}}^{[g]}$-module} is a  $\hat{\mathfrak{g}}^{[g]}$-module 
$W$ with a $\C$-grading $W=\coprod_{n\in \C}W_{[n]}$ such that $\hat{\g}_{[m]}^{[g]}W_{[n]}
\subset W_{[m+n]}$ for 
$m\in P_{\g}+\Z$
 and $n\in \C$, where $\hat{\g}_{[m]}^{[g]}=\mathfrak{g}^{[\alpha]}
\otimes t^{m}$ for $m\in  (P_{\g}+\Z)\setminus \{0\}$
and $\hat{\g}_{[0]}^{[g]}=\g\otimes \C t^{0}\oplus \C\mathbf{k}$. 
A {\it graded $\hat{\mathfrak{g}}^{[g]}$-module of level $\ell$} is a graded $\hat{\mathfrak{g}}^{[g]}$-module
such that $\mathbf{k}$ acts as $\ell\in \C$. 
A {\it lower-bounded $\hat{\mathfrak{g}}^{[g]}$-module}  is a graded $\hat{\mathfrak{g}}^{[g]}$-module
$W=\coprod_{n\in \C}W_{[n]}$ such that  $W_{[n]}=0$ when $\Re(n)$ is sufficiently negative. A {\it grading-restricted 
$\hat{\mathfrak{g}}^{[g]}$-module} is 
a lower-bounded $\hat{\mathfrak{g}}^{[g]}$-module $W=\coprod_{n\in \C}W_{[n]}$
such that $\dim W_{[n]}<\infty$ for $n\in \C$. A {\it $\hat{\mathfrak{g}}^{[g]}$-module with a 
compatible action of $g$} or simply a {\it $\hat{\mathfrak{g}}^{[g]}$-module with an 
 action of $g$}  is a $\hat{\mathfrak{g}}^{[g]}$-module $W$ with actions of $g$, $\mathcal{S}_{g}$
and $\mathcal{N}_{g}$ satisfying the following conditions: (i) $W$ is a direct sum of 
generalized eigenspaces of $g$. (ii) 
$g=e^{2\pi i\mathcal{L}_{g}}$, where $\mathcal{L}_{g}$ is the operator on $W$ such that 
$\mathcal{S}_{g}$ and $\mathcal{N}_{g}$ on $W$ are the semisimple and nilpotent parts of $\mathcal{L}_{g}$.
(iii) $g(uw)=g(u)g(w)$ for $u\in \hat{\mathfrak{g}}^{[g]}$ and $w\in W$.}
\end{defn}

In this paper, $\hat{\mathfrak{g}}^{[g]}$-modules are always assumed to be graded and with 
compatible $g$ actions. So we shall call them simply $\hat{\mathfrak{g}}^{[g]}$-modules.
In particular, in this paper,  lower-bounded $\hat{\mathfrak{g}}^{[g]}$-modules 
and grading-restricted $\hat{\mathfrak{g}}^{[g]}$-modules
are always with compatible $g$ actions.

\renewcommand{\theequation}{\thesection.\arabic{equation}}
\renewcommand{\thethm}{\thesection.\arabic{thm}}
\setcounter{equation}{0}
\setcounter{thm}{0}
\section{Vertex operator algebras associated to affine Lie algebras and 
their automorphisms}

We recall the vertex operator algebras constructed from suitable modules 
for the affine Lie algebra $\hat{\g}$ and their automorphsims in this section.

Let $M$ be a $\mathfrak{g}$-module and let $\ell \in \C$.
Let $\hat{\mathfrak{g}}_{+}$ act on $M$ trivially and let
$\mathbf{k}$ act as the scalar multiplication by $\ell$. 
Then $M$ becomes a $\mathfrak{g} \oplus \C{\bf k} \oplus \hat{\mathfrak{g}}_{+}$-module
and we have an induced $\hat{\mathfrak{g}}$-module
$$\widehat{M}_{\ell} = U(\hat{\mathfrak{g}})
\otimes_{U(\mathfrak{g}\oplus \C\mathbf{k}\oplus \hat{\mathfrak{g}}_{+})}M,$$

Let $M=\C$ and let $\g$ act on $\C$ trivially. The corresponding $\hat{\mathfrak{g}}$-module
$\widehat{\C}_{\ell}$ is denoted by $M(\ell, 0)$. Let $J(\ell, 0)$ be 
the maximal proper submodule of $M(\ell, 0)$ 
and $L(\ell, 0) = M(\ell, 0)/J(\ell, 0)$. 
Then $L(\ell, 0)$ is the unique irreducible graded $\hat{\mathfrak{g}}$-module 
such that $\mathbf{k}$ acts as $\ell$ and the space of all elements annihilated by 
$\hat{\mathfrak{g}}_{+}$ is isomorphic to the trivial $\mathfrak{g}$-module $\C$.

Frenkel and Zhu  \cite{FZ} gave both $M(\ell, 0)$ and $L(\ell, 0)$ natural structures 
of  vertex operator algebras (see also \cite{LL}). In particular, $M(\ell, 0)$ and $L(\ell, 0)$
are grading-restricted vertex  algebras. 
We shall apply the results in 
\cite{H-const-twisted-mod} to construct lower-bounded and grading-restricted generalized twisted 
$M(\ell, 0)$- and $L(\ell, 0)$-modules. Since 
\cite{H-const-twisted-mod} needs the first construction of grading-restricted 
vertex algebras in \cite{H-2-const}, we describe
the grading-restricted vertex algebra structures on $M(\ell, 0)$ and $L(\ell, 0)$
using the construction in Section 3 in  \cite{H-2-const}. 

We discuss $M(\ell, 0)$ first. Note that $U(\hat{\mathfrak{g}}_{-})$ is 
linearly isomorphic to $M(\ell, 0)$. The $\Z_{+}$-grading on $\hat{\mathfrak{g}}_{-}$ induces an $\N$-grading 
on $M(\ell, \lambda)$. We denote the homogeneous subspace of $M(\ell, 0)$
of degree (conformal weight) $n$ by $M_{(n)}(\ell, 0)$ for $n\in \N$. We denote 
the action of $a\otimes t^{n}$ on $M(\ell, 0)$ by $a(n)$ for $a\in \g$ and $n\in \Z$. 
We also denote $1\in M(\ell, 0)$ by $\one_{M(\ell, 0)}$. Then $M(\ell, 0)$ is spanned 
by elements of the form $a_{1}(-n_{1})\cdots a_{k}(-n_{k})\one_{M(\ell, 0)}$ 
for $a_{1}, \dots, a_{k}\in \g$ 
and $n_{1}, \dots, n_{k}\in -\Z_{+}$. 
For $a\in \g$, let $a(x)=\sum_{n\in \Z}a(n)x^{-n-1}$. In particular, 
$z\mapsto a(z)$ for $z\in \C^{\times}$ is an analytic map from $\C^{\times}$ to 
$\hom(M(\ell, 0), \overline{M(\ell, 0)})$. 

Let $L_{M(\ell, 0)}(0)$ be the operator on $M(\ell, 0)$ giving the grading on $M(\ell, 0)$, that is,
$L_{M(\ell, 0)}(0)v=nv$ for $v\in (M_{(n)}(\ell, 0))$. 
We define an operator $L_{M(\ell, 0)}(-1)$ on $M(\ell, 0)$ by 
\begin{align*}
L&_{M(\ell, 0)}(-1)a_{1}(-n_{1})\cdots a_{k}(-n_{k})\one_{M(\ell, 0)}\nn
&=\sum_{i=1}^{k}n_{i}a_{1}(-n_{1})\cdots a_{i-1}(-n_{i-1})
a_{i}(-n_{i}-1)a_{i+1}(-n_{i+1})\cdots a_{k}(-n_{k})\one_{M(\ell, 0)}.
\end{align*}

It is easy to verify that the series $a(x)$ for $a\in \g$
and the operators $L_{M(\ell, 0)}(0)$ and $L_{M(\ell, 0)}(-1)$ have the following 
properties:

\begin{enumerate}

\item  For $a\in \g$, 
$[L_{M(\ell, 0)}(0), a(x)]=x\frac{d}{dx}a(x)+a(x)$.

\item $L_{M(\ell, 0)}(-1)\one=0$, 
$[L_{M(\ell, 0)}(-1), a(x)]=\frac{d}{dz}a(x)$
for $a\in \g$.

\item For $a\in \g$,  $a(x)\one_{M(\ell, 0)}\in M(\ell, 0)[[x]]$.
Moreover,  $\lim_{x\to 0}a(x)\one=a(-1)\one_{M(\ell, 0)}$.

\item The vector space $M(\ell, 0)$ is spanned by elements of the form 
$a_{1}(n_{1})\cdots a_{k}(n_{k})\one_{M(\ell, 0)}$ for $a_{1}, \dots, a_{k}\in \g$ 
and $n_{1}, \dots, n_{k}\in \Z$.

\item For $a, b \in \g$, 
$$(x_{1}-x_{2})^{2}a(x_{1})b(x_{2})=(x_{1}-x_{2})^{2}b(x_{2})a(x_{1}).$$

\end{enumerate}

Then by Proposition 3.3 in \cite{H-2-const},  
$\langle v', a_{1}(z_{1})\cdots a_{k}(z_{k})
v\rangle$ for $a_{1}, \dots, a_{k}\in \g$, $v\in M(\ell, 0)$ and $v'\in M(\ell, 0)'$
is absolutely convergent in the region $|z_{1}|>\cdots >|z_{k}|>0$ 
to a rational function, denoted by
$R(\langle v', a_{1}(z_{1})\cdots a_{k}(z_{k})
v\rangle)$,
in $z_{1}, \dots, z_{k}$ with the only possible 
poles at $z_{i}=0$ for $i=1, \dots k$ and $z_{i}=z_{j}$ for $i<j$, $i, j=1, \dots, k$. 

By Theorem 3.5 in  \cite{H-2-const}, 
the vector space $M(\ell, 0)$ equipped with the vertex operator map 
$$Y_{M(\ell, 0)}: M(\ell, 0)\otimes M(\ell, 0)\to M(\ell, 0)[[x, x^{-1}]]$$
defined by 
\begin{align*}
\langle v', &Y_{M(\ell, 0)}(a_{1}(n_{1})\cdots a_{k}(n_{k})\one_{M(\ell, 0)}, z)
v\rangle\nn
&=\res_{\xi_{1}=0}\cdots\res_{\xi_{k}=0}
\xi_{1}^{n_{1}}\cdots\xi_{k}^{n_{k}}
R(\langle v', a_{1}(\xi_{1}+z)\cdots a_{k}(\xi_{k}+z)
v\rangle)
\end{align*}
for $z\in \C^{\times}$, $a_{1}, \dots, a_{k}\in \g$,  $n_{1}, \dots, n_{k}\in \Z$, $v\in M(\ell, 0)$ 
and $v'\in M(\ell, 0)'$
and the vacuum 
$\one_{M(\ell, 0)}$ is a grading-restricted vertex algebra.
Moreover, this is the unique grading-restricted vertex algebra structure on 
$M(\ell, 0)$ with the vacuum $\one$
such that $Y(a(-1)\one, x)=a(x)$
for $a\in \g$. In particular, this grading-restricted vertex algebra structure on $M(\ell, 0)$
is the same as the one  constructed  in \cite{FZ} (see also \cite{LL}), that is, the graded space 
$M(\ell, 0)$, 
the vertex operator map $Y_{M(\ell, 0)}$, the operator $L_{M(\ell, 0)}(-1)$ and the vacuum one 
are equal to the graded space, the vertex operator maps, the operator $L(-1)$ and the vacuum in \cite{FZ}. 

Since $J(\ell, 0)$ is a $\hat{\g}$-module, we can define the action of $a(x)$ 
for $a\in \g$ on $L(\ell, 0)=M(\ell, 0)/J(\ell, 0)$. Similarly $L_{M(\ell, 0)}(0)$
and $L_{M(\ell, 0)}(-1)$ induce operators $L_{L(\ell, 0)}(0)$ and $L_{L(\ell, 0)}(-1)$.
We also have an element $\one_{L(\ell, 0)}=
\one+J(\ell, 0)\in L(\ell, 0)$. It is clear that the space $L(\ell, 0)$,
these series,  operators and the
element also satisfy the five properties above. Thus by 
Theorem 3.5 in  \cite{H-2-const}, $L(\ell, 0)$ equipped with the 
vertex operator map $Y_{L(\ell, 0)}$
defined by 
\begin{align*}
\langle v', &Y_{L(\ell, 0)}(a_{1}(n_{1})\cdots a_{k}(n_{k})\one_{L(\ell, 0)}, z)
v\rangle\nn
&=\res_{\xi_{1}=0}\cdots\res_{\xi_{k}=0}
\xi_{1}^{n_{1}}\cdots\xi_{k}^{n_{k}}
R(\langle v', a_{1}(\xi_{1}+z)\cdots a_{k}(\xi_{k}+z)
v\rangle)
\end{align*}
for $a_{1}, \dots, a_{k}\in \g$,  $n_{1}, \dots, n_{k}\in \Z$, $v\in L(\ell, 0)$ 
and $v'\in L(\ell, 0)'$
and the vacuum 
$\one_{L(\ell, 0)}$  is a grading-restricted vertex algebra.
Moreover, this is the unique grading-restricted vertex algebra structure on 
$L(\ell, 0)$ with the vacuum $\one$
such that $Y(a(-1)\one, x)=a(x)$
for $a\in \g$. In particular, this grading-restricted vertex algebra structure on $L(\ell, 0)$
is the same as the one  constructed first in \cite{FZ} (see also \cite{LL}). 

Let $g$ be an automorphism of $\g$ as is discussed in the preceding section. 
The actions of $g$, $\mathcal{L}_{g}$, $\mathcal{S}_{g}$ and
 $\mathcal{N}_{g}$ on
$\hat{\mathfrak{g}}$ further induce their actions, still denoted by $g$, $\mathcal{L}_{g}$, 
$\mathcal{S}_{g}$ and  $\mathcal{N}_{g}$, on $M(\ell, 0)$ and $L(\ell, 0)$. Moreover,
$g$, $e^{2\pi i\mathcal{S}_{g}}$ and $e^{2\pi i\mathcal{N}_{g}}$ are
all automorphisms of the grading-restricted vertex algebras $M(\ell, 0)$ and $L(\ell, 0)$. 

In the case that $\g$ is simple, we shall always take $(\cdot, \cdot)$ be the normalized 
Killing form such that $(\alpha, \alpha)=2$ for a long root $\alpha$.
Let $h^{\vee}$ be dual Coxeter number of $\mathfrak{g}$.
In the case that $\ell+h^{\vee}\ne 0$, the grading-restricted vertex algebra $M(\ell, 0)$ has a conformal element 
$$\omega_{M(\ell, 0)}=\frac{1}{2(\ell+h^{\vee})}\sum_{i=1}^{\dim \g}(a^{i})'(-1)a^{i}(-1)\one,$$
where $\{a^{i}\}_{i=1}^{\dim \g}$ is a basis of $\g$ and $\{(a^{i})'\}_{i=1}^{\dim \g}$ is 
its dual basis with respect to $(\cdot, \cdot)$. See \cite{FZ} 
and also \cite{LL}. The grading-restricted vertex algebra $M(\ell, 0)$ with this conformal
element is a vertex operator algebra (or a grading-restricted conformal vertex algebra). 
Moreover, $L_{M(\ell, 0)}(0)$ and $L_{M(\ell, 0)}(-1)$ above are in fact the 
coefficients of $x^{-2}$ and $x^{-1}$ in $Y_{M(\ell, 0)}(\omega_{M(\ell, 0)}, x)$.
Since $\omega_{M(\ell, 0)}$ is not in $J(\ell, 0)$, we see that $L(\ell, 0)$
has a conformal element $\omega_{L(\ell, 0)}=\omega_{M(\ell, 0)}+J(\ell, 0)$. 
So the grading-restricted vertex algebra $L(\ell, 0)$ with this conformal element 
is also a vertex operator algebra and $L_{L(\ell, 0)}(0)$ and $L_{L(\ell, 0)}(-1)$ are in fact the 
coefficients of $x^{-2}$ and $x^{-1}$ in $Y_{L(\ell, 0)}(\omega_{L(\ell, 0)}, x)$.
Again see \cite{FZ} and also \cite{LL} for details. 

Since the Killing form on $\g$ is invariant under
the action of $g$, the conformal element $\omega_{M(\ell, 0)}$
and $\omega_{L(\ell, 0)}$ are also invariant under $g$. 
Thus $g$, $e^{2\pi i\mathcal{S}_{g}}$ and $e^{2\pi i\mathcal{N}_{g}}$ are
 in fact  automorphisms of the vertex operator algebras
$M(\ell, 0)$ and $L(\ell, 0)$. Since $\mathcal{N}_{g}$ is a nilpotent operator on $\g$,
we must have $\mathcal{N}_{g}^{\dim \g}=0$ on $M(\ell, 0)$ and $L(\ell, 0)$. 

The decomposition 
$$\hat{\mathfrak{g}}=\coprod_{\alpha\in P_{\g}}
\hat{\mathfrak{g}}^{[\alpha]}$$
induced from the decomposition
$$\mathfrak{g}=\coprod_{\alpha\in P_{\g}}
\mathfrak{g}^{[\alpha]}$$
further induces decompositions
\begin{align*}
M(\ell, 0)&=\coprod_{\alpha\in P_{\g}}M^{[\alpha]}(\ell, 0),\\
L(\ell, 0)&=\coprod_{\alpha\in P_{\g}}L^{[\alpha]}(\ell, 0),
\end{align*}
where $M^{[\alpha]}(\ell, 0)$ and $L^{[\alpha]}(\ell, 0)$ for $\alpha\in P_{\g}$ 
are the generalized eigenspaces of $g$  (or the eigenspaces of 
$e^{2\pi i\mathcal{S}_{g}}$)  on $M(\ell, 0)$ and 
$L(\ell, 0)$, respectively, with the eigenvalue $e^{2\pi i\alpha}$. 
%These $\hat{\g}$-module $M(\ell, 0)$ also has a filtration given by the 
%operator $\mathcal{N}_{g}$ as follows: Let $M(\ell, 0)_{p}$ for $p\in \Z_{+}$ 
%be the subspace of $M(\ell, 0)$ spanned by elements of the form 
%$a_{1}(-n_{1})\cdots a_{k}(-n_{k})\one_{M(\ell, 0)}$ satisfying the following property:
%There exist $p_{1}, \dots, p_{k}\in \Z_{+}$ such that $p=p_{1}+\cdots +p_{k}$
%and $\mathcal{N}_{g}^{p_{1}}a_{1}=\cdots =\mathcal{N}_{g}^{p_{k}}a_{k}=0$. 
%Similarly we have a filtration for $L(\ell, 0)$. 
%We shall need the filtration 
%for $M(\ell, 0)$ later.  CHECK

We now choose suitable generating fields $a(x)$ such 
 that Assumption 2.1 in \cite{H-const-twisted-mod} is satisfied for
the grading-restricted vertex algebra $M(\ell, 0)$. 
Since $\mathcal{L}_{g}$ is an operator on a finite-dimensional vector space $\g$, 
we can find a Jordan basis $\{a^{i}\}_{i=1}^{\dim \g}$ for $g$, that is,
a basis $\{a^{i}\}_{i=1}^{\dim \g}$ of $\g$ such that under this basis, the matrix
representation of $\mathcal{L}_{g}$ is a Jordan canonical form. We use 
$I$ to denote the set $\{1, \dots, \dim \g\}$. Then the Jordan basis can be written 
as $\{a^{i}\}_{i\in I}$.
Since $\{a^{i}\}_{i\in I}$ is a basis of $\g$ and $M(\ell, 0)$
as a grading-restricted vertex algebra is generated by fields of the form 
$a(x)$ for $a\in \g$, $M(\ell, 0)$ is also generated by 
the fields $a^{i}(x)$ for $i\in I$. Since for $i\in I$,
$a^{i}$ is an element of a Jordan  basis, 
there exist an $\alpha_{i}\in P_{\g}$ and $n_{i}\in \Z$ such that 
$a^{i}$ is a generalized eigenvector of $\mathcal{L}_{g}$ with the eigenvalue
$\alpha_{i}+n_{i}$, or equivalently, a generalized  eigenvector of $g$ with the eigenvalue
$e^{2\pi i\alpha_{i}}$. Thus $a^{i}(-1)\one$ is also a generalized eigenvector of 
$\mathcal{L}_{g}$ on $M(\ell, 0)$ with the eigenvalue $\alpha_{i}+n_{i}$, or equivalently, a 
generalized eigenvector of $g$ on $M(\ell, 0)$ with the eigenvalue $e^{2\pi i\alpha_{i}}$.
Also since $a^{i}$ is an element of a Jordan  basis, 
either $\mathcal{N}_{g}a^{i}=0$ or $\mathcal{N}_{g}a^{i}$ is another element 
in the basis $\{a^{i}\}_{i\in I}$. Therefore there exists $\mathcal{N}_{g}(i)\in I$ 
such that $\mathcal{N}_{g}a^{i}=a^{\mathcal{N}_{g}(i)}$. Thus we also have 
$\mathcal{N}_{g}a^{i}(-1)\one=a^{\mathcal{N}_{g}(i)}(-1)\one$. 
These also hold for $L(\ell, 0)$.
In summary, these discussions give the following:

\begin{prop}\label{assum-va}
Assumption 
2.1 in \cite{H-const-twisted-mod} is satisfied by $M(\ell, 0)$ and $L(\ell, 0)$ with the 
set of generating fields 
$a^{i}(x)$ for $i\in I$.
\end{prop}

\renewcommand{\theequation}{\thesection.\arabic{equation}}
\renewcommand{\thethm}{\thesection.\arabic{thm}}
\setcounter{equation}{0}
\setcounter{thm}{0}
\section{Lower-bounded and grading-restricted generalized \\
twisted $M(\ell, 0)$-modules}

In this section, we first construct universal lower-bounded and grading-restricted 
generalized twisted modules for $M(\ell, 0)$ viewed as a grading-restricted vertex algebra in 
Subsection \ref{4.1}. Then in the case 
that $\g$ is simple and $\ell+h^{\vee}\ne 0$, we construct 
universal lower-bounded and grading-restricted 
generalized twisted modules for
$M(\ell, 0)$ viewed as a vertex operator algebra in Subsection \ref{4.2}. 
We also discuss their basic properties such as their universal properties and so on in 
Subsection \ref{4.3}.  

\subsection{The constructions when $M(\ell, 0)$ is viewed as a grading-restricted vertex algebra}\label{4.1}

Before we give constructions of universal lower-bounded and grading-restricted 
generalized twisted $M(\ell, 0)$-modules, we first show that 
such twisted modules must be a module for the twisted affine Lie algebra $\hat{\g}^{[g]}$. 
Let $W$ be a lower-bounded generalized 
$g$-twisted $M(\ell, 0)$-module. For $\alpha\in P_{\g}$,
$a\in \g^{[\alpha]}$ and $n\in \alpha+\Z$, we write
$$Y_{W}^{g}(a(-1)\one, x)=\sum_{k=0}^{K}\sum_{n\in \alpha+\Z}
(a_{W})_{n, k}x^{-n-1}(\log x)^{k}.$$
From (2.10) in \cite{HY},
$$Y_{W}^{g}(a(-1)\one, x)=(Y_{W}^{g})_{0}(x^{\mathcal{N}_{g}}
a(-1)\one, x)=\sum_{k=0}^{K}\frac{1}{k!}(Y_{W}^{g})_{0}((\mathcal{N}_{g}^{k}
a)(-1)\one, x)(\log x)^{k},$$
where $(Y_{W}^{g})_{0}(x^{\mathcal{N}_{g}}
a(-1)\one, x)$ is the constant term when $Y_{W}^{g}(a(-1)\one, x)$ is viewed as 
a polynomial in $\log x$. So in our notation, 
\begin{align*}
(Y_{W}^{g})_{0}(x^{\mathcal{N}_{g}}
a(-1)\one, x)&=\sum_{n\in \alpha+\Z}(a_{W})_{n, 0}x^{-n-1},\\ 
Y_{W}^{g}(a(-1)\one, x)
&=\sum_{k=0}^{K}\sum_{n\in \alpha+\Z}
\frac{1}{k!}((\mathcal{N}_{g}^{k}a)_{W})_{n, 0}x^{-n-1}(\log x)^{k}.
\end{align*}
We shall need the following version ((3.24)  in \cite{H-twist-vo}) 
of the Jacobi identity  for $(Y_{W}^{g})_{0}$
(obtained in \cite{B} and proved to be equivalent to the duality property for $(Y_{W}^{g})$  in \cite{HY}):
\begin{align}\label{jacobi}
x_0^{-1}&\delta\left(\frac{x_1 - x_2}{x_0}\right)
Y_{W}^{g}(u, x_1)
Y_{W}^{g}(v, x_2)-  x_0^{-1}\delta\left(\frac{- x_2 + x_1}{x_0}\right)
Y_{W}^{g}(v, x_2)
Y_{W}^{g}(u, x_1)\nn
&= x_1^{-1}\delta\left(\frac{x_2+x_0}{x_1}\right)
Y_{W}^{g}\left(Y_{V}\left(\left(\frac{x_2+x_0}{x_1}\right)^{\mathcal{L}_{g}}
u, x_0\right)v, x_2\right).
\end{align}
This Jacobi identity holds for lower-bounded generalized twisted modules and even more 
general types of twisted modules for an arbitrary grading-restricted vertex algebra,
including, in particular, $M(\ell, 0)$ or  $L(\ell, 0)$. 
Using the commutator formula obtained from this Jacobi identity, 
we have the following result of Bakalov in \cite{B} and, 
for reader's convenience, we give a proof:

\begin{prop}\label{twisted-affine-mod}
Let $W$ be a lower-bounded generalized $g$-twisted $M(\ell, 0)$-module.
Then $W$, with the action of $\hat{\g}^{[g]}$ given by $a\otimes t^{n}\mapsto 
(a_{W})_{n, 0}$ and $\mathbf{k}\mapsto \ell 1_{W}$ for $\alpha\in P_{\g}$,
$a\in \g^{[\alpha]}$ and $n\in \alpha+\Z$ and 
with the existing action of $g$, $\mathcal{S}_{g}$ and $\mathcal{N}_{g}$ on $W$,
 is a lower-bounded $\hat{\g}^{[g]}$-module of level $\ell$.
\end{prop}
\pf
Taking $\res_{x_{0}}$ on both sides of 
(\ref{jacobi}) and taking $u=a^{i}(-1)\one$ and 
$v=a^{j}(-1)\one$, we obtain the commutator formula
\begin{align}\label{commutator}
Y&_{W}^{g}(a^{i}(-1)\one, x_1)
Y_{W}^{g}(a^{j}(-1)\one, x_2)
 - Y_{W}^{g}(a^{j}(-1)\one, x_2)Y_{W}^{g}
(a^{i}(-1)\one, x_1)\nn
&= \res_{x_{0}}x_1^{-1}\delta\left(\frac{x_2+x_0}{x_1}\right)
Y_{W}^{g}\left(Y_{M(\ell, 0)}\left(\left(\frac{x_2+x_0}{x_1}\right)^{\mathcal{L}_{g}}a^{i}(-1)\one, x_0\right)a^{j}(-1)\one, x_2\right).
\end{align}
By definition,  
the left-hand side of (\ref{commutator}) is equal to 
\begin{align}\label{commutator-l}
&\sum_{m\in \alpha_{i}+\Z}\sum_{k\in \N}\sum_{n\in \alpha_{j}+\Z}
\sum_{l\in \N}\frac{(-1)^{k}}{k!}\frac{(-1)^{l}}{l!}
((\mathcal{N}_{g}^{k}a^{i})_{W})_{m, 0}
((\mathcal{N}_{g}^{l}a^{j})_{W})_{n, 0}x_{1}^{-m-1}
x_{1}^{-n-1}(\log x_{1})^{k}(\log x_{2})^{l}\nn
& -\sum_{m\in \alpha_{i}+\Z}\sum_{k\in \N}\sum_{n\in \alpha_{j}+\Z}
\sum_{l\in \N}\frac{(-1)^{k}}{k!}\frac{(-1)^{l}}{l!}
((\mathcal{N}_{g}^{l}
a^{j})_{W})_{n, 0}((\mathcal{N}_{g}^{k}a^{i})_{W})_{m, 0}x_{1}^{-m-1}
x_{1}^{-n-1}(\log x_{1})^{k}(\log x_{2})^{l}.
\end{align}
On the other hand,  by straightforward calculations, we see that  the right-hand side of (\ref{commutator}) is equal to 
\begin{align}\label{commutator-r}
& \res_{x_{0}}e^{x_{0}\frac{\partial}{\partial x_{2}}}
x_1^{-1}\delta\left(\frac{x_2}{x_1}\right)
Y_{W}^{g}\left(Y_{M(\ell, 0)}
\left(\left(\frac{x_2}{x_1}\right)^{\mathcal{L}_{g}}a^{i}(-1)\one, 
x_0\right)a^{j}(-1)\one, x_2\right)\nn
&=\sum_{n\in \Z}\sum_{k\in \N}\sum_{l\in \N}
\sum_{p\in \alpha_{i}+\alpha_{j}+\Z}
\frac{(-1)^{k}}{k!}\frac{(-1)^{l}}{l!}
(([\mathcal{N}_{g}^{k}a^{i}, \mathcal{N}_{g}^{l}a^{j}])_{W})_{p, 0}%\cdot\nn
%&\quad\quad\quad\cdot 
x_{1}^{-n-\alpha_{i}-1}x_{2}^{n+\alpha_{i}-p-1}
(\log x_{1})^{k}(\log x_{2})^{l}
\nn
&\quad +\sum_{n\in \Z}\sum_{k\in \N}\sum_{l\in \N}
\frac{(-1)^{k}}{k!}\frac{(-1)^{l}}{l!}
\ell(\mathcal{N}_{g}^{k}a^{i}, \mathcal{N}_{g}^{l}a^{j})%\cdot\nn
%&\quad\quad\quad\cdot 
\frac{\partial}{\partial x_{2}}
x_{1}^{-n-\alpha_{i}-1}x_{2}^{n+\alpha_{i}}
(\log x_{1})^{k}(\log x_{2})^{l}.
\end{align}
Taking coefficients of $x_{1}^{-m-1}x_{2}^{-n-1}(\log x_{1})^{k}(\log x_{2})^{l}$
for $m\in \alpha_{i}+\Z$, $n\in \alpha_{j}+\Z$, $k, l\in \N$
in both sides of (\ref{commutator}),
using (\ref{commutator-l}) and (\ref{commutator-r}), dividing
the both  results by $\frac{(-1)^{k}}{k!}\frac{(-1)^{l}}{l!}$ and then using 
Corollary \ref{N-g-inv},  we obtain
\begin{align}\label{commutator-comp}
&[((\mathcal{N}_{g}^{k}a^{i})_{W})_{m, 0}, 
((\mathcal{N}_{g}^{l}
a^{j})_{W})_{n, 0}]\nn
&\quad =(([\mathcal{N}_{g}^{k}a^{i}, \mathcal{N}_{g}^{l}a^{j}])_{W})_{m+n, 0}
+m (\mathcal{N}_{g}^{k}a^{i}, \mathcal{N}_{g}^{l}a^{j})\delta_{m+n, 0}\ell
+(\mathcal{N}_{g}^{k+1}a^{i}, \mathcal{N}_{g}^{l}a^{j})\delta_{m+n, 0}\ell.
\end{align}
Let $a=\mathcal{N}_{g}^{k}a^{i}$ and $b=\mathcal{N}_{g}^{l}
a^{j}$. Also note that $\g$ is certainly spanned by such $a$ and $b$. 
Then (\ref{commutator-comp}) is exactly what we can also obtain by replacing
$a\otimes t^{n}$ and $\mathbf{k}$ in (\ref{bracket-1})  by 
$(a_{W})_{n, 0}$ and $\ell 1_{W}$ for
$a\in \g^{[\alpha]}$ and $n\in \alpha+\Z$. 
Thus (\ref{commutator-comp}) gives $W$ a 
structure of a lower-bounded $\hat{\g}^{[g]}$-module  of level $\ell$.
\epfv

We first construct and identify explicitly universal lower-bounded generalized 
$g$-twisted $M(\ell, 0)$-modules 
generated by a space annihilated by $\hat{\g}^{[g]}_{+}$
using the results in Section 5 of \cite{H-const-twisted-mod}
when $M(\ell, 0)$ is viewed as a grading-restricted vertex algebra. 

Let $M$ be a vector space. Assume that $g$ acts on $M$ and there is an operator 
$L_{M}(0)$ on $M$. If $M$ is finite dimensional, then there exist operators $\mathcal{L}_{g}$,
$\mathcal{S}_{g}$, $\mathcal{N}_{g}$ such that on $M$, $g=e^{2\pi i\mathcal{L}_{g}}$
and $\mathcal{S}_{g}$ and $\mathcal{N}_{g}$ are the semisimple and nilpotent, respectively,
parts of $\mathcal{L}_{g}$. In this case, $M$ is also a direct sum of generalized 
eigenspaces for the operator $L_{M}(0)$ and $L_{M}(0)$ can be decomposed as 
the sum of its semisimple part $L_{M}(0)_{S}$ and nilpotent part $L_{M}(0)_{N}$. 
Moreover, the real parts of the eigenvalues of 
$L_{M}(0)$ has a lower bound. In the case that $M$ is infinite dimensional, we 
assume that all of these properties for $g$ and $L_{M}(0)$ hold. 
We call the eigenvalue of a generalized eigenvector 
$w\in M$ for $L_{M}(0)$ the {\it (conformal) weight} of $w$ and denote 
it by $\wt w$. We first assume that $M$ is itself a generalized eigenspace of 
$L_{M}(0)$ with eigenvalue $h$.

Let $\{w^{a}\}_{a\in A}$ be a basis of $M$ consisting of vectors homogeneous in
$g$-weights (eigenvalues of $g$)  such that for $a\in A$, 
either $L_{M}(0)_{N}w^{a}=0$ or there exists $L_{M}(0)_{N}(a)\in A$ 
such that $L_{M}(0)_{N}w^{a}=w^{L_{M}(0)_{N}(a)}$. For simplicity, when 
$L_{M}(0)_{N}w^{a}=0$, we shall use $w^{L_{M}(0)_{N}(a)}$ to denote $0$. 
Then for $a\in A$, we always have $L_{M}(0)_{N}w^{a}=w^{L_{M}(0)_{N}(a)}$.
For $a\in A$, let $\alpha^{a}\in \C$ such that $\Re(\alpha^{a})
\in [0, 1)$ and $e^{2\pi i\alpha^{a}}$ is the eigenvalue of $g$ for the generalized 
eigenvector $w^{a}$.

Taking the grading-restricted vertex algebra $V$, the space $M$
and $B\in \R$  in Section 5 of \cite{H-const-twisted-mod}
to be $M(\ell, 0)$, the space $M$ above and $\Re(h)$, respectively, 
we obtain the universal lower-bounded  generalized $g$-twisted $M(\ell, 0)$-module
$\widehat{M}^{[g]}_{h}$, which we shall denote by $\widehat{M}^{[g]}_{\ell, h}$
to exhibit explicit the dependence on $\ell$. The twisted generating fields and generator twist fields 
for $\widehat{M}^{[g]}_{\ell, h}$ are denoted by
$$a^{i}_{\widehat{M}^{[g]}_{\ell, h}}(x)=\sum_{n\in \alpha^{i}+\Z}\sum_{k=0}^{K^{i}}
(a^{i}_{\widehat{M}^{[g]}_{\ell, h}})_{n, k}x^{-n-1}(\log x)^{k}$$
for $i\in I$ and 
$$\psi^{a}_{\widehat{M}^{[g]}_{\ell, h}}(x)=\sum_{n\in \alpha^{i}+\Z}\sum_{k=0}^{K^{i}}
(\psi^{a}_{\widehat{M}^{[g]}_{\ell, h}})_{n, k}x^{-n-1}(\log x)^{k}$$ 
for $a\in A$. 
For simplicity, we shall denote $a^{i}_{\widehat{M}^{[g]}_{\ell, h}}(x)$ and  
$(a^{i}_{\widehat{M}^{[g]}_{\ell, h}})_{n, k}$ by 
$a^{i}_{[g], \ell}(x)$, $(a^{i}_{[g], \ell})_{n, k}$, respectively, since their 
commutators involve $\ell$ and denote 
$\psi^{a}_{\widehat{M}^{[g]}_{\ell, h}}(x)$ and $(\psi^{a}_{\widehat{M}^{[g]}_{\ell, h}})_{n, k}$  by
$\psi^{a}_{[g]}(x)$ and $(\psi^{a}_{[g]})_{n, k}$, respectively. 
For a general element $a\in \g^{[\alpha]}$ and $w\in M^{[\beta]}$, we shall use the 
similar notations to denote the twisted and twist fields associated to $a$ and $w$, respectively,
and similarly for their components.

The construction above is based on the assumption that $M$ is itself a generalized eigenspace of 
$L_{M}(0)$ with eigenvalue $h$. In the general case, 
$M=\coprod_{h\in Q_{M}}M_{[h]}$, where $Q_{M}$ is the set of all eigenvalues of $L_{M}(0)$
and $M_{[h]}$ is the generalized eigenspace
of $L_{M}(0)$ with the eigenvalue $h$. In this case, 
we have the 
lower-bounded generalized $g$-twisted $M(\ell, 0)$-module
$\coprod_{h\in Q_{M}}\widehat{(M_{[h]})}_{\ell, h}^{[g]}$, which we shall denote 
by $\widehat{M}^{[g]}_{\ell}$. For $h\in Q_{M}$, we have a basis $\{w^{a}\}_{a\in A_{h}}$
of $M_{[h]}$ satisfying the condition $L_{M}(0)_{N}w^{a}=w^{L_{M}(0)_{N}(a)}$ for 
$a\in A_{h}$. Let $A=\sqcup_{h\in Q_{M}}A_{h}$. Then we have a basis 
$\{w^{a}\}_{a\in A}$ of $M$ satisfying the same condition for all $a\in A$. 

We now construct 
a  lower-bounded $\hat{\g}^{[g]}$-module that we will prove to be equivalent 
to $\widehat{M}^{[g]}_{\ell}$ viewed as a lower-bounded $\hat{\g}^{[g]}$-module. 
Let $L_{-1}$ be a basis of a one-dimensional vector space $\C L_{-1}$. Let $T(\C L_{-1})$ be the tensor algebra
of the one-dimensional space $\C L_{-1}$. Consider the vector space 
$\Lambda(M)=T(\C L_{-1})\otimes M$. We define actions of $g$, $\mathcal{L}_{g}$, $\mathcal{S}_{g}$ and 
$\mathcal{N}_{g}$ on $\Lambda(M)$ by acting only on the second tensor factor $M$. 
We define an operator $L_{\Lambda(M)}(0)$  on $\Lambda(M)$
by $L_{\Lambda(M)}(0)(L_{-1}^{m}\otimes w)=
m(L_{-1}^{m}\otimes w)+L_{-1}^{m}\otimes L_{M}(0)w$
for $m\in \N$ and $w\in M$.  We also define operators $L_{\Lambda(M)}(0)_{N}$ and 
$L_{\Lambda(M)}(0)_{S}$  on $\Lambda(M)$ by 
$L_{\Lambda(M)}(0)_{S}(L_{-1}^{m}\otimes w)=
m(L_{-1}^{m}\otimes w)+L_{-1}^{m}\otimes L_{M}(0)_{S}w$ and 
$L_{\Lambda(M)}(0)_{N}(L_{-1}^{m}\otimes w)=
L_{-1}^{m}\otimes L_{M}(0)_{N}w$, respectively, for $m\in \N$ and $w\in M$. 
Then $L_{\Lambda(M)}(0)_{N}$ and 
$L_{\Lambda(M)}(0)_{S}$  are the semisimple and nilpotent, respectively, parts of 
$L_{\Lambda(M)}(0)$. The space $\Lambda(M)$ is graded by the eigenvalues of $L_{\Lambda(M)}(0)$. 
We define another operator $L_{\Lambda(M)}(-1)$ on $\Lambda(M)$ by 
$L_{\Lambda(M)}(-1)(L_{-1}^{m}\otimes w)=L_{-1}^{m+1}\otimes w$. Then 
$\Lambda(M)$ is spanned by elements of the form $L_{\Lambda(M)}(-1)^{m}(1\otimes w)$ for $m\in \N$ and 
$w\in M$. For simplicity, we shall identify $1\otimes w$ with $w\in M$ and hence embed
$M$ as a subspace of $\Lambda(M)$. Thus $\Lambda(M)$ is spanned by elements of the form
$L_{\Lambda(M)}(-1)^{m}w$ for $w\in M$. 

Let $\hat{\g}^{[g]}_{+}$ act on $M$ as $0$. We define an action of $\hat{\g}^{[g]}_{+}$
on $\Lambda(M) $ by the commutator formula 
\begin{equation}\label{comm-a-m-L-1}
[a(m), L_{\Lambda(M)}(-1)]=ma(m-1)+(\mathcal{N}_{g}a)(m-1)
\end{equation}
 for $a\in \g^{[\alpha]}$ and $m\in \alpha+\N$ when $\Re(\alpha)>0$ and 
$m\in \alpha+\Z_{+}$ when $\Re(\alpha)=0$. Let $\mathbf{k}$ act on $\Lambda(M)$ 
as $\ell$. Then it is clear that $\Lambda(M)$ is a $U(\hat{\g}^{[g]}_{+}\oplus \C\mathbf{k})$-module.
and we have the induced lower-bounded $\hat{\g}^{[g]}$-module
$U(\hat{\g}^{[g]})\otimes _{U(\hat{\g}^{[g]}_{+}\oplus \C\mathbf{k})}\Lambda(M)$
(recalling that by a
lower-bounded $\hat{\g}^{[g]}$-module we mean  a
lower-bounded $\hat{\g}^{[g]}$-module with a compatible $g$ action). 
Using the commutator formula (\ref{comm-a-m-L-1}), we can extend the operator 
$L_{\Lambda(M)}(-1)$ to an operator on 
$U(\hat{\g}^{[g]})\otimes _{U(\hat{\g}^{[g]}_{+}\oplus \C\mathbf{k})}\Lambda(M)$.
For simplicity, we shall still denote this extension of $L_{\Lambda(M)}(-1)$ by
the same notation $L_{\Lambda(M)}(-1)$. But note that $L_{\Lambda(M)}(-1)$
acts on $U(\hat{\g}^{[g]})\otimes _{U(\hat{\g}^{[g]}_{+}\oplus \C\mathbf{k})}\Lambda(M)$ now. 

\begin{thm}\label{ulbtm-va}
As a lower-bounded $\hat{\g}^{[g]}$-module,
$\widehat{M}_{\ell}^{[g]}$ 
is equivalent    to $U(\hat{\g}^{[g]})
\otimes _{U(\hat{\g}^{[g]}_{+}\oplus \C\mathbf{k})}\Lambda(M)$.
\end{thm}
\pf
Consider the subspace $\widehat{\Lambda}(M)$ of $\widehat{M}_{\ell}^{[g]}$
spanned by elements of the form 
\begin{equation}\label{ulbtm-va-1}
L_{\widehat{M}^{[g]}_{\ell}}(-1)^{k}(\psi_{[g]}^{b})_{-1, 0}\one
\end{equation}
 for $k\in \N$ and $b\in A$. Then we have 
a linear map $\rho: \Lambda(M)\to \widehat{\Lambda}(M)$ defined by 
$$\rho(L_{\Lambda(M)}(-1)^{k}w^{b})=L_{\widehat{M}^{[g]}_{\ell}}(-1)^{k}
(\psi_{[g]}^{b})_{-1, 0}\one$$
for $k\in \N$ and $b\in A$. In particular, 
$\rho(w^{b})=(\psi_{[g]}^{b})_{-1, 0}\one$ for $b\in A$. So
$\rho(M)$ is the subspace of $\widehat{M}_{\ell}^{[g]}$ spanned by 
$(\psi_{[g]}^{b})_{-1, 0}\one$ for $b\in A$. From the 
$\hat{\g}^{[g]}$-module structure on $\widehat{M}_{\ell}^{[g]}$, we see that 
$\hat{\g}^{[g]}_{+}$ acts on $\rho(M)$ as $0$. From the 
commutator formula 
$$[L_{\widehat{M}^{[g]}_{\ell, h}}(-1), a_{[g], \ell}(x)]=\frac{d}{dx}a_{[g], \ell}(x),$$
we obtain 
$$[L_{\widehat{M}^{[g]}_{\ell, h}}(-1), (a_{[g], \ell})_{m, 0}]=m(a_{[g], \ell})_{m-1, 0}
+((\mathcal{N}_{g}a)_{[g], \ell})_{m-1, 0}$$
for $a\in \g^{[\alpha]}$ and $m\in \alpha+\N$ when $\Re(\alpha)>0$ and 
$m\in \alpha+\Z_{+}$ when $\Re(\alpha)=0$. Thus we also have an action of 
$\hat{\g}^{[g]}_{+}$ on $\widehat{\Lambda}(M)$. From the 
$\hat{\g}^{[g]}$-module structure on $\widehat{M}_{\ell}^{[g]}$ again, we see that 
$\mathbf{k}$ acts on $\widehat{M}_{\ell}^{[g]}$ as $\ell$. These actions give $\widehat{\Lambda}(M)$
a $\hat{\g}^{[g]}_{+}\oplus \C\mathbf{k}$-module structure. From the definitions of $\rho$ and the 
$\hat{\g}^{[g]}_{+}\oplus \C\mathbf{k}$-module structures on $\Lambda(M)$ and 
$\widehat{\Lambda}(M)$, we see that $\rho$ is in fact a $\hat{\g}^{[g]}_{+}\oplus \C\mathbf{k}$-module
map. Moreover, by Theorem 2.4 in \cite{H-exist-twisted-mod}, for $h\in Q_{M}$,
$L_{\widehat{M}^{[g]}_{\ell, h}}(-1)^{k}
(\psi_{[g]}^{b})_{-1, 0}\one$ for $k\in \N$ and $b\in A_{h}$ are linearly independent and thus form a basis of 
$\widehat{\Lambda}(M_{h})$, which is the subspace of $\widehat{\Lambda}(M)$
spanned by elements of the form (\ref{ulbtm-va-1}) for $k\in \N$ and $b\in A_{h}$.
Then $L_{\widehat{M}^{[g]}_{\ell, h}}(-1)^{k}
(\psi_{[g]}^{b})_{-1, 0}\one$ for $k\in \N$ and $b\in A$ form a basis of $\widehat{\Lambda}(M)$.
So $\rho$ is in fact an equivalence of $\hat{\g}^{[g]}_{+}\oplus \C\mathbf{k}$-modules and commutes
with the actions of $L_{\widehat{M}^{[g]}_{\ell}}(-1)$ and $L_{\widehat{M}^{[g]}_{\ell}}(-1)$.

Now by the universal property of the induced module
$U(\hat{\g}^{[g]})
\otimes _{U(\hat{\g}^{[g]}_{+}\oplus \C\mathbf{k})}\Lambda(M)$, there exists
a unique $\hat{\g}^{[g]}$-module map 
$$\hat{\rho}: U(\hat{\g}^{[g]})\otimes _{U(\hat{\g}^{[g]}_{+}\oplus \C\mathbf{k})}\Lambda(M)\to 
\widehat{M}^{[g]}_{\ell}$$
such that $\hat{\rho}|_{\Lambda(M)}=\rho$. Since $\widehat{M}^{[g]}_{\ell}$ as a 
 $\hat{\g}^{[g]}$-module is generated by $\Lambda(M)$, 
$\hat{\rho}$ is surjective. We need only prove that $\hat{\rho}$ is injective. 

For $h\in Q_{M}$,  by Theorem 2.3 in \cite{H-exist-twisted-mod}, $\widehat{M}_{\ell, h}^{[g]}$ is spanned by 
elements of the form 
$$(a^{i_{1}}_{[g], \ell})_{n_{1}, k_{1}}
\cdots (a^{i_{l}}_{[g], \ell})_{n_{l}, k_{l}}L_{\widehat{M}^{[g]}_{\ell, h}}(-1)^{k}
(\psi_{[g]}^{b})_{-1, 0}\one$$
for $n_{j}\in \alpha^{i_{j}}+\Z$, $0\le k_{j}\le K_{j}$, $k\in \N$ and $b\in A_{h}$. 
On the other hand, from 
$$a^{i}_{[g], \ell}(x)=Y^{g}_{\widehat{M}^{[g]}_{\ell, h}}(a^{i}(-1)\one, x)
=(Y^{g}_{\widehat{M}^{[g]}_{\ell, h}})_{0}(x^{-\mathcal{N}_{g}}a^{i}(-1)\one, x),$$
we obtain 
$$(a^{i}_{[g], \ell})_{n_{i}, k_{i}}=\frac{(-1)^{k_{i}}}{k_{i}!}(\mathcal{N}_{g}^{k_{i}}
a^{i}_{[g], \ell})_{n_{i},0}.$$
Moreover, $\mathcal{N}_{g}^{k_{i}}
a^{i}_{[g], \ell}$ is a linear combination of $a^{j}$ for $j\in I$ since $a^{j}$ for $j\in I$
form a basis of $\g$. Thus $\widehat{M}_{\ell, h}^{[g]}$ is spanned by elements of the form
\begin{equation}\label{element-form-0}
(a^{i_{1}}_{[g], \ell})_{n_{1},0}
\cdots (a^{i_{l}}_{[g], \ell})_{n_{l}, 0}L_{\widehat{M}^{[g]}_{\ell, h}}(-1)^{k}
(\psi_{[g]}^{b})_{-1, 0}\one
\end{equation}
for $i_{j}\in I$,
$n_{j}\in \alpha^{i_{j}}+\Z$ for $j=1, \dots, l$, $k\in \N$ and
$b\in A_{h}$. Therefore $\widehat{M}^{g]}_{\ell}$ is spanned by elements of the form 
(\ref{element-form-0}) with $L_{\widehat{M}^{[g]}_{\ell, h}}(-1)$ replaced by 
$L_{\widehat{M}^{[g]}_{\ell}}(-1)$ for $i_{j}\in I$,
$n_{j}\in \alpha^{i_{j}}+\Z$ for $j=1, \dots, l$, $k\in \N$ and
$b\in A$. 

On the other hand, $U(\hat{\g}^{[g]})
\otimes _{U(\hat{\g}^{[g]}_{+}\oplus \C\mathbf{k})}\Lambda(M)$ is spanned by elements of the form
\begin{equation}\label{element-form-0-1}
a^{i_{1}}(n_{1})
\cdots a^{i_{l}}(n_{l})
L_{\Lambda(M)}(-1)^{k}w^{b}
\end{equation}
for $i_{j}\in I$,
$n_{j}\in \alpha^{i_{j}}+\Z$ for $j=1, \dots, l$, $k\in \N$, 
$b\in A$. Since $\hat{\rho}$ is a $\g^{[g]}$-module map, we have 
$$\hat{\rho}(a^{i_{1}}(n_{1})
\cdots a^{i_{l}}(n_{l})L_{\Lambda(M)}(-1)
w^{b})=(a^{i_{1}}_{[g], \ell})_{n_{1},0}
\cdots (a^{i_{l}}_{[g], \ell})_{n_{l}, 0}L_{\widehat{M}^{[g]}_{\ell}}(-1)
(\psi_{[g]}^{b})_{-1, 0}\one$$
for $i_{j}\in I$,
$n_{j}\in \alpha^{i_{j}}+\Z$ for $j=1, \dots, l$, $k\in \N$ and
$b\in A$. 
To prove that $\hat{\rho}$ is injective, we prove that if we replace
$a^{i}(n)$, $L_{\Lambda(M)}(-1)$ and $w^{b}$  by 
$(a_{[g], \ell}^{i})_{n, 0}$, $L_{\widehat{M}^{[g]}_{\ell}}(-1)$ and 
$(\psi_{[g]}^{b})_{-1, 0}\one$, respectively, for $i\in I$, $n\in \alpha^{i}+\Z$
and $b\in A$, 
the relations satisfied by elements of the spanning sets
(\ref{element-form-0-1}) of $U(\hat{\g}^{[g]})
\otimes _{U(\hat{\g}^{[g]}_{+}\oplus \C\mathbf{k})}\Lambda(M)$
must be satisfied by 
elements of the spanning sets (\ref{element-form-0}) of $\widehat{M}_{\ell}^{[g]}$.

To prove this, we first list all the relations satisfied  by
 (\ref{element-form-0}). 
From the construction of $\widehat{M}^{[g]}_{\ell, h}$ for $h\in Q_{M}$ 
given in Section 5 of \cite{H-const-twisted-mod} and Theorem 2.4 in \cite{H-exist-twisted-mod},
we see that the only relations satisfied by elements of the form  (\ref{element-form-0})  
for $i_{j}\in I$,
$n_{j}\in \alpha^{i_{j}}+\Z$ for $j=1, \dots, l$ and
$b\in A$  are generated by the following:
(i) A homogeneous element of the form (\ref{element-form-0})  with 
$i_{j}\in I$,
$n_{j}\in \alpha^{i_{j}}+\Z$ for $j=1, \dots, l$, $k\in \N$ and
$b\in A_{h}$ satisfying $-n_{1}-\cdots -n_{l}<\Re(h)$ is equal to $0$.
(ii) The relations induced from the coefficients of the weak commutativity for the 
generating $g$-twisted fields $a_{[g], \ell}^{i}(x)$ for $i\in I$.
(iii) The commutator relations 
between $(a_{[g], \ell}^{i})_{n, 0}$ and $L_{\widehat{M}^{[g]}_{h}}(-1)$
for $a\in \g^{[\alpha]}$ and $m\in \alpha+\N$ (when $\Re(\alpha)>0$) and 
$m\in \alpha+\Z_{+}$ (when $\Re(\alpha)=0$).
The other relations given in Section 5 of \cite{H-const-twisted-mod} involve elements that are 
not of the form  (\ref{element-form-0}). 

We need only prove that elements of the form
 (\ref{element-form-0-1}) also satisfy  the relations corresponding to 
the relations (i), (ii) and (iii). By the definitions of the actions of $a^{i}(n)$
on $U(\hat{\g}^{[g]})
\otimes _{U(\hat{\g}^{[g]}_{+}\oplus \C\mathbf{k})}\Lambda(M)$ and the fact that
the weights of $w^{b}$ for $b\in A_{h}$ are $h$,
elements of the form  (\ref{element-form-0-1}) satisfy the relations corresponding to (i). 
Since $a^{i}(n)$ and $(a_{[g], \ell}^{i})_{n, 0}$ 
for $i\in I$  and $n\in\alpha^{i}+\Z$ satisfy the same commutator formula, 
$a^{i}(x)=\sum_{n\in \alpha^{i}+\Z}a^{i}(n)x^{-n-1}$ and 
$a_{[g], \ell}^{i}(x)$ for $i\in I$ also satisfy the same commutator formula.
Since weak commutativity follows from the commutator formula for generating twisted fields, 
we see that elements of the form   (\ref{element-form-0-1}) satisfy 
the relations corresponding to (ii). Since  $\rho$ is in fact an equivalence of 
$\hat{\g}^{[g]}_{+}\oplus \C\mathbf{k}$-modules and commutes
with the actions of $L_{\widehat{M}^{[g]}_{\ell}}(-1)$ and $L_{\widehat{M}^{[g]}_{\ell}}(-1)$, 
elements of the form  (\ref{element-form-0-1}) satisfy 
the relations corresponding to (iii). This finishes the proof.
\epfv

\begin{rema}\label{elemenet-order-1}
{\rm By the Poincar\'{e}-Birkhoff-Witt theorem, the induced lower-bounded $\hat{\g}^{[g]}$-module
$U(\hat{\g}^{[g]}) \otimes _{U(\hat{\g}^{[g]}_{+}\oplus \C\mathbf{k})}\Lambda(M)$
is linearly isomorphic to $U(\hat{\g}^{[g]}_{-})\otimes U(\hat{\g}^{[g]}_{\mathbb{I}}) 
\otimes\Lambda(M)$. In particular, $U(\hat{\g}^{[g]}) \otimes _{U(\hat{\g}^{[g]}_{+}
\oplus \C\mathbf{k})}\Lambda(M)$ is spanned by elements of the form
$$%\begin{equation}\label{element-form-0-2}
a^{i_{1}}(n_{1})
\cdots a^{i_{l}}(n_{l})a^{j_{1}}(\alpha^{j_{1}})
\cdots a^{j_{m}}(\alpha^{j_{m}})
L_{\Lambda(M)}(-1)^{k}w^{b}
$$%\end{equation}
for $i_{p}, j_{q}\in I$,
$n_{p}\in \alpha^{i_{p}}-\Z_{+}$, $\Re(\alpha^{j_{q}})=0$ for $p=1, \dots, l$, 
$q=1, \dots, m$, $k\in \N$ and
$b\in A$. Using the commutator formula between $a^{j}(\alpha^{j})$ and $L_{\Lambda(M)}(-1)$
for $j\in I_{\I}$, we see that $U(\hat{\g}^{[g]}) \otimes _{U(\hat{\g}^{[g]}_{+}
\oplus \C\mathbf{k})}\Lambda(M)$ is also spanned by elements of the form
\begin{equation}\label{element-form-0-2}
a^{i_{1}}(n_{1})
\cdots a^{i_{l}}(n_{l})
L_{\Lambda(M)}(-1)^{k}a^{j_{1}}(\alpha^{j_{1}})
\cdots a^{j_{m}}(\alpha^{j_{m}})w^{b}
\end{equation}
for $i_{p}, j_{q}\in I$,
$n_{p}\in \alpha^{i_{p}}-\Z_{+}$, $\Re(\alpha^{j_{q}})=0$ for $p=1, \dots, l$, 
$q=1, \dots, m$, $k\in \N$ and
$b\in A$. 
By Theorem \ref{ulbtm-va}, we see that $\widehat{M}_{\ell, h}^{[g]}$  
is spanned by elements of the form
$$%\begin{equation}\label{element-form-0-3}
(a^{i_{1}}_{[g], \ell})_{n_{1},0}
\cdots (a^{i_{l}}_{[g], \ell})_{n_{l}, 0}
(a^{j_{1}}_{[g], \ell})_{\alpha^{j_{1}},0}
\cdots (a^{j_{m}}_{[g], \ell})_{\alpha^{j_{m}}, 0}L_{\widehat{M}^{[g]}_{\ell, h}}(-1)^{k}
(\psi_{[g]}^{b})_{-1, 0}\one
$$%\end{equation}
for $i_{p}, j_{q}\in I$,
$n_{p}\in \alpha^{i_{p}}-\Z_{+}$, $\Re(\alpha^{j_{q}})=0$ for $p=1, \dots, l$, 
$q=1, \dots, m$, $k\in \N$ and
$b\in A$. Using the commutator formula between 
$(a^{j}_{[g], \ell})_{\alpha^{j},0}$ and $L_{\widehat{M}^{[g]}_{\ell, h}}(-1)$,
we see that $\widehat{M}_{\ell, h}^{[g]}$  
is also spanned by elements of the form
\begin{equation}\label{element-form-0-3}
(a^{i_{1}}_{[g], \ell})_{n_{1},0}
\cdots (a^{i_{l}}_{[g], \ell})_{n_{l}, 0}
L_{\widehat{M}^{[g]}_{\ell, h}}(-1)^{k}
(a^{j_{1}}_{[g], \ell})_{\alpha^{j_{1}},0}
\cdots (a^{j_{m}}_{[g], \ell})_{\alpha^{j_{m}}, 0}(\psi_{[g]}^{b})_{-1, 0}\one
\end{equation}
for $i_{p}, j_{q}\in I$,
$n_{p}\in \alpha^{i_{p}}-\Z_{+}$, $\Re(\alpha^{j_{q}})=0$ for $p=1, \dots, l$, 
$q=1, \dots, m$, $k\in \N$ and
$b\in A$.}
\end{rema}

Next we construct universal grading-restricted generalized $g$-twisted $M(\ell, 0)$-modules. 
From (\ref{element-form-0-3}), we see that it is impossible for the 
homogeneous subspaces of $U(\hat{\g}^{[g]})
\otimes _{U(\hat{\g}^{[g]}_{+}\oplus \C\mathbf{k})}\Lambda(M)$ to be finite dimensional 
since $a^{i}(0)$ for $a^{i}\in \hat{\g}^{[0]}$ act on $w^{b}$ generate an infinite-dimensional 
homogeneous subspace. But if 
$M$ is a finite-dimensional $\hat{\g}_{\mathbb{I}}$-module,  a quotient of 
$U(\hat{\g}^{[g]})
\otimes _{U(\hat{\g}^{[g]}_{+}\oplus \C\mathbf{k})}\Lambda(M)$ might be grading restricted. 
Since $U(\hat{\g}^{[g]})
\otimes _{U(\hat{\g}^{[g]}_{+}\oplus \C\mathbf{k})}\Lambda(M)$ 
as a lower-bounded $\hat{\g}^{[g]}$-module is equivalent to $\widehat{M}_{\ell}^{[g]}$,
the same discussion applies to $\widehat{M}_{\ell}^{[g]}$. 

Now we assume that  $M$ is in addition a finite-dimensional $\hat{\g}_{\mathbb{I}}$-module
with a compatible action of $g$. Here by $M$ has a 
compatible action of $g$ we mean
$g(a(n)w)=(g(a))(n)g(w)$  for $a\in \g^{[\alpha]}$ such that $\Re(\alpha)=0$, 
$n\in \alpha+\Z$ and $w\in M$.
We have a universal lower-bounded generalized $g$-twisted $M(\ell, 0)$-module
$\widehat{M}^{[g]}_{\ell}$. 
Since $M$ in our case is a $\hat{\g}^{[g]}_{\I}$-module but the construction of $\widehat{M}^{[g]}_{\ell}$ 
above does not use such a structure on $M$, to 
incorporate such an action on $M$, 
we need to take a further quotient. Let $I_{\I}$ be the set of elements $\alpha^{i}$ of $I$ such that 
$\Re(\alpha^{i})=0$. Then $\hat{\g}^{[g]}_{\I}$ is spanned by elements of the form 
$a^{i}(\alpha^{i})$ for $i\in I_{\I}$. Since $\{w^{a}\}_{a\in 
A}$ is a basis of $M$, there exist $\lambda_{ic}^{a}\in \C$
for  $i\in I_{\I}$ and $b, c\in A$
such that 
$$a^{i}(\alpha^{i})w^{b}=\sum_{c\in A}\lambda_{ic}^{b}w^{c}$$
for $i\in I_{\I}$ and $b\in A$. 

Consider the lower-bounded generalized $g$-twisted $M(\ell, 0)$-submodule
of $\widehat{M}^{[g]}_{\ell, h}$ generated by elements of the form
\begin{align}\label{widebreve-relation}
(a^{i}_{[g], \ell})_{\alpha^{i}, 0}
(\psi^{b}_{[g], h^{b}})_{-1, 0}\one
-\sum_{c\in A}\lambda_{ic}^{b}(\psi^{c}_{[g], h^{c}})_{-1, 0}\one
\end{align}
for $i\in I_{\I}$ and $b\in A$.
We denote the quotient of $\widehat{M}^{[g]}_{\ell}$ by this submodule 
by $\widebreve{M}^{[g]}_{\ell}$. Then $\widebreve{M}^{[g]}_{\ell}$ is also a 
lower-bounded generalized $g$-twisted $M(\ell, 0)$-module. 
We shall use the same notations for the generating twisted fields, generator
twist fields and their coefficients for $\widehat{M}^{[g]}_{\ell}$ to denotes the corresponding 
fields and their coefficients for $\widebreve{M}^{[g]}_{\ell}$. 

On the $\hat{\g}^{[g]}_{\I}$-module $M$,
we define 
an action of $\hat{\g}^{[g]}_{+}$ to be $0$. Then we use 
the commutator formula (\ref{comm-a-m-L-1}) 
 for $a\in \g^{[\alpha]}$ and $m\in \alpha+\N$ to define an action of 
$\hat{\g}^{[g]}_{+}\oplus \hat{\g}^{[g]}_{\I}$ on $\Lambda(M)$. 
Let $\mathbf{k}$ act on $\Lambda(M)$ 
as $\ell$. Then $\Lambda(M)$ becomes a $\hat{\g}^{[g]}_{+}\oplus \hat{\g}^{[g]}_{0}$-module
and  we have the induced lower-bounded 
$\hat{\g}^{[g]}$-module $U(\hat{\g}^{[g]})
\otimes_{\hat{\g}^{[g]}_{+}\oplus \hat{\g}^{[g]}_{0}}\Lambda(M)$.
From the construction, we see that the $\hat{\g}^{[g]}$-module $U(\hat{\g}^{[g]})
\otimes_{\hat{\g}^{[g]}_{+}\oplus \hat{\g}^{[g]}_{0}}\Lambda(M)$
is in fact the quotient of $\hat{\g}^{[g]}$-module  $U(\hat{\g}^{[g]})
\otimes _{U(\hat{\g}^{[g]}_{+}\oplus \C\mathbf{k})}\Lambda(M)$ by the 
submodule generated by elements of the form 
\begin{equation}\label{widebreve-relation-2}
a^{i}(\alpha^{i})\otimes w^{b}-\sum_{c\in A}\lambda_{ic}^{b}w^{c}
\end{equation}
for $i\in I_{\I}$ and $b\in A$. 

\begin{thm}\label{widebreve-equiv}
As a lower-bounded $\hat{\g}^{[g]}$-module, $\widebreve{M}^{[g]}_{\ell}$ is equivalent to 
the induced lower-bounded $\hat{\g}^{[g]}$-module $U(\hat{\g}^{[g]})
\otimes_{\hat{\g}^{[g]}_{+}\oplus \hat{\g}^{[g]}_{0}}\Lambda(M)$.
\end{thm}
\pf
By Theorem \ref{ulbtm-va}, $\widehat{M}^{[g]}_{\ell}$ as a 
lower-bounded $\hat{\g}^{[g]}$-module is equivalent to 
$\hat{\g}^{[g]}$-module $U(\hat{\g}^{[g]})
\otimes_{\hat{\g}^{[g]}_{+}\oplus \C\mathbf{k}}\Lambda(M)$.
It is also clear that the submodule of $\widehat{M}^{[g]}_{\ell}$
generated by elements of the form (\ref{widebreve-relation})
for $i\in I_{\I}$ and $b\in A$ and the submodule of 
$\hat{\g}^{[g]}$-module $U(\hat{\g}^{[g]})
\otimes_{\hat{\g}^{[g]}_{+}\oplus \C\mathbf{k}}\Lambda(M)$ 
generated by elements of the form (\ref{widebreve-relation-2}) 
are equivalent under the equivalence from $\widehat{M}^{[g]}_{\ell}$
to $\hat{\g}^{[g]}$-module $U(\hat{\g}^{[g]})
\otimes_{\hat{\g}^{[g]}_{+}\oplus \C\mathbf{k}}\Lambda(M)$. 
Thus their quotients 
$\widebreve{M}^{[g]}_{\ell}$ and $\hat{\g}^{[g]}$-module $U(\hat{\g}^{[g]})
\otimes_{\hat{\g}^{[g]}_{+}\oplus \hat{\g}^{[g]}_{0}}\Lambda(M)$
are equivalent. 
\epfv

\begin{rema}\label{elemenet-order-2}
{\rm From the construction of $U(\hat{\g}^{[g]})
\otimes_{\hat{\g}^{[g]}_{+}\oplus \hat{\g}^{[g]}_{0}}\Lambda(M)$, 
it is spanned by elements of the form
\begin{equation}\label{element-form-0-2.5}
a^{i_{1}}(n_{1})
\cdots a^{i_{l}}(n_{l})
L_{\Lambda(M)}(-1)^{k}w^{b}
\end{equation}
for $i_{p}\in I$,
$n_{p}\in \alpha^{i_{p}}-\Z_{+}$ for $p=1, \dots, l$, $k\in \N$ and
$b\in A$. Similarly, from the construction of $\widebreve{M}^{[g]}_{\ell}$,
it is spanned by elements of the form 
\begin{equation}\label{element-form-0-3.5}
(a^{i_{1}}_{[g], \ell})_{n_{1},0}
\cdots (a^{i_{l}}_{[g], \ell})_{n_{l}, 0}
L_{\widehat{M}^{[g]}_{\ell, h}}(-1)^{k}
(\psi_{[g]}^{b})_{-1, 0}\one
\end{equation}
for $i_{p}, j_{q}\in I$,
$n_{p}\in \alpha^{i_{p}}-\Z_{+}$ for $p=1, \dots, l$, $k\in \N$ and
$b\in A$.}
\end{rema}

We are ready to prove that $\widebreve{M}^{[g]}_{\ell}$ is in fact grading restricted now. 

\begin{thm}\label{grad-rest}
The lower-bounded generalized $g$-twisted $M(\ell, 0)$-module
$\widebreve{M}^{[g]}_{\ell}$ is in fact grading restricted. 
\end{thm}
\pf
By Theorem \ref{widebreve-equiv}, we need only prove that 
$U(\hat{\g}^{[g]})
\otimes_{\hat{\g}^{[g]}_{+}\oplus \hat{\g}^{[g]}_{0}}\Lambda(M)$
 is grading restricted. 
By Remark \ref{elemenet-order-2},  $U(\hat{\g}^{[g]})
\otimes_{\hat{\g}^{[g]}_{+}\oplus \hat{\g}^{[g]}_{0}}\Lambda(M)$
is spanned by elements of the form 
(\ref{element-form-0-2.5})
for $i_{p}\in I$,
$n_{p}\in \alpha^{i_{p}}-\Z_{+}$ for $p=1, \dots, l$, $k\in \N$ and
$b\in A$.  The weight of such an element is $-n_{1}-\cdots -n_{l}+k+\wt w^{a}$.
For fixed $n\in \C$, elements of weight $n$ of the form (\ref{element-form-0-2.5}) must 
satisfy $n=-n_{1}-\cdots -n_{l}+k+\wt w^{a}$. So we have
$$n_{1}+\cdots +n_{l}-k=-n+\wt w^{a}.$$
Since $M$ is finite dimensional, 
there are only finitely many $w^{a}$ and thus finitely many $\wt w^{a}$. 
Let $N\in \R$ such that $\Re(\wt w^{a})\ge N$ for $a\in A$. Then 
$$\Re(n_{1})+\cdots +\Re(n_{l})-k=-\Re(n)+\Re(\wt w^{a})\ge -\Re(n)+N.$$
On the other hand, since $n_{j}\in \alpha^{i_{j}}-\Z_{+}$, 
$\Re(n_{j})<0$ and we obtain
\begin{equation}\label{grad-rest-1}
0>\Re(n_{1})+\cdots +\Re(n_{l})-k\ge -\Re(n)+N.
\end{equation}
Let $P=\max_{i\in I}\{\Re(\alpha^{i})-1\}$. Then 
$P\in [-1, 0)$. Since $n_{j}=\alpha^{i_{j}}-\Z_{+}=\alpha^{i_{j}}-1-\N$,
we have $\Re(n_{j})\le \Re(\alpha^{i_{j}})-1\le P<0$. So
$\Re(n_{1})+\cdots +\Re(n_{l})\le lP$. If $lP< -\Re(n)+N$,
we have $\Re(n_{1})+\cdots +\Re(n_{l})-k<-\Re(n)+N-k\le -\Re(n)+N$. 
Contradiction to (\ref{grad-rest-1}). Thus we must have 
$lP\ge  -\Re(n)+N$ or equivalently, $l\le \frac{1}{P}(-\Re(n)+N)$ (note that $P<0$).
Since $\Re(n_{j})<0$ for $j=1, \dots, l$, from (\ref{grad-rest-1}) and $-k\le 0$,
we obtain also $\Re(n_{j})\ge -\Re(n)+N$ and $-k\ge -\Re(n)+N$. 
From $0>\Re(n_{j})\ge -\Re(n)+N$ for $j=1, \dots, l$ and $0\ge -k\ge  -\Re(n)+N$, 
we see that for fixed $l$,
there are only finitely many 
possible choices of $a^{i_{j}}$, $n_{j}$ and $k$. 
Thus for fixed $n\in \C$, 
there are only finitely many elements of weight $n$ of the form (\ref{element-form-0-2.5}).
So $U(\hat{\g}^{[g]})
\otimes_{\hat{\g}^{[g]}_{+}\oplus \hat{\g}^{[g]}_{\I}}M$, or equivalently,
$\widebreve{M}^{[g]}_{\ell}$  is grading restricted.
\epfv

\subsection{The constructions when $M(\ell, 0)$ is viewed as a vertex operator algebra}\label{4.2}

Assume that $\g$ is simple and $\ell+h^{\vee}\ne 0$. Then $M(\ell, 0)$ has a conformal vector
$\omega_{M(\ell, 0)}$ and thus is a vertex operator algebra.
Now we want to construct and identify explicitly universal lower-bounded generalized $g$-twisted 
modules for $M(\ell, 0)$ viewed as a vertex operator algebra. 
Since in \cite{H-const-twisted-mod}, we give only the construction for a grading-restricted 
vertex algebra or a M\"{o}bius vertex algebra, here we first give a construction of
universal lower-bounded generalized twisted  module
for a general vertex operator algebra. 

Let $V$ be a vertex operator algebra, that is, a grading-restricted vertex algebra $V$
with a conformal element $\omega$, and $g$ an automorphism of $V$ as a vertex operator 
algebra (meaning in particular that $g$ fixes $\omega$). 
Let $M$ be a vector space with actions of 
$g$, $\mathcal{S}_{g}$, $\mathcal{N}_{g}$, $L_{M}(0)$, 
$L_{M}(0)_{S}$ and $L_{M}(0)_{N}$ and $B$ a real number such that $M$ is a direct sum of 
generalized eigenspaces of $L_{M}(0)$ and the real parts of 
the eigenvalues of 
$L_{M}(0)$ are larger than or equal to $B$. From Section 5 of \cite{H-const-twisted-mod},
we have a universal lower-bounded generalized $g$-twisted $V$-module $\widehat{M}_{B}^{[g]}$. 
Since $g$ fix $\omega$, the coefficients of $Y_{\widehat{M}_{B}^{[g]}}^{g}(\omega, x)$
satisfy the Virasoro commutator relations. Note that for a lower-bounded 
generalized $g$-twisted $V$-module $W$, the operator $L_{W}(0)$
and $L_{W}(-1)$
must be equal to 
the coefficients of $x^{-2}$ and $x^{-1}$, respectively,  in 
the vertex operator $Y_{W}(\omega, x)$. But $L_{\widehat{M}_{B}^{[g]}}(0)$ 
and $L_{\widehat{M}_{B}^{[g]}}(-1)$ for $\widehat{M}_{B}^{[g]}$ are not equal to 
the coefficients of $x^{-2}$ and $x^{-1}$, respectively. To obtain a lower-bounded 
generalized $g$-twisted module for $V$ viewed as a vertex operator algebra, 
we have to take the quotient by a submodule generated by the difference of these operators acting on 
elements of $\widehat{M}^{[g]}_{B}$.  

Consider the lower-bounded generalized $g$-twisted $V$-submodule
of $\widehat{M}_{B}^{[g]}$
generated by elements of the form
\begin{align*}
&L_{\widehat{M}^{[g]}_{B}}(0)w-\res_{x}xY_{\widehat{M}^{[g]}_{B}}^{g}(\omega, x)w,\\
&L_{\widehat{M}^{[g]}_{B}}(-1)w-\res_{x}Y_{\widehat{M}^{[g]}_{B}}^{g}(\omega, x)w
\end{align*}
for $w\in \widehat{M}^{[g]}_{B}$. We shall denote the quotient of $\widehat{M}^{[g]}_{B}$
by this submodule by  
$$\overarc{M}^{[g]}_{B}$$
and call this quotient module the 
lower-bounded generalized $g$-twisted $V$-module for the vertex operator algebra $V$, not the underlying
grading-restricted vertex algebra $V$. 
By Theorem 5.2 and the construction of $\widehat{M}^{[g]}_{B}$
in \cite{H-const-twisted-mod}, 
we immediately obtain the following result:

\begin{thm}\label{universal-voa}
Let $V$ be a vertex operator algebra and 
$(W, Y^{g}_{W})$ a lower-bounded generalized $g$-twisted $V$-module
and $M^{0}$ a  subspace of $W$ invariant under the actions of 
$g$, $\mathcal{S}_{g}$, $\mathcal{N}_{g}$, $L_{W}(0)=
\res_{x}xY_{W}^{g}(\omega, x)$, 
$L_{W}(0)_{S}$ and
$L_{W}(0)_{N}$. Let $B\in \R$ such that $W_{[n]}=0$ when $\Re(n)<B$. 
Assume that there is a linear map $f: M\to M^{0}$ 
 commuting with the actions of 
$g$, $\mathcal{S}_{g}$, $\mathcal{N}_{g}$, $L_{W}(0)|_{M^{0}}$ and $L_{M}(0)$ , 
$L_{W}(0)_{S}|_{M^{0}}$ and $L_{M}(0)_{S}$ and $L_{W}(0)_{N}|_{M^{0}}$ and
$L_{M}(0)_{N}$. Then there exists a unique module 
map $\overarc{f}: \overarc{M}^{[g]}_{B}\to W$ such that $\overarc{f}|_{M}=f$. 
If $f$ is surjective and $(W, Y^{g}_{W})$ is generated by the 
coefficients of $(Y^{g})_{WV}^{W}(w, x)v$ for $w\in M_{0}$ and $v\in V$, 
where $(Y^{g})_{WV}^{W}$ 
is the twist vertex operator map obtained from $Y_{W}^{g}$, then 
$\overarc{f}$ is surjective. 
\end{thm}

We now assume that $\g$ is simple and $\ell+h^{\vee}\ne 0$. Then $M(\ell, 0)$ is a vertex operator
algebra. Take the vertex operator algebra $V$ above to
be $M(\ell, 0)$ and $g$ an automorphism of $M(\ell, 0)$ induced 
from an automorphism of $\g$ as 
discussed above.  Let $M$, as above, be a vector space with actions of 
$g$,  $\mathcal{L}_{g}$,
$\mathcal{S}_{g}$, $\mathcal{N}_{g}$, $L_{M}(0)$, $L_{M}(0)_{S}$ and $L_{M}(0)_{N}$.
We assume that 
$M=\coprod_{h\in Q_{M}}M_{[h]}$ as above, where 
$M_{[h]}$ is the generalized eigenspace of $L_{M}(0)$ with eigenvalue $h$ and 
$Q_{M}$ is the set of all eigenvalues of $L_{M}(0)$. 
For $h\in Q_{h}$, take $V$, $g$, $M$ and $B$ in the construction above to be 
$M(\ell, 0)$, $g$, $M_{[h]}$, $\Re(h)$. Then we have a universal
lower-bounded generalized $g$-twisted $M(\ell, 0)$-module
$\overarc{(M_{[h]})}_{h}^{[g]}$. To exhibit its dependence on $\ell$ explicitly, we 
denote it by $\overarc{(M_{[h]})}_{\ell, h}^{[g]}$.
Adding them together, we obtain a 
lower-bounded generalized $g$-twisted $M(\ell, 0)$-module
$\coprod_{h\in Q_{M}}\overarc{(M_{[h]})}_{\ell, h}^{[g]}$, which we shall denote 
by $\overarc{M}^{[g]}_{\ell}$. 
We shall use the same notations $a^{i}_{[g], \ell}(x)$, 
$(a^{i}_{[g], \ell})_{n, k}$, $\psi^{b}_{[g]}(x)$
and $(\psi^{b}_{[g]})_{n, k}$ and so on 
as those for $\widehat{M}^{[g]}_{\ell}$ and 
$\widebreve{M}^{[g]}_{\ell}$
to denote the generating twisted fields, their coefficients,
the generator twist fields and their coefficients for $\overarc{M}^{[g]}_{\ell}$. 

We need to identify $\overarc{M}^{[g]}_{\ell}$ with a suitable $\hat{\g}^{[g]}$-module.
We first need to identify $L_{M}(0)$ with the
action of an element  of $U(\hat{\g}^{[g]})$.
Recall the Jordan basis $\{a^{i}\}_{i\in I}$ of $\g$ that we have chosen in the end of the preceding 
section. Let  $\{(a^{i})'\}_{i\in I}$ be the dual basis of  $\{a^{i}\}_{i\in I}$ with respect to the 
nondegenerate bilinear form $(\cdot, \cdot)$. For simplicity (but with an abuse of the 
notation), we shall denote this dual 
basis by $\{a^{i'}\}_{i\in I}$. Then 
$$(e^{2\pi i \mathcal{S}_{g}}a^{i'}, a^{j})=(a^{i'}, e^{-2\pi i \mathcal{S}_{g}}a^{j})
=(a^{i'}, e^{-2\pi i \alpha^{j}}a^{j})=e^{-2\pi i \alpha^{j}}\delta_{ij}=e^{-2\pi i \alpha^{i}}\delta_{ij}.$$
This means that 
$$e^{2\pi i \mathcal{S}_{g}}a^{i'}=e^{-2\pi i \alpha^{i}}a^{i'}.$$
So $a^{i'}\in \g^{[1-\alpha^{i}]}$ when $\Re(\alpha^{i})>0$ or $i\in I\setminus I_{\I}$ and 
$a^{i'}\in \g^{[-\alpha^{i}]}$ when $\Re(\alpha^{i})=0$ or $i\in I_{\I}$. By abuse of 
notation, let 
$$\alpha^{i'}=\left\{\begin{array}{ll}1-\alpha^{i}&i\in I\setminus I_{\I},\\
-\alpha^{i}&i\in I_{\I}.\end{array}\right.$$
By definition, the conformal element of $M(\ell, 0)$ is
$$\omega_{M(\ell, 0)}=\sum_{i\in I}a^{i'}(-1)a^{i}(-1)\one
\in M^{[0]}(\ell, 0),$$
where $M^{[0]}(\ell, 0)$ is the fixed-point subalgebra of $M(\ell, 0)$. 

We need to recall the Virasoro operators on $M(\ell, 0)$.
Since $\omega_{M(\ell, 0)}$ is in the fixed-point subalgebra of $M(\ell, 0)$,
$\mathcal{N}_{g}\omega_{M(\ell, 0)}=0$. Hence 
$$Y^{g}_{\overarc{M}_{\ell}^{[g]}}(\omega_{M(\ell, 0)}, x)
=(Y^{g}_{\overarc{M}_{\ell}^{[g]}})_{0}(x^{-\mathcal{N}_{g}}\omega_{M(\ell, 0)}, x)
=(Y^{g}_{\overarc{M}_{\ell}^{[g]}})_{0}(\omega_{M(\ell, 0)}, x)$$
and from the equivariance property of the twisted vertex operators, 
$Y^{g}_{\overarc{(M_{[h]})}_{h}^{[g]}}(\omega, x)$ 
or equivalently $(Y^{g}_{\overarc{(M_{[h]})}_{h}^{[g]}})_{0}(\omega, x)$
must have only integral powers of $x$. In particular, 
$$Y^{g}_{\overarc{M}_{\ell}^{[g]}}(\omega, x)
=(Y^{g}_{\overarc{M}_{\ell}^{[g]}})_{0}(\omega, x)
=\sum_{n\in \Z}L_{\overarc{M}_{\ell}^{[g]}}(n)x^{-n-2}.$$
where $L_{\overarc{M}_{\ell}^{[g]}}(n)$ for $n\in \Z$ are the Virasoro operators on 
$\overarc{M}_{\ell}^{[g]}$ satisfying the Virasoro commutator relations 
with central charge $\frac{\ell \dim \g}{\ell+h^{\vee}}$. In particular, we have the operators 
$L_{\overarc{M}_{\ell}^{[g]}}(0)$ and $L_{\overarc{M}_{\ell}^{[g]}}(-1)$.

\begin{prop}
For $n\in \Z$, 
\begin{align}\label{twisted-vir-comm-2}
L&_{\overarc{M}_{\ell}^{[g]}}(n)\nn
&=\sum_{i\in I}\sum_{p\in \alpha^{i}+\Z_{+}}\frac{1}{2(\ell+h^{\vee})}
a^{i'}_{[g], \ell}(-p)a^{i}_{[g], \ell}(p+n) + \sum_{i\in I}\sum_{p\in \alpha^{i}-\N}\frac{1}{2(\ell+h^{\vee})}
a^{i}_{[g], \ell}(p+n)a^{i'}_{[g], \ell}(-p)\nn
&\quad -\sum_{i\in I}\frac{1}{2(\ell+h^{\vee})}
[(\mathcal{N}_{g}-\alpha^{i})
a^{i'}, a^{i}]_{[g], \ell}(n) -\sum_{i\in I}\frac{\ell\delta_{n, 0}}{4(\ell+h^{\vee})}
((\mathcal{N}_{g}-\alpha^{i})
(\mathcal{N}_{g}-\alpha^{i}-1)
a^{i'}, a^{i}).
\end{align}
\end{prop}
\pf
We take $W=\overarc{M}_{\ell}^{[g]}$,
$u=x_{1}^{\mathcal{N}_{g}}
a^{i'}(-1)\one$ and $v=x_{2}^{\mathcal{N}_{g}}a^{i}(-1)\one$ 
in the Jacobi identity (\ref{jacobi}).
Then $Y^{g}_{\overarc{M}_{\ell}^{[g]}}(u, x_{1})=(Y^{g}_{\overarc{M}_{\ell}^{[g]}})_{0}(a^{i'}(-1)\one, x_{1})$,
$Y^{g}_{\overarc{M}_{\ell}^{[g]}}(v, x_{2})=(Y^{g}_{\overarc{M}_{\ell}^{[g]}})_{0}(a^{i}(-1)\one, x_{2})$
 and 
$\mathcal{S}_{g}a^{i'}(-1)\one=\alpha^{i'}a^{i'}(-1)\one$. Then 
(\ref{jacobi}) becomes
\begin{align}\label{jacobi-a-i}
x&_0^{-1}\delta\left(\frac{x_1 - x_2}{x_0}\right)
(Y_{\overarc{M}^{[g]}_{\ell}}^{g})_{0}(a^{i'}(-1)\one, x_1)
(Y_{\overarc{M}^{[g]}_{\ell}}^{g})_{0}(a^{i}(-1)\one, x_2)\nn
&\quad -  x_0^{-1}\delta\left(\frac{- x_2 + x_1}{x_0}\right)
(Y_{\overarc{M}^{[g]}_{\ell}}^{g})_{0}(a^{i}(-1)\one, x_2)
(Y_{\overarc{M}^{[g]}_{\ell}}^{g})_{0}(a^{i'}(-1)\one, x_1)\nn
&= x_1^{-1}\delta\left(\frac{x_2+x_0}{x_1}\right)
\left(\frac{x_2+x_0}{x_1}\right)^{-\alpha^{i}}\cdot\nn
&\quad\quad\quad\cdot 
(Y_{\overarc{M}^{[g]}_{\ell}}^{g})_{0}\left(Y_{M(\ell, 0)}
\left(\left(1+\frac{x_0}{x_2}\right)^{\mathcal{N}_{g}}
a^{i'}(-1)\one, x_0\right)a^{i}(-1)\one, x_2\right),
\end{align}
where we have used 
$$x_1^{-1}\delta\left(\frac{x_2+x_0}{x_1}\right)
\left(\frac{x_2+x_0}{x_1}\right)^{\alpha^{i'}}=x_1^{-1}\delta\left(\frac{x_2+x_0}{x_1}\right)
\left(\frac{x_2+x_0}{x_1}\right)^{-\alpha^{i}}.$$
Multiplying both sides of 
(\ref{jacobi-a-i}) by $x_{1}^{-\alpha^{i}}$ and 
then take $\res_{x_{1}}$, rewriting $(x_2+x_0)^{-\alpha^{i}}$ as 
$x_{2}^{-\alpha^{i}}\left(1+\frac{x_0}{x_2}\right)^{-\alpha^{i}}$ and then multiplying
both sides by $x_{2}^{\alpha^{i}}$, we obtain 
\begin{align}\label{jacobi-a-i-1}
\res&_{x_{1}}x_{1}^{-\alpha^{i}}x_{2}^{\alpha^{i}}
x_0^{-1}\delta\left(\frac{x_1 - x_2}{x_0}\right)
(Y_{\overarc{M}^{[g]}_{\ell}}^{g})_{0}(a^{i'}(-1)\one, x_1)
(Y_{\overarc{M}^{[g]}_{\ell}}^{g})_{0}(a^{i}(-1)\one, x_2)\nn
&\quad -  \res_{x_{1}}x_{1}^{-\alpha^{i}}x_{2}^{\alpha^{i}}
x_0^{-1}\delta\left(\frac{- x_2 + x_1}{x_0}\right)
(Y_{\overarc{M}^{[g]}_{\ell}}^{g})_{0}(a^{i}(-1)\one, x_2)
(Y_{\overarc{M}^{[g]}_{\ell}}^{g})_{0}(a^{i'}(-1)\one, x_1)\nn
&=\left(1+\frac{x_0}{x_2}\right)^{-\alpha^{i}}
(Y_{\overarc{M}^{[g]}_{\ell}}^{g})_{0}\left(Y_{M(\ell, 0)}
\left(\left(1+\frac{x_0}{x_2}\right)^{\mathcal{N}_{g}}
a^{i'}(-1)\one, x_0\right)a^{i}(-1)\one, x_2\right).
\end{align}
Using the definition, we have 
$$(Y_{\overarc{M}^{[g]}_{\ell}}^{g})_{0}(a(-1)\one, x_1)=(a_{[g], \ell})_{0}(x),$$
where 
$$(a_{[g], \ell})_{0}(x)=\sum_{n\in \alpha+\Z}(a_{[g], \ell})_{n, 0}x^{-n-1}$$
for $a\in \g^{[\alpha]}$. 
Then the constant term in $x_{0}$ (or equivalently, 
the result of applying $\res_{x_{0}}x_{0}^{-1}$) 
of the right-hand side of (\ref{jacobi-a-i-1}) is equal to 
\begin{align*}
&\res_{x_{0}}x_{0}^{-1}(Y_{\overarc{M}^{[g]}_{\ell}}^{g})_{0}\left(Y_{M(\ell, 0)}
\left(\left(1+\frac{x_0}{x_2}\right)^{\mathcal{N}_{g}-\alpha^{i}}
a^{i'}(-1)\one, x_0\right)a^{i}(-1)\one, x_2\right)\nn
&\quad=\sum_{m\in \N}\res_{x_{0}}x_{0}^{-1}\left(\frac{x_{0}}{x_{2}}\right)^{m}
(Y_{\overarc{M}^{[g]}_{\ell}}^{g})_{0}\left(Y_{M(\ell, 0)}
\left(\binom{\mathcal{N}_{g}-\alpha^{i}}{m}
a^{i'}(-1)\one, x_0\right)a^{i}(-1)\one, x_2\right)\nn
\end{align*}
\begin{align}\label{jacobi-a-i-2}
&\quad=\sum_{m\in \N}\res_{x_{0}}x_{0}^{-1}\left(\frac{x_{0}}{x_{2}}\right)^{m}
(Y_{\overarc{M}^{[g]}_{\ell}}^{g})_{0}\left(\left(\binom{\mathcal{N}_{g}-\alpha^{i}}{m}
a^{i'}\right)(x_0)a^{i}(-1)\one, x_2\right)\nn
&\quad=\sum_{m\in \N}\sum_{n\in \Z}
\res_{x_{0}}x_{0}^{-1}\left(\frac{x_{0}}{x_{2}}\right)^{m}x_{0}^{-n-1}
(Y_{\overarc{M}^{[g]}_{\ell}}^{g})_{0}\left(\left(\binom{\mathcal{N}_{g}-\alpha^{i}}{m}
a^{i'}\right)(n)a^{i}(-1)\one, x_2\right)\quad\;\;\nn
&\quad=\sum_{m\in \N}
x_{2}^{-m}(Y_{\overarc{M}^{[g]}_{\ell}}^{g})_{0}\left(\left(\binom{\mathcal{N}_{g}-\alpha^{i}}{m}
a^{i'}\right)(m-1)a^{i}(-1)\one, x_2\right)\nn
&\quad=(Y_{\overarc{M}^{[g]}_{\ell}}^{g})_{0}(a^{i'}(-1)a^{i}(-1)\one, x_2)\nn
&\quad\quad  +x_{2}^{-1}
(Y_{\overarc{M}^{[g]}_{\ell}}^{g})_{0}(((\mathcal{N}_{g}-\alpha^{i})
a^{i'})(0)a^{i}(-1)\one, x_2)\nn
&\quad\quad  +\frac{x_{2}^{-2}}{2}
(Y_{\overarc{M}^{[g]}_{\ell}}^{g})_{0}(((\mathcal{N}_{g}-\alpha^{i})
(\mathcal{N}_{g}-\alpha^{i}-1)
a^{i'})(1)a^{i}(-1)\one, x_2)\nn
&\quad=(Y_{\overarc{M}^{[g]}_{\ell}}^{g})_{0}(a^{i'}(-1)a^{i}(-1)\one, x_2)
+x_{2}^{-1}
([(\mathcal{N}_{g}-\alpha^{i})
a^{i'}, a^{i}]_{[g], \ell})_{0}(x_2)\nn
&\quad\quad  +\frac{\ell x_{2}^{-2}}{2}
((\mathcal{N}_{g}-\alpha^{i})
(\mathcal{N}_{g}-\alpha^{i}-1)
a^{i'}, a^{i}).
\end{align}

Applying $\res_{x_{0}}x_{0}^{-1}$ to both sides of 
(\ref{jacobi-a-i-1}), using (\ref{jacobi-a-i-2}), taking sum over 
$i\in I$ on both sides and dividing both sides by $2(\ell+h^{\vee})$, we obtain 
\begin{align*}
&\sum_{n\in \Z}L_{\overarc{M}^{[g]}_{\ell}}(n)x_{2}^{-n-2}\nn
&\quad =(Y^{g}_{\overarc{M}^{[g]}_{\ell}})_{0}(\omega, x_{2})\nn
&\quad =\sum_{i\in I}\frac{1}{2(\ell+h^{\vee})}
(Y_{\overarc{M}^{[g]}_{\ell}}^{g})_{0}(a^{i'}(-1)a^{i}(-1)\one, x_2)\nn
&\quad =\sum_{i\in I}\frac{1}{2(\ell+h^{\vee})}
\res_{x_{0}}x_{0}^{-1}\res_{x_{1}}x_{1}^{-\alpha^{i}}x_{2}^{\alpha^{i}}
x_0^{-1}\delta\left(\frac{x_1 - x_2}{x_0}\right)
(a^{i'}_{[g], \ell})_{0}(x_1)
(a^{i}_{[g], \ell})_{0}(x_2)\nn
&\quad\quad  - \sum_{i\in I}\frac{1}{2(\ell+h^{\vee})}
\res_{x_{0}}x_{0}^{-1} \res_{x_{1}}x_{1}^{-\alpha^{i}}x_{2}^{\alpha^{i}}
x_0^{-1}\delta\left(\frac{- x_2 + x_1}{x_0}\right)
(a^{i}_{[g], \ell})_{0}(x_2)
(a^{i'}_{[g], \ell})_{0}(x_1)\nn
&\quad\quad -\sum_{i\in I}\frac{x_{2}^{-1}}{2(\ell+h^{\vee})}
([(\mathcal{N}_{g}-\alpha^{i})
a^{i'}, a^{i}]_{[g], \ell})_{0}(x_2)\nn
&\quad\quad 
-\sum_{i\in I}\frac{\ell x_{2}^{-2}}{4(\ell+h^{\vee})}
((\mathcal{N}_{g}-\alpha^{i})
(\mathcal{N}_{g}-\alpha^{i}-1)
a^{i'}, a^{i})\nn
&\quad =\sum_{i\in I}\frac{1}{2(\ell+h^{\vee})}
\res_{x_{1}}x_{1}^{-\alpha^{i}}x_{2}^{\alpha^{i}}(x_1 - x_2)^{-1}
(a^{i'}_{[g], \ell})_{0}(x_1)
(a^{i}_{[g], \ell})_{0}(x_2)\nn
&\quad\quad  - \sum_{i\in I}\frac{1}{2(\ell+h^{\vee})}
 \res_{x_{1}}x_{1}^{-\alpha^{i}}x_{2}^{\alpha^{i}}
(- x_2 + x_1)^{-1}
(a^{i}_{[g], \ell})_{0}(x_2)
(a^{i'}_{[g], \ell})_{0}(x_1)
\end{align*}
\begin{align}\label{twisted-vir-comm-1}
&\quad\quad -\sum_{i\in I}\frac{x_{2}^{-1}}{2(\ell+h^{\vee})}
([(\mathcal{N}_{g}-\alpha^{i})
a^{i'}, a^{i}]_{[g], \ell})_{0}(x_2)\nn
&\quad\quad 
-\sum_{i\in I}\frac{\ell x_{2}^{-2}}{4(\ell+h^{\vee})}
((\mathcal{N}_{g}-\alpha^{i})
(\mathcal{N}_{g}-\alpha^{i}-1)
a^{i'}, a^{i}).\quad\quad\quad\quad\quad\quad\quad\quad\quad\quad\quad
\end{align}
Taking the coefficients of $x_{2}^{-n-2}$ of (\ref{twisted-vir-comm-1}),
we obtain
\begin{align*}
L&_{\overarc{M}_{\ell}^{[g]}}(n)\nn
&=\sum_{i\in I}\sum_{m\in \Z_{+}}\sum_{k\in -\alpha^{i}+\Z}\sum_{l\in \alpha^{i}+\Z}
\frac{1}{2(\ell+h^{\vee})}\res_{x_{1}}\res_{x_{2}}x_{1}^{-\alpha^{i}-m-k-2}x_{2}^{\alpha^{i}+m+n-l}
a^{i'}_{[g], \ell}(k)a^{i}_{[g], \ell}(l)\nn
&\quad +\sum_{i\in I}\sum_{m\in \Z_{+}}\sum_{k\in -\alpha^{i}+\Z}\sum_{l\in \alpha^{i}+\Z}
\frac{1}{2(\ell+h^{\vee})} \res_{x_{1}}\res_{x_{2}}x_{1}^{-\alpha^{i}+m-k}x_{2}^{\alpha^{i}+n-m-l-1}
a^{i}_{[g], \ell}(l)
a^{i'}_{[g], \ell}(k)\nn
&\quad -\sum_{i\in I}\frac{1}{2(\ell+h^{\vee})}
[(\mathcal{N}_{g}-\alpha^{i})
a^{i'}, a^{i}]_{[g], \ell}(n) -\sum_{i\in I}\frac{\ell\delta_{n, 0}}{4(\ell+h^{\vee})}
((\mathcal{N}_{g}-\alpha^{i})
(\mathcal{N}_{g}-\alpha^{i}-1)
a^{i'}, a^{i})\nn
&=\sum_{i\in I}\sum_{p\in \alpha^{i}+\Z_{+}}\frac{1}{2(\ell+h^{\vee})}
a^{i'}_{[g], \ell}(-p)a^{i}_{[g], \ell}(p+n) + \sum_{i\in I}\sum_{p\in \alpha^{i}-\N}\frac{1}{2(\ell+h^{\vee})}
a^{i}_{[g], \ell}(p+n)a^{i'}_{[g], \ell}(-p)\nn
&\quad -\sum_{i\in I}\frac{1}{2(\ell+h^{\vee})}
[(\mathcal{N}_{g}-\alpha^{i})
a^{i'}, a^{i}]_{[g], \ell}(n) -\sum_{i\in I}\frac{\ell\delta_{n, 0}}{4(\ell+h^{\vee})}
((\mathcal{N}_{g}-\alpha^{i})
(\mathcal{N}_{g}-\alpha^{i}-1)
a^{i'}, a^{i}),
\end{align*}
proving  (\ref{twisted-vir-comm-2}).
\epfv

From (\ref{twisted-vir-comm-2}), we obtain
\begin{align}\label{twisted-L-0}
L&_{\overarc{M}^{[g]}_{\ell}}(0)\nn
&=\sum_{i\in I}\sum_{p\in \alpha^{i}+\Z_{+}}\frac{1}{2(\ell+h^{\vee})}
a^{i'}_{[g], \ell}(-p)a^{i}_{[g], \ell}(p) +\sum_{i\in I}\sum_{p\in \alpha^{i}-\N}\frac{1}{2(\ell+h^{\vee})}
a^{i}_{[g], \ell}(p)a^{i'}_{[g], \ell}(-p)\nn
&\quad -\sum_{i\in I}\frac{1}{2(\ell+h^{\vee})}
[(\mathcal{N}_{g}-\alpha^{i})
a^{i'}, a^{i}]_{[g], \ell}(0) -\sum_{i\in I}\frac{\ell}{4(\ell+h^{\vee})}
((\mathcal{N}_{g}-\alpha^{i})
(\mathcal{N}_{g}-\alpha^{i}-1)
a^{i'}, a^{i})
\end{align}
and 
\begin{align}\label{twisted-L--1}
L&_{\overarc{M}^{[g]}_{\ell}}(-1)\nn
&=\sum_{i\in I}\sum_{p\in \alpha^{i}+\Z_{+}}\frac{1}{2(\ell+h^{\vee})}
a^{i'}_{[g], \ell}(-p)a^{i}_{[g], \ell}(p-1) + \sum_{i\in I}\sum_{p\in \alpha^{i}-\N}\frac{1}{2(\ell+h^{\vee})}
a^{i}_{[g], \ell}(p-1)a^{i'}_{[g], \ell}(-p)\nn
&\quad -\sum_{i\in I}\frac{1}{2(\ell+h^{\vee})}
[(\mathcal{N}_{g}-\alpha^{i})
a^{i'}, a^{i}]_{[g], \ell}(-1).
\end{align}
Note that for $b\in A$, 
\begin{align}\label{twisted-L-0-psi}
L_{\overarc{M}^{[g]}_{\ell}}(0)(\psi^{b}_{[g]})_{-1, 0}\one
&=\sum_{i\in I_{\I}}\frac{1}{2(\ell+h^{\vee})}
a^{i}_{[g], \ell}(\alpha^{i})a^{i'}_{[g], \ell}(-\alpha^{i})(\psi^{b}_{[g]})_{-1, 0}\one\nn
&\quad  -\sum_{i\in I}
\frac{1}{2(\ell+h^{\vee})}[(\mathcal{N}_{g}-\alpha^{i})
a^{i'}, a^{i}]_{[g], \ell}(0)(\psi^{b}_{[g]})_{-1, 0}\one \nn
&\quad -\sum_{i\in I}\frac{\ell}{4(\ell+h^{\vee})}
((\mathcal{N}_{g}-\alpha^{i})
(\mathcal{N}_{g}-\alpha^{i}-1)
a^{i'}, a^{i})(\psi^{b}_{[g]})_{-1, 0}\one,
\end{align}
that is, as an operator on the subspace of $\overarc{M}^{[g]}_{\ell}$ spanned by $(\psi^{b}_{[g]})_{-1, 0}\one$,
$L_{\overarc{M}^{[g]}_{\ell}}(0)$ is equal to 
\begin{align*}
&\sum_{i\in I_{\I}}\frac{1}{2(\ell+h^{\vee})}
a^{i}_{[g], \ell}(\alpha^{i})a^{i'}_{[g], \ell}(-\alpha^{i}) -\sum_{i\in I}\frac{1}{2(\ell+h^{\vee})}
[(\mathcal{N}_{g}-\alpha^{i})
a^{i'}, a^{i}]_{[g], \ell}(0) \nn
&\quad\quad \quad\quad -\sum_{i\in I}\frac{\ell}{4(\ell+h^{\vee})}
((\mathcal{N}_{g}-\alpha^{i})
(\mathcal{N}_{g}-\alpha^{i}-1)
a^{i'}, a^{i}).
\end{align*}
For $h\in Q_{M}$ and $b\in A_{h}$, by definition, 
$$L_{\overarc{M}^{[g]}_{\ell}}(0)(\psi^{b}_{[g]})_{-1, 0}\one=h(\psi^{b}_{[g]})_{-1, 0}\one.$$
Together with (\ref{twisted-L-0-psi}), we obtain the relation
\begin{align}\label{twisted-L-0-relation}
h(\psi^{b}_{[g]})_{-1, 0}\one
&=\sum_{i\in I_{\I}}\frac{1}{2(\ell+h^{\vee})}
a^{i}_{[g], \ell}(\alpha^{i})a^{i'}_{[g], \ell}(-\alpha^{i})(\psi^{b}_{[g]})_{-1, 0}\one\nn
&\quad -\sum_{i\in I}
\frac{1}{2(\ell+h^{\vee})}[(\mathcal{N}_{g}-\alpha^{i})
a^{i'}, a^{i}]_{[g], \ell}(0)(\psi^{b}_{[g]})_{-1, 0}\one \nn
&\quad -\sum_{i\in I}\frac{\ell}{4(\ell+h^{\vee})}
((\mathcal{N}_{g}-\alpha^{i})
(\mathcal{N}_{g}-\alpha^{i}-1)
a^{i'}, a^{i})(\psi^{b}_{[g]})_{-1, 0}\one.
\end{align}

Let 
\begin{align}\label{twisted-affine-L-0}
\Omega^{[g]}
&= \sum_{i\in I_{\I}}
a^{i}(\alpha^{i})a^{i'}(-\alpha^{i})
 -\sum_{i\in I}
[(\mathcal{N}_{g}-\alpha^{i})
a^{i'}, a^{i}](0)\nn
&\quad -\sum_{i\in I}\frac{\ell}{2}
((\mathcal{N}_{g}-\alpha^{i})
(\mathcal{N}_{g}-\alpha^{i}-1)
a^{i'}, a^{i})\nn
&\in U(\hat{\g}^{[g]})^{[0]},
\end{align}
where $U(\hat{\g}^{[g]})^{[0]}$ is the fixed-point subspace of $U(\hat{\g}^{[g]})$ under $g$. 
Let $\Omega^{[g]}$, $\hat{\g}_{+}^{[g]}$ and $\mathbf{k}$ act on $M$ as $L_{M}(0)$, $0$ and $\ell$,
respectively. Let $G(\Omega^{[g]}, \hat{\g}_{+}^{[g]}, \mathbf{k})$ be the 
subalgebra of $U(\hat{\g}^{[g]})$ generated by $\Omega^{[g]}$, $\hat{\g}_{+}^{[g]}$ and $\mathbf{k}$.
Then we have the induced lower-bounded $\hat{\g}_{+}^{[g]}$-module 
$U(\hat{\g}^{[g]})\otimes_{G(\Omega^{[g]}, \hat{\g}_{+}^{[g]}, \mathbf{k})}M$. Note that 
this induced $\hat{\g}_{+}^{[g]}$-module is a quotient of $U(\hat{\g}^{[g]})
\otimes _{U(\hat{\g}^{[g]}_{+}\oplus \C\mathbf{k})}M$. 

\begin{thm}\label{ulbtm-voa}
The universal lower-bounded generalized $g$-twisted $M(\ell, 0)$-module
$\overarc{M}_{\ell}^{[g]}$ 
is equivalent  as a lower-bounded $\hat{\g}^{[g]}$-module  to 
$U(\hat{\g}^{[g]})\otimes_{G(\Omega^{[g]}, \hat{\g}_{+}^{[g]}, \mathbf{k})}M$.
\end{thm}
\pf
We know that $\overarc{M}_{\ell}^{[g]}$ as a quotient of  $\widehat{M}_{\ell}^{[g]}$  is also spanned by 
elements of the form (\ref{element-form-0}) for $i_{j}\in I$,
$n_{j}\in \alpha^{i_{j}}+\Z$ for $j=1, \dots, l$, $k\in \N$ and
$b\in A$. Using (\ref{twisted-L--1}), we see that 
elements of the form 
(\ref{element-form-0}) in $\overarc{M}_{\ell}^{[g]}$ 
can be written as linear combinations of elements of the form
\begin{equation}\label{element-form-0-4}
(a^{i_{1}}_{[g], \ell})_{n_{1},0}
\cdots (a^{i_{l}}_{[g], \ell})_{n_{l}, 0}
(\psi_{[g]}^{b})_{-1, 0}\one
\end{equation}
for $i_{j}\in I$,
$n_{j}\in \alpha^{i_{j}}+\Z$ for $j=1, \dots, l$ and
$b\in A$.
On the other hand,
by definition, $U(\hat{\g}^{[g]})\otimes_{G(\Omega^{[g]}, \hat{\g}_{+}^{[g]}, \mathbf{k})}M$
is spanned by elements of the form
\begin{equation}\label{element-form-0-5}
a^{i_{1}}(n_{1})
\cdots a^{i_{l}}(n_{l})w^{b}
\end{equation}
for $i_{j}\in I$,
$n_{j}\in \alpha^{i_{j}}+\Z$ for $j=1, \dots, l$ and
$b\in A$. 

By Theorem \ref{ulbtm-va}, we have an invertible $\hat{\g}^{[g]}$-module map $\hat{\rho}: 
U(\hat{\g}^{[g]})
\otimes _{U(\hat{\g}^{[g]}_{+}\oplus \C\mathbf{k})}\Lambda(M) \to \widehat{M}_{\ell}^{[g]}$
such that $\hat{\rho}$ maps the element  (\ref{element-form-0}) 
to the element  (\ref{element-form-0-1}). In particular, 
$\hat{\rho}$ maps the  element  of  $U(\hat{\g}^{[g]})
\otimes _{U(\hat{\g}^{[g]}_{+}\oplus \C\mathbf{k})}\Lambda(M)$
of the same form as (\ref{element-form-0-5})  to the
element (\ref{element-form-0-4}). We want to use the map $\hat{\rho}$ restricted to elements of the same form 
as (\ref{element-form-0-5})  to obtain an invertible  
$\hat{\g}^{[g]}$-module map from $U(\hat{\g}^{[g]})
\otimes_{G(\Omega^{[g]}, \hat{\g}_{+}^{[g]}, \mathbf{k})}M$ to $\overarc{M}_{\ell}^{[g]}$. 

To do this, we need only prove that the relations among elements of 
$U(\hat{\g}^{[g]})
\otimes _{U(\hat{\g}^{[g]}_{+}\oplus \C\mathbf{k})}\Lambda(M)$ of the same form as 
(\ref{element-form-0-5})  and the relations among elements of the form 
(\ref{element-form-0-4}) in $\overarc{M}_{\ell}^{[g]}$ are the same. 
The relations among elements of $U(\hat{\g}^{[g]})
\otimes _{U(\hat{\g}^{[g]}_{+}\oplus \C\mathbf{k})}\Lambda(M)$ of the same form as 
(\ref{element-form-0-5}) are generated by the 
following two types: The first type of relations are
induced from the relations in $U(\hat{\g}^{[g]})
\otimes _{U(\hat{\g}^{[g]}_{+}\oplus \C\mathbf{k})}\Lambda(M)$. The second type is 
the additional relations $\frac{1}{2(\ell+h^{\vee})}\Omega^{[g]}w^{b}=L_{M}(0)w^{b}$ 
for $b\in A$. For $h\in Q_{M}$ and $b\in A_{h}$, 
this additional relations become
\begin{align}\label{t-affine-L-0-relation}
h w^{b}
&=\sum_{i\in I_{\I}}\frac{1}{2(\ell+h^{\vee})}
a^{i}(\alpha^{i})a^{i'}(-\alpha^{i})w^{b}
 -\sum_{i\in I}\frac{1}{2(\ell+h^{\vee})}
[(\mathcal{N}_{g}-\alpha^{i})
a^{i'}, a^{i}](0)w^{b}\nn
&\quad -\sum_{i\in I}\frac{\ell}{4(\ell+h^{\vee})}
((\mathcal{N}_{g}-\alpha^{i})
(\mathcal{N}_{g}-\alpha^{i}-1)
a^{i'}, a^{i})w^{b}.
\end{align}

The first type of relations are the same as the corresponding 
type of relations 
in $\overarc{M}_{\ell}^{[g]}$ by Theorem \ref{ulbtm-va}. The second type of relations
(\ref{t-affine-L-0-relation}) correspond exactly 
to the relations (\ref{twisted-L-0-psi}) in $\overarc{M}_{\ell}^{[g]}$. 
The relations (\ref{twisted-L-0-psi}) are also the only relations 
in $\overarc{M}^{[g]}_{\ell}$ in addition to the
relations induced from the relations in $\widehat{M}^{[g]}_{\ell}$.  Thus the theorem is proved. 
\epfv

\begin{rema}\label{element-order-2}
{\rm We have seen in the proof of Theorem \ref{ulbtm-voa} that
$U(\hat{\g}^{[g]})
\otimes_{G(\Omega^{[g]}, \hat{\g}_{+}^{[g]}, \mathbf{k})}M$ is spanned by 
elements of the form  (\ref{element-form-0-5}). Using the commutator relations 
for $a^{i}(n)$ for $i\in I$ and $n\in \alpha^{i}+\Z$, we see that it is in fact spanned by elements of the form
\begin{equation}\label{element-form-0-6}
a^{i_{1}}(n_{1})
\cdots a^{i_{l}}(n_{l})a^{j_{1}}(\alpha^{j_{1}})
\cdots a^{j_{m}}(\alpha^{j_{m}})w^{b}
\end{equation}
for $i_{p}\in I$,
$n_{p}\in \alpha^{i_{p}}-\Z_{+}$ for $p=1, \dots, l$, $j_{q}\in I_{\I}$ for 
$q=1, \dots, m$ and
$b\in A$. By Theorem \ref{ulbtm-va}, we also see that $\overarc{M}_{\ell}^{[g]}$  
is spanned by elements of the form
\begin{equation}\label{element-form-0-7}
(a^{i_{1}}_{[g], \ell})_{n_{1},0}
\cdots (a^{i_{l}}_{[g], \ell})_{n_{l}, 0}
(a^{j_{1}}_{[g], \ell})_{\alpha^{j_{1}},0}
\cdots (a^{j_{m}}_{[g], \ell})_{\alpha^{j_{m}}, 0}
(\psi_{[g], h}^{b})_{-1, 0}\one
\end{equation}
for $i_{p}\in I$,
$n_{p}\in \alpha^{i_{p}}-\Z_{+}$ for $p=1, \dots, l$, $j_{q}\in I_{\I}$ for 
$q=1, \dots, m$ and
$b\in A$.}
\end{rema}

We now construct and identify explicitly grading-restricted generalized $g$-twisted modules
for $M(\ell, 0)$ viewed as a vertex operator algebra. 
We assume that  $M$ is in addition a finite-dimensional $\hat{\g}_{\mathbb{I}}$-module
with a compatible action of $g$
such that the action of $\Omega^{[g]}$, or equivalently, the operator 
$L_{M}(0)=\frac{1}{2(\ell+h^{\vee})}
\Omega^{[g]}$,  on $M$ is induced from this $\hat{\g}_{\mathbb{I}}$-module
structure.
In particular, $M$ is a direct sum of generalized eigenspaces of $L_{M}(0)$ as above. 
We have a universal lower-bounded generalized $g$-twisted $M(\ell, 0)$-module
$\overarc{M}^{[g]}_{\ell}$. 
As in the preceding subsection, since $\hat{\g}^{[g]}_{\I}$ is spanned by elements of the form 
$a^{i}(\alpha^{i})$ for $i\in I_{\I}$ and $\{w^{a}\}_{a\in 
A}$ is a basis of $M$, there exist $\lambda_{ic}^{a}\in \C$
for  $i\in I_{\I}$ and $b, c\in A$
such that 
$$a^{i}(\alpha^{i})w^{b}=\sum_{c\in A}\lambda_{ic}^{b}w^{c}$$
for $i\in I_{\I}$ and $b\in A$. 
Consider the lower-bounded generalized $g$-twisted $M(\ell, 0)$-submodule
of $\overarc{M}^{[g]}_{\ell}$ generated by elements of the form
\begin{align*}
(a^{i}_{[g], \ell})_{\alpha^{i}, 0}
(\psi^{b}_{[g]})_{-1, 0}\one
-\sum_{c\in A}\lambda_{ic}^{b}(\psi^{c}_{[g]})_{-1, 0}\one
\end{align*}
for $i\in I_{\I}$ and $b\in A$.
We denote the quotient of $\overarc{M}^{[g]}_{\ell}$ by this submodule 
by $\widetilde{M}^{[g]}_{\ell}$. Then $\widetilde{M}^{[g]}_{\ell}$ is also a 
lower-bounded generalized $g$-twisted $M(\ell, 0)$-module. 
We shall use the same notations for the generating twisted fields, generator
twist fields and their coefficients for $\overarc{M}^{[g]}_{\ell}$ to denote the corresponding 
fields and coefficients for $\widetilde{M}^{[g]}_{\ell}$. 

On the $\hat{\g}^{[g]}_{\I}$-module $M$,
we define 
actions of $\hat{\g}^{[g]}_{+}$ and $\mathbf{k}$ to be $0$ and $\ell$, respectively. 
Then we have the induced lower-bounded 
$\hat{\g}^{[g]}$-module $U(\hat{\g}^{[g]})
\otimes_{U(\hat{\g}^{[g]}_{+}\oplus \hat{\g}^{[g]}_{0})}M$.

\begin{thm}\label{first-main}
As a lower-bounded $\hat{\g}^{[g]}$-module, $\widetilde{M}^{[g]}_{\ell}$ is equivalent to 
the induced lower-bounded $\hat{\g}^{[g]}$-module $U(\hat{\g}^{[g]})
\otimes_{U(\hat{\g}^{[g]}_{+}\oplus \hat{\g}^{[g]}_{0})}M$.
\end{thm}
\pf
By Theorem \ref{ulbtm-voa}, we have an invertible  $\hat{\g}^{[g]}$-module map $\overarc{\rho}$ from 
$U(\hat{\g}^{[g]})
\otimes_{G(\Omega^{[g]}, \hat{\g}_{+}^{[g]}, \mathbf{k})}M$ to $\overarc{M}^{[g]}_{\ell}$
which maps  (\ref{element-form-0-6}) to 
 (\ref{element-form-0-7}). Since $\widetilde{M}^{[g]}_{\ell}$ is a quotient of 
$\overarc{M}^{[g]}_{\ell}$, we have a surjective $\hat{\g}^{[g]}$-module map $\tilde{\varrho}$  from 
$U(\hat{\g}^{[g]})
\otimes_{G(\Omega^{[g]}, \hat{\g}_{+}^{[g]}, \mathbf{k})}M$ to $\widetilde{M}^{[g]}_{\ell}$.

Note that by Poincar\'{e}-Birkhoff-Witt theorem,
$U(\hat{\g}^{[g]})
\otimes_{U(\hat{\g}^{[g]}_{+}\oplus \hat{\g}^{[g]}_{0})}M$ as a graded vector space is isomorphic to 
$U(\hat{\g}^{[g]}_{-})
\otimes M$.
In particular, the $\hat{\g}^{[g]}$-module $U(\hat{\g}^{[g]})
\otimes_{U(\hat{\g}^{[g]}_{+}\oplus \hat{\g}^{[g]}_{0})}M$ is spanned by elements of the form
\begin{equation}\label{element-form-0-8}
a^{i_{1}}(n_{1})\cdots a^{i_{l}}(n_{l})w^{b}
\end{equation}
for $i_{p}\in I$,
$n_{p}\in \alpha^{i_{p}}-\Z_{+}$ for $p=1, \dots, l$ and
$b\in A$.

By Remark \ref{element-order-2}, $\overarc{M}^{[g]}_{\ell}$
is spanned by elements of the form (\ref{element-form-0-7}) 
for $i_{p}\in I$,
$n_{p}\in \alpha^{i_{p}}-\Z_{+}$ for $p=1, \dots, l$, $j_{q}\in I_{\I}$ for
$q=1, \dots, m$ and
$b\in A$. Then $\widetilde{M}^{[g]}_{\ell}$ is spanned by elements of the form
\begin{equation}\label{element-form}
(a^{i_{1}}_{[g], \ell})_{n_{1},0}
\cdots (a^{i_{l}}_{[g], \ell})_{n_{l}, 0}
(\psi_{[g]}^{b})_{-1, 0}\one
\end{equation}
for $i_{p}\in I$,
$n_{p}\in \alpha^{i_{p}}-\Z_{+}$ for $p=1, \dots, l$ and
$b\in A$.

Since elements of $U(\hat{\g}^{[g]})
\otimes_{G(\Omega^{[g]}, \hat{\g}_{+}^{[g]}, \mathbf{k})}M$ of 
the same form as (\ref{element-form-0-8}) are sent under $\overarc{\rho}$ to 
elements of $\overarc{M}^{[g]}_{\ell}$ of the same form as (\ref{element-form}),
$\tilde{\varrho}$ maps elements of $U(\hat{\g}^{[g]})
\otimes_{G(\Omega^{[g]}, \hat{\g}_{+}^{[g]}, \mathbf{k})}M$ of 
the same form as (\ref{element-form-0-8})  to elements of $\widetilde{M}^{[g]}_{\ell}$ of the form 
(\ref{element-form}). But the only relations among 
elements of $U(\hat{\g}^{[g]})
\otimes_{G(\Omega^{[g]}, \hat{\g}_{+}^{[g]}, \mathbf{k})}M$ of
 the same form as (\ref{element-form-0-8}) are generated by 
the commutator relations for $a^{i}(n)$ for $i\in I$ and $n\in \alpha^{i}-\Z_{+}$.
Since $(a^{i}_{[g], \ell})_{n,0}$ for $i\in I$ and $n\in \alpha^{i}-\Z_{+}$ satisfy the 
same commutator relations as $a^{i}(n)$ and the only relations among 
elements of $\widetilde{M}^{[g]}_{\ell}$ of the form (\ref{element-form}) 
are generated by these commutator relations, the surjective $\hat{\g}^{[g]}$-module map $\tilde{\varrho}$
induces a bijective $\hat{\g}^{[g]}$-module map $\tilde{\rho}$ from $U(\hat{\g}^{[g]})
\otimes_{U(\hat{\g}^{[g]}_{+}\oplus \hat{\g}^{[g]}_{0})}M$ to $\widetilde{M}^{[g]}_{\ell}$.
Thus $U(\hat{\g}^{[g]})
\otimes_{U(\hat{\g}^{[g]}_{+}\oplus \hat{\g}^{[g]}_{0})}M$ is equivalent to $\widetilde{M}^{[g]}_{\ell}$.
\epfv

\begin{rema}
{\rm From the proof of Theorem \ref{first-main}, we see that 
$U(\hat{\g}^{[g]})
\otimes_{U(\hat{\g}^{[g]}_{+}\oplus \hat{\g}^{[g]}_{\I})}M$ and
$\widetilde{M}^{[g]}_{\ell}$ are spanned by elements of the form 
(\ref{element-form-0-8}) and (\ref{element-form}), respectively, 
for $i_{p}\in I$,
$n_{p}\in \alpha^{i_{p}}-\Z_{+}$ for $p=1, \dots, l$ and
$b\in A$. }
\end{rema}

\begin{thm}\label{grad-restr}
The lower-bounded generalized $g$-twisted $M(\ell, 0)$-module
$\widetilde{M}^{[g]}_{\ell}$ is in fact grading restricted. 
\end{thm}
\pf
By Theorem \ref{first-main}, we need only prove that 
$U(\hat{\g}^{[g]})
\otimes_{U(\hat{\g}^{[g]}_{+}\oplus \hat{\g}^{[g]}_{0})}M$ is grading restricted. 
But $U(\hat{\g}^{[g]})
\otimes_{U(\hat{\g}^{[g]}_{+}\oplus \hat{\g}^{[g]}_{0})}M$ is a graded subspace of 
$U(\hat{\g}^{[g]})
\otimes_{U(\hat{\g}^{[g]}_{+}\oplus \hat{\g}^{[g]}_{0})}\Lambda(M)$. By Theorem \ref{grad-rest},
$U(\hat{\g}^{[g]})
\otimes_{U(\hat{\g}^{[g]}_{+}\oplus \hat{\g}^{[g]}_{0})}\Lambda(M)$ is grading restricted, 
$U(\hat{\g}^{[g]})
\otimes_{U(\hat{\g}^{[g]}_{+}\oplus \hat{\g}^{[g]}_{0})}M$ is also grading restricted. 
\epfv

\subsection{Basic properties}\label{4.3}

The lower-bounded or grading-restricted $g$-twisted $M(\ell, 0)$-modules 
$\widehat{M}^{[g]}_{\ell}$, $\widebreve{M}^{[g]}_{\ell}$, 
$\overarc{M}^{[g]}_{\ell}$, $\widetilde{M}^{[g]}_{\ell}$ constructed above 
all have their own universal properties and other basic properties. 
We first give the universal properties of $\widehat{M}^{[g]}_{\ell}$ and $\overarc{M}^{[g]}_{\ell}$.

\begin{thm}\label{general-univ-va}
Let $(W, Y^{g}_{W})$ be a lower-bounded generalized $g$-twisted module
for $M(\ell, 0)$ viewed
as a grading-restricted vertex algebra (vertex operator algebra when $\g$ is simple and $\ell+h^{\vee}\ne 0$)
and $M^{0}$ a  subspace of $W$ invariant under the actions of 
$g$, $\mathcal{S}_{g}$, $\mathcal{N}_{g}$, $L_{W}(0)$, 
$L_{W}(0)_{S}$ and
$L_{W}(0)_{N}$.  Assume that $\hat{\g}^{[g]}_{+}$ acts on $M^{0}$ as $0$. 
If there is a linear map $f: M\to M^{0}$ 
 commuting with the actions of 
$g$, $\mathcal{S}_{g}$, $\mathcal{N}_{g}$, $L_{W}(0)|_{M^{0}}$ and $L_{M}(0)$, 
$L_{W}(0)_{S}|_{M^{0}}$ and $L_{M}(0)_{S}$ and
$L_{W}(0)_{N}|_{M^{0}}$ and $L_{M}(0)_{N}$, then there exists a unique module 
map $\hat{f}: \widehat{M}^{[g]}_{\ell}\to W$ $(\overarc{f}: \overarc{M}^{[g]}_{\ell}\to W)$
such that  $\hat{f}|_{M}=f$ $(\overarc{f}|_{M}=f)$. 
If $f$ is surjective and $(W, Y^{g}_{W})$ is generated by the 
coefficients of $(Y^{g})_{WM(\ell, 0)}^{W}(w, x)\one$ for $w\in M^{0}$, 
where $(Y^{g})_{WM(\ell, 0)}^{W}$ 
is the twist vertex operator map obtained from $Y_{W}^{g}$ (see \cite{H-twist-vo}), then 
$\hat{f}$ $(\overarc{f})$ is surjective. 
\end{thm}
\pf
Since $f$ commutes with the action of $L_{W}(0)|_{M^{0}}$ and $L_{\widehat{M}^{[g]}_{\ell}}(0)$, 
we have $f(M_{[h]})\subset M^{0}_{[h]}$ for $h\in Q_{M}$. 
Since $\hat{\g}^{[g]}_{+}$ acts on $M^{0}_{[h]}$ as $0$,
no nonzero elements of the submodule of $W$ generated by $M^{0}_{[h]}$ have weights less than $h$.
Then by the universal property of 
$\widehat{(M_{[h]})}_{\ell, h}^{[g]}$ (given by Theorem 5.2 in \cite{H-const-twisted-mod}), 
there exists a unique module map $\hat{f}_{h}: \widehat{(M_{[h]})}_{\ell, h}^{[g]}
\to W$ such that $\hat{f}_{h}|_{M_{[h]}}=f|_{M_{[h]}}$. Let 
$\hat{f}: \widehat{M}^{[g]}_{\ell}\to W$  be defined to be $\hat{f}_{h}$ on 
$\widehat{(M_{[h]})}_{\ell, h}^{[g]}$. Then $\hat{f}|_{M}=f$. The uniqueness of $\hat{f}$
follows from the uniqueness of $\hat{f}_{h}$ for $h\in Q_{M}$. The second conclusion also follows
from the property of $\hat{f}_{h}$ (see Theorem 5.2 in \cite{H-const-twisted-mod}) and the fact that 
$\widehat{(M_{[h]})}_{\ell, h}^{[g]}$ is generated by the subspace spanned by 
$(\psi^{b}_{[g]})_{-1, 0}$ (see Theorem 2.3 in \cite{H-exist-twisted-mod}). 

In the case that $M(\ell, 0)$ is viewed as a vertex operator algebra, we have a module map 
$\hat{f}: \widehat{M}^{[g]}_{\ell}\to W$. Since on $W$, 
$L_{W}(0)=\res_{x}xY_{W}^{g}(\omega_{M(\ell, 0)}, x)$ and
$L_{W}(-1)=\res_{x}Y_{W}^{g}(\omega_{M(\ell, 0)}, x)$ and $\overarc{f}|_{M}=f$ is obtained 
from $\widehat{M}^{[g]}_{\ell}$ by taking the quotient by a submodule 
generated by exactly these relations, we see that $\hat{f}$ induces a module map
$\overarc{f}: \overarc{M}^{[g]}_{\ell}\to W$. The other conclusions follow from the properties 
of $\hat{f}$ which we have proved.
\epfv

We have the following immediate consequence whose proof is omitted:

\begin{cor}\label{univ-cor-lb}
Let $W$ be a lower-bounded generalized $g$-twisted module for $M(\ell, 0)$ viewed
as a grading-restricted vertex algebra (vertex operator algebra  
when $\g$ is simple and $\ell+h^{\vee}\ne 0$)  generated by a 
subspace $M$ invariant under $g$, $\mathcal{L}_{g}$, $\mathcal{S}_{g}$,
$\mathcal{N}_{g}$, $L_{W}(0)$, $L_{W}(0)_{S}$ and $L_{W}(0)_{N}$ and annihilated
 by $\hat{\g}^{[g]}_{+}$. 
Then 
$W$ as a lower-bounded $\hat{\g}^{[g]}$-module is equivalent to a quotient of 
$U(\hat{\g}^{[g]})
\otimes _{U(\hat{\g}^{[g]}_{+}\oplus \C\mathbf{k})}\Lambda(M)$
$(U(\hat{\g}^{[g]})\otimes_{G(\Omega^{[g]}, \hat{\g}_{+}^{[g]}, \mathbf{k})}M)$.
Conversely, let  $M$ be a vector space with actions of  $g$, $\mathcal{L}_{g}$, $\mathcal{S}_{g}$,
$\mathcal{N}_{g}$, $L_{W}(0)$, $L_{W}(0)_{S}$ and $L_{W}(0)_{N}$ satisfying the conditions
given above. Then 
a quotient module of 
$U(\hat{\g}^{[g]})
\otimes _{U(\hat{\g}^{[g]}_{+}\oplus \C\mathbf{k})}\Lambda(M)$
$(U(\hat{\g}^{[g]})\otimes_{G(\Omega^{[g]}, \hat{\g}_{+}^{[g]}, \mathbf{k})}M)$ has a natural structure 
of a lower-bounded generalized $g$-twisted $M(\ell, 0)$-module when $M(\ell, 0)$ is viewed
as a grading-restricted vertex algebra (vertex operator algebra  when $\g$ is simple and $\ell+h^{\vee}\ne 0$). 
\end{cor}

Now we discuss the universal properties of $\widebreve{M}^{[g]}_{\ell}$ and 
$\widetilde{M}^{[g]}_{\ell}$. 

\begin{thm}\label{univ-property}
Let $(W, Y^{g}_{W})$ be a lower-bounded generalized $g$-twisted module
for $M(\ell, 0)$ viewed as a grading-restricted vertex algebra (vertex operator algebra
 when $\g$ is simple and $\ell+h^{\vee}\ne 0$).
Let $M^{0}$ a finite-dimensional $\hat{\g}^{[g]}_{\I}$-submodule of $W$ invariant also under the actions of 
$g$, $\mathcal{S}_{g}$, $\mathcal{N}_{g}$, $L_{W}(0)$, $L_{W}(0)_{S}$ and
$L_{W}(0)_{N}$ and annihilated by $\hat{\g}^{[g]}_{+}$. 
Assume that there is a $\hat{\g}^{[g]}_{\I}$-module map $f: M\to M^{0}$ 
commuting with the actions of 
$g$, $\mathcal{S}_{g}$, $\mathcal{N}_{g}$, $L_{W}(0)|_{M^{0}}$ and $L_{M}(0)$, 
$L_{W}(0)_{S}|_{M^{0}}$ and $L_{M}(0)_{S}$ and
$L_{W}(0)_{N}|_{M^{0}}$ and $L_{M}(0)_{N}$. Then there exists a unique module 
map $\breve{f}: \widebreve{M}^{[g]}_{\ell}\to W$ $(\tilde{f}: \widetilde{M}^{[g]}_{\ell}\to W)$
such that $\breve{f}|_{M}=f$ $(\tilde{f}|_{M}=f)$. 
If $f$ is surjective and $(W, Y^{g}_{W})$ is generated by the 
coefficients of $(Y^{g})_{WM(\ell, 0)}^{W}(w, x)\one$ for $w\in M^{0}$, 
where $(Y^{g})_{WM(\ell, 0)}^{W}$ 
is the twist vertex operator map obtained from $Y_{W}^{g}$, then 
$\breve{f}$ $(\tilde{f})$ is surjective and thus $W$ is grading restricted. 
\end{thm}
\pf
By Theorem \ref{general-univ-va}, we have a unique module map $\hat{f}: \widehat{M}^{[g]}_{\ell}
\to W$ such that $\hat{f}|_{M}=f$. Since $f$ is a $\hat{\g}^{[g]}_{\I}$-module map,
we have 
$$a_{W}^{i}(\alpha^{i})f(w^{b})-\sum_{c\in A}\lambda_{ic}^{b}f(w^{c})
=f\left(a^{i}(\alpha^{i})w^{b}-\sum_{c\in A}\lambda_{ic}^{b}w^{c}\right)=f(0)=0$$
for $i\in I_{\I}$, $b\in A$. Since $\hat{f}$ is a module map, we have
$$\hat{f}((\psi^{b}_{[g], h^{b}})_{-1, 0}\one)=
\hat{f}(((Y^{g})_{\widehat{M}^{[g]}_{\ell}M(\ell, 0)}^{\widehat{M}^{[g]}_{\ell}})_{-1, 0}
(w^{b})\one)=((Y^{g})_{WM(\ell, 0)}^{W})_{-1, 0}(f(w^{b}))\one
=f(w^{b})$$ 
for $b\in A$. 
Thus we obtain 
$$\hat{f}\left((a^{i}_{[g], \ell})_{\alpha^{i}, 0}
(\psi^{b}_{[g], h^{b}})_{-1, 0}\one
-\sum_{c\in A}\lambda_{ic}^{b}(\psi^{c}_{[g], h^{c}})_{-1, 0}\one\right)
=a_{W}^{i}(\alpha^{i})f(w^{b})-\sum_{c\in A}\lambda_{ic}^{b}f(w^{c})=0.$$
So 
$$(a^{i}_{[g], \ell})_{\alpha^{i}, 0}
(\psi^{b}_{[g], h^{b}})_{-1, 0}\one
-\sum_{c\in A}\lambda_{ic}^{b}(\psi^{c}_{[g], h^{c}})_{-1, 0}\one$$
is in the kernel of $\hat{f}$. In particular, we have a module map
$\breve{f}: \widebreve{M}^{[g]}_{\ell, h}\to W$. The uniqueness of 
$\breve{f}$ and the  surjectivity of $\breve{f}$ when $f$ is surjective 
 follow from the uniqueness of $\hat{f}$ and  the  surjectivity of $\hat{f}$.
Since $\widebreve{M}^{[g]}_{\ell}$ is grading restricted, 
$W$ is grading restricted when $\breve{f}$ is surjective. 

The proof for  $\widetilde{M}^{[g]}_{\ell}$ is the same except that we
use $\overarc{M}^{[g]}_{\ell}$ instead of $\widehat{M}^{[g]}_{\ell}$. 
\epfv

We also have the following immediate consequence whose proof is also omitted:

\begin{cor}\label{univ-cor-gr}
Let $W$ be a lower-bounded generalized $g$-twisted module for $M(\ell, 0)$ viewed
as a grading-restricted vertex algebra (vertex operator algebra
 when $\g$ is simple and $\ell+h^{\vee}\ne 0$) generated by a finite-dimensional 
$\hat{\g}_{I_{\I}}$-submodule $M$ invariant under $g$, $\mathcal{L}_{g}$, $\mathcal{S}_{g}$,
$\mathcal{N}_{g}$, $L_{W}(0)$, $L_{W}(0)_{S}$ and $L_{W}(0)_{N}$ and  annihilated by
$\hat{\g}^{[g]}_{+}$. 
Then 
$W$ as a lower-bounded $\hat{\g}^{[g]}$-module is equivalent to a quotient of 
$U(\hat{\g}^{[g]})
\otimes _{U(\hat{\g}^{[g]}_{+}\oplus \hat{\g}^{[g]}_{0})}\Lambda(M)$
$(U(\hat{\g}^{[g]})
\otimes _{U(\hat{\g}^{[g]}_{+}\oplus \hat{\g}^{[g]}_{0})}M)$ and, in particular, 
$W$ is grading restricted.
Conversely, let  $M$ be a finite-dimensional 
$\hat{\g}_{I_{\I}}$-module with compatible actions of  $g$, $\mathcal{L}_{g}$, $\mathcal{S}_{g}$,
$\mathcal{N}_{g}$, $L_{W}(0)$, $L_{W}(0)_{N}$ and $L_{W}(0)_{N}$ satisfying the conditions
we discussed above. Then 
a quotient module of 
$U(\hat{\g}^{[g]})
\otimes _{U(\hat{\g}^{[g]}_{+}\oplus \hat{\g}^{[g]}_{0})}\Lambda(M)$
$(U(\hat{\g}^{[g]})
\otimes _{U(\hat{\g}^{[g]}_{+}\oplus \hat{\g}^{[g]}_{0})}M)$ has a natural structure 
of a grading-restricted generalized $g$-twisted $M(\ell, 0)$-module when $M(\ell, 0)$ is viewed
as a grading-restricted vertex algebra (vertex operator algebra  when $\g$ is simple and $\ell+h^{\vee}\ne 0$). 
\end{cor}

\begin{rema}
{\rm In \cite{H-exist-twisted-mod}, only  the existence of a grading-restricted  generalized
$g$-twisted $V$-module for a grading-restricted vertex algebra $V$ and an
automorphism $g$ of $V$ is proved under suitable conditions. 
But no construction is given in that paper. Theorems \ref{grad-rest}, \ref{grad-restr} and Corollary 
\ref{univ-cor-gr} give explicit constructions of grading-restricted  generalized
$g$-twisted $M(\ell, 0)$-modules.}
\end{rema}

\renewcommand{\theequation}{\thesection.\arabic{equation}}
\renewcommand{\thethm}{\thesection.\arabic{thm}}
\setcounter{equation}{0}
\setcounter{thm}{0}
\section{Lower-bounded and grading-restricted generalized \\
twisted $L(\ell, 0)$-modules}

In this section, we construct lower-bounded and grading-restricted generalized $g$-twisted $L(\ell, 0)$-modules.
We shall mainly discuss the case  that $\g$ is simple, $\ell\in \Z_{+}$ and 
$L(\ell, 0)$ is viewed as a vertex operator algebra.

We first give some straightforward general results. 

\begin{prop}\label{L-mod-t-aff-mod}
Let $W$ be a lower-bounded generalized $g$-twisted module for
$L(\ell, 0)$ viewed as a grading-restricted vertex algebra (vertex operator algebra 
 when $\g$ is simple and $\ell+h^{\vee}\ne 0$)
generated by a 
subspace $M$ invariant under $g$, $\mathcal{L}_{g}$, $\mathcal{S}_{g}$,
$\mathcal{N}_{g}$, $L_{W}(0)$, $L_{W}(0)_{N}$ and $L_{W}(0)_{N}$ and annihilated by
$\hat{\g}^{[g]}_{+}$.  Then $W$ 
is a lower-bounded generalized $g$-twisted module for $M(\ell, 0)$ viewed as 
a grading-restricted vertex algebra (vertex operator algebra). In particular, $W$
is a lower bounded $\hat{\g}^{[g]}$-module and is a quotient of 
$U(\hat{\g}^{[g]})
\otimes _{U(\hat{\g}^{[g]}_{+}\oplus \C\mathbf{k})}\Lambda(M)$
$(U(\hat{\g}^{[g]})\otimes_{G(\Omega^{[g]}, \hat{\g}_{+}^{[g]}, \mathbf{k})}M)$.
If  $M$ is in addition a finite-dimensional $\hat{\g}_{\I}$-module, then $W$
 is a quotient of 
$U(\hat{\g}^{[g]})
\otimes _{U(\hat{\g}^{[g]}_{+}\oplus \hat{\g}^{[g]}_{0})}\Lambda(M)$
$(U(\hat{\g}^{[g]})
\otimes _{U(\hat{\g}^{[g]}_{+}\oplus \hat{\g}^{[g]}_{0})}M)$ and in particular, 
$W$ is grading restricted.
\end{prop}
\pf
Since $L(\ell, 0)$ is a quotient of $M(\ell, 0)$, 
$W$ must be a lower-bounded generalized $g$-twisted $M(\ell, 0)$-module.
By Proposition \ref{twisted-affine-mod}, $W$ is a lower-bounded $\hat{\g}^{[g]}$-module.
By Corollary \ref{univ-cor-lb}, $W$ as a lower-bounded $\hat{\g}^{[g]}$-module is a quotient of 
$U(\hat{\g}^{[g]})
\otimes _{U(\hat{\g}^{[g]}_{+}\oplus \C\mathbf{k})}\Lambda(M)$
($U(\hat{\g}^{[g]})\otimes_{G(\Omega^{[g]}, \hat{\g}_{+}^{[g]}, \mathbf{k})}M$).

If  $M$ is in addition a finite-dimensional $\hat{\g}_{\I}$-module, then by Corollary \ref{univ-cor-gr},
$W$ as a lower-bounded $\hat{\g}^{[g]}$-module is a quotient of 
$U(\hat{\g}^{[g]})
\otimes _{U(\hat{\g}^{[g]}_{+}\oplus \hat{\g}^{[g]}_{0})}\Lambda(M)$
($U(\hat{\g}^{[g]})
\otimes _{U(\hat{\g}^{[g]}_{+}\oplus \hat{\g}^{[g]}_{0})}M$) and, in particular, 
$W$ is grading restricted.
\epfv

%We say that a lower-bounded generalized 
%$g$-twisted $M(\ell, 0)$-module $(W, Y_{W}^{g})$ is faithful if $Y_{W}^{g}(u, x)=0$ 
%implies $u=0$. 
Using the structure of $L(\ell, 0)$ as a $\hat{\g}$-module and 
properties of  lower-bounded or grading-restricted generalized $g$-twisted 
$M(\ell, 0)$-modules, we have the following result:

\begin{prop}\label{L-mod-t-aff-mod-1}
Assume that $\g$ is simple with a given Cartan subalgebra and 
a given set of simple roots and $\ell\in \Z_{+}$. Let $(W, Y_{W}^{g})$
be a  lower-bounded (grading-restricted) generalized $g$-twisted 
$M(\ell, 0)$-module% such that $W$ as a $\hat{\g}^{[g]}$-module is faithful
. 
Then $(W, Y_{W}^{g})$
 is a  lower-bounded (grading-restricted) generalized $g$-twisted 
$L(\ell, 0)$-module 
if and only if  $Y^{g}_{W}(e_{\theta}(-1)^{\ell+1}\one, x)=0$ , where $\theta$ is the 
highest root of $\g$ and $e_{\theta}\in \g_{\theta}\setminus \{0\}$.
\end{prop}
\pf
From  \cite{K} and Proposition 6.6.17 in \cite{LL}, we know that 
$L(\ell, 0)=M(\ell, 0)/I(\ell, 0)$ where $I(\ell, 0)=U(\hat{\g})e_{\theta}(-1)^{\ell+1}\one$. 
Then $W$  is a lower-bounded (grading-restricted) generalized $g$-twisted 
$L(\ell, 0)$-module 
if and only if $Y^{g}_{W}(I(\ell, 0), x)=0$. 
We need only prove that $Y^{g}_{W}(I(\ell, 0), x)=0$ if and only if $Y^{g}_{W}(e_{\theta}(-1)^{\ell+1}\one, x)=0$. 

If  $Y^{g}_{W}(I(\ell, 0), x)=0$, then certainly
$Y^{g}_{W}(e_{\theta}(-1)^{\ell+1}\one, x)=0$.
If 
$Y^{g}_{W}(e_{\theta}(-1)^{\ell+1}\one, x)=0$, then $Y^{g}_{W}(e_{\theta}(-1)^{\ell+1}\one, x)w=0$
for $w\in W$. Therefore we have $(Y^{g})_{WM(\ell, 0)}^{W}(w, x)e_{\theta}(-1)^{\ell+1}\one=0$ 
for $w\in W$, where
$(Y^{g})_{WM(\ell, 0)}^{W}$ is the twist vertex operator map introduced and studied in 
\cite{H-twist-vo}. But by Corollary 4.3 in \cite{H-twist-vo} (commutativity for twisted and twist vertex operators), 
for $v\in M(\ell, 0)$ and $w'\in W'$, 
\begin{align*}
F^{p}&(\langle w', (Y^{g})_{WM(\ell, 0)}^{W}(w, z_{2})
Y_{M(\ell, 0)}(v, z_{1})e_{\theta}(-1)^{\ell+1}\one\rangle)\nn
&=F^{p}(\langle w', Y_{W}^{g}(v, z_{1})(Y^{g})_{WM(\ell, 0)}^{W}(w, z_{2})
e_{\theta}(-1)^{\ell+1}\one\rangle)\nn
&=0,
\end{align*}
where the first two lines are different expressions of
the $p$-th branch of a multivalued analytic function that 
$$\langle w', Y_{W}^{g}(v, z_{1})(Y^{g})_{WM(\ell, 0)}^{W}(w, z_{2})e_{\theta}(-1)^{\ell+1}\one\rangle$$
and 
$$\langle w', (Y^{g})_{WM(\ell, 0)}^{W}(w, z_{2})Y_{M(\ell, 0)}(v, z_{1})e_{\theta}(-1)^{\ell+1}\one\rangle$$
converge to in the regions $|z_{1}|>|z_{2}|>0$ and $|z_{2}|>|z_{1}|>0$, respectively. 
Thus we have
$$(Y^{g})_{WM(\ell, 0)}^{W}(w, x)Y_{M(\ell, 0)}(v, x_{1})e_{\theta}(-1)^{\ell+1}\one=0$$
for $w\in W$. From the definition of the twist vertex operator map $(Y^{g})_{WV}^{W}$,
we obtain
$$Y^{g}_{W}(Y_{M(\ell, 0)}(v, x_{1})e_{\theta}(-1)^{\ell+1}\one, x)=0.$$
But the coefficients of $Y_{M(\ell, 0)}(v, x_{1})$ span $U(\hat{\g})$. So we obtain
$Y^{g}_{W}(I(\ell, 0), x)=0$. 
\epfv

Our goal is to give and identify explicitly universal lower-bounded and grading-restricted generalized $g$-twisted 
$L(\ell, 0)$-modules using the results we obtained in the preceding section and some conditions on 
the corresponding $\hat{\g}^{[g]}$-modules.
To do this, we first prove some results on $\mathcal{S}_{g}$ and $\mathcal{N}_{g}$. 

For $\mathcal{S}_{g}$, or equivalently, the semisimple automorphism 
$\sigma=e^{2\pi i\mathcal{S}_{g}}$, we have the 
following generalization of Proposition 8.1 in 
 \cite{K} on automorphisms of finite orders of a simple Lie algebra:

\begin{prop}\label{autom-decomp}
Assume that $\g$ is simple  with a given Cartan subalgebra $\mathfrak{h}$
and 
a given set $\Delta$ of simple roots. Let $\sigma=e^{2\pi i\mathcal{S}_{g}}$ be the semisimple
part of $g$. 
Then there exists an automorphism $\tau_{\sigma}$ of $\g$ such that 
$\sigma=\tau_{\sigma} e^{{\rm ad}\;h} \mu \tau_{\sigma}^{-1}$,
where $\mu$ is a diagram automorphism of $\g$ preserving $\mathfrak{h}$ 
and $\Delta$ and $h$ is an element  of the fixed-point subspace $\mathfrak{h}^{[0]}$ 
of  $\mathfrak{h}$.
\end{prop}
\pf
Let $\tilde{\mathfrak{h}}^{[0]}\subset \g^{[0]}$ be a maximal toral subalgebra (that is,
maximal ad-diagonalizable subalgebra) of $\g^{[0]}$ and  
$C_{\g}(\tilde{\mathfrak{h}}^{[0]})$  the centralizer of $\tilde{\mathfrak{h}}^{[0]}$ in $\g$. 
Extend $\tilde{\mathfrak{h}}^{[0]}$ to a maximal toral subalgebra $\tilde{\mathfrak{h}}$ of $\g$. Then 
$\tilde{\mathfrak{h}}$ is a Cartan subalgebra of $\g$ and
$C_{\g}(\tilde{\mathfrak{h}}^{[0]})=\tilde{\mathfrak{h}}+\sum_{\xi}\tilde{\g}_{\xi}$, where the sum is over the 
roots $\xi$ in $\tilde{\mathfrak{h}}$
such that $\xi$ restricted to $\tilde{\mathfrak{h}}^{[0]}$ is $0$ and $\tilde{\g}_{\xi}$ is the 
root space associated to $\xi$. We first prove that $C_{\g}(\tilde{\mathfrak{h}}^{[0]})=\tilde{\mathfrak{h}}$.

Let $\mathfrak{s}=\sum_{\xi}\tilde{\g}_{\xi}$. By definition,
$\mathfrak{s}$ is a subalgebra of $\g$ invariant under $g$ such that $\mathfrak{s}\cap \g^{[0]}=0$.
Moreover, the restriction of the bilinear form $(\cdot, \cdot)$ to $\mathfrak{s}$ is 
nondegenerate because $\tilde{\g}_{\xi_{1}}$ and $\tilde{\g}_{\xi_{2}}$ are orthogonal if $\xi_{1}+\xi_{2}\ne 0$
and the restriction of $(\cdot, \cdot)$ to $\tilde{\g}_{\xi}\times \tilde{\g}_{-\xi}$ is nondegenerate. 
Since $\mathfrak{s}$ is invariant under $g$, we have $\mathfrak{s}=\coprod_{\alpha\in P_{\mathfrak{s}}}
\mathfrak{s}^{[\alpha]}$
where $P_{\mathfrak{s}}$ is the set of $\alpha\in [0, 1)+i\R$ such that 
$e^{2\pi i\alpha}$ is an eigenvalue of $\sigma$ (or equivalently, of $g$) and 
for $\alpha\in \mathfrak{s}$, $\mathfrak{s}^{[\alpha]}=\mathfrak{s}\cap  \g^{[\alpha]}$ is the 
eigenspace of $\sigma$ (or equivalently, the generalized eigenspace of $g$)
in $\mathfrak{s}$ with the eigenvalue $e^{2\pi i \alpha}$.
For $\alpha\in ([0, 1)+i\R)\setminus P_{\mathfrak{s}}$, let $\mathfrak{s}^{[\alpha]}=0$. Then we have 
$\mathfrak{s}=\coprod_{\alpha\in [0, 1)+i\R}\mathfrak{s}^{[\alpha]}$. Moreover, by Lemma 
\ref{brac-eigenvalue} and the fact that $\mathfrak{s}$ is a subalgebra of $\g$, we have 
$[\mathfrak{s}^{[\alpha]}, \mathfrak{s}^{[\beta]}]\subset \mathfrak{s}^{[s(\alpha, \beta)]}$
for $\alpha, \beta\in [0, 1)+i\R$ (recall $s(\alpha, \beta)$ defined before Lemma \ref{brac-eigenvalue}). 
We need only prove that $\mathfrak{s}=0$.

Since $\g$ is finite dimensional, $\mathfrak{s}$ is finite dimensional and hence $P_{\mathfrak{s}}$ is a finite 
set. We use induction on the finitely many real parts of the elements of $P_{\mathfrak{s}}$.
 First, we know 
that $\mathfrak{s}^{[0]}=\mathfrak{s}\cap \g^{[0]}=0$. For $\alpha\in (P_{\mathfrak{s}}\setminus \{0\})\cap (i\R)$ 
and $a\in \mathfrak{s}^{[\alpha]}$, we know that
$({\rm ad}\;a)^{r}\mathfrak{s}^{[\beta]}\in \mathfrak{s}^{[r\alpha+\beta]}$
for $\beta\in P_{g}$. Since $P_{\mathfrak{s}}$ is a finite set, $\mathfrak{s}^{[r\alpha+\beta]}$ must be $0$ 
when $r$ is sufficiently large. So ${\rm ad}\;a$ is nilpotent on $\g$. Applying Lemma \ref{S-g-ortho}
to $\mathfrak{s}$ and the restriction of $\sigma$ to $\mathfrak{s}$,
we see that $(\cdot, \cdot)$ restricted to $\mathfrak{s}^{[\alpha]}\times \mathfrak{s}^{[-\alpha]}$ 
is nondegenerate. In particular, if $\mathfrak{s}^{[\alpha]}\ne 0$, then $\mathfrak{s}^{[-\alpha]}\ne 0$. 
If $\mathfrak{s}^{[\alpha]}\ne 0$, let $a\in \mathfrak{s}^{[\alpha]}\setminus \{0\}$
and $b\in \mathfrak{s}^{[-\alpha]}\setminus \{0\}$. Then both ${\rm ad}\;a$ and ${\rm ad}\;b$ are nilpotent 
on  $\mathfrak{s}$. Therefore the eigenvalues of ${\rm ad}\;a$ and ${\rm ad}\;b$ are all $0$. 
Since $[\mathfrak{s}^{[\alpha]}, \mathfrak{s}^{[-\alpha]}]\subset \mathfrak{s}^{[0]}=0$,
${\rm ad}\;a$ and ${\rm ad}\;b$ commute and hence can be diagonalized simultaneously. 
In particular, the trace of $({\rm ad}\;a)({\rm ad}\;b)$ is $0$. But this contradicts  the nondegneracy 
of  $(\cdot, \cdot)$ restricted to $\mathfrak{s}^{[\alpha]}\times \mathfrak{s}^{[-\alpha]}$ because 
$(a, b)$ is proportional to this trace. Thus $\mathfrak{s}^{[\alpha]}= 0$.

Note that for $\alpha\in P_{\mathfrak{s}}$ such that $\Re(\alpha)>0$,
if $\mathfrak{s}^{[\alpha]}\ne 0$, then $\mathfrak{s}^{[1-\alpha]}\ne 0$ since by Lemma 2.3,
the restriction of $(\cdot, \cdot)$ to $\mathfrak{s}$
is nondegenerate and $\mathfrak{s}^{[\alpha]}$ is orthogonal to $\mathfrak{s}^{[\beta]}$
for $\beta\ne 1-\alpha$. 
We now assume that for $\alpha\in P_{\mathfrak{s}}$ with $\Re(\alpha)>0$,
$\mathfrak{s}^{[\alpha']}=0$ for $\alpha'\in P_{g}$ and $\Re(\alpha')<\Re(\alpha)$ and 
for $\alpha'\in P_{g}$ and $\Re(\alpha')>\Re(1-\alpha)$. Then for $a\in \mathfrak{s}^{[\alpha]}$,
$({\rm ad}\;a)^{r}\mathfrak{s}^{[\beta]}\in \mathfrak{s}^{[s(\alpha, r, \beta)]}$,
where $s(\alpha, r, \beta)\equiv r\alpha+\beta \mod \Z$ satisfying 
$0\le \Re(s(\alpha, r, \beta))<1$. Since $0\le \beta<1$, there exists $r\in \Z_{+}$ such that 
$\Re((r-1)\alpha+\beta)<1$ but $\Re(r\alpha+\beta)\ge 1$. From $\Re((r-1)\alpha+\beta)<1$,
we obtain $0\le s(\alpha, r, \beta)=\Re(r\alpha+\beta)-1<\Re(\alpha)$.
By the  induction assumption, $\mathfrak{s}^{[s(\alpha, r, \beta)]}=0$. 
So we obtain $({\rm ad}\;a)^{r}\mathfrak{s}^{[\beta]}=0$ for $\beta\in P_{\mathfrak{s}}$. Thus ${\rm ad}\;a$ 
is nilpotent on $\mathfrak{s}$. Similarly, for $b\in  \mathfrak{s}^{[1-\alpha]}$,
$({\rm ad}\;b)^{r}\mathfrak{s}^{[\beta]}\in \mathfrak{s}^{[s(1-\alpha, r, \beta)]}$,
where $s(1-\alpha, r, \beta)\equiv r(1-\alpha)+\beta \mod \Z$ satisfying 
$0\le \Re(s(1-\alpha, r, \beta))<1$. When $\Re(\beta)=0$, since we have 
proved $\mathfrak{s}^{[\beta]}=0$, $({\rm ad}\;b)^{r}\mathfrak{s}^{[\beta]}=0$ for $r\in \Z_{+}$. 
When $\Re(\beta)\ne 0$, there exists $r\in \Z_{+}$ such that 
$\Re(-r\alpha+\beta)\ge -1$ but $\Re(-(r+1)\alpha+\beta)< -1$. Then 
$\Re(r(1-\alpha)+\beta)\ge r-1$ but $\Re((r-1)(1-\alpha)+\beta)< r$. Since $0<\Re(\beta)<1$, we obtain
$0< \Re(1-\alpha) < \Re(1-\alpha+\beta)=\Re(s(1-\alpha, r, \beta))<1$. By the induction assumption,
$\mathfrak{s}^{[s(1-\alpha, r, \beta)]}=0$. So we obtain $({\rm ad}\;b)^{r}\mathfrak{s}^{[\beta]}=0$.
Thus ${\rm ad}\;b$ is also nilpotent on $\mathfrak{s}$. 
If $\mathfrak{s}^{[\alpha]}\ne 0$, let $a\in \mathfrak{s}^{[\alpha]}\setminus \{0\}$
and $b\in \mathfrak{s}^{[1-\alpha]}\setminus \{0\}$. Then we have proved that 
both ${\rm ad}\;a$ and ${\rm ad}\;b$ are nilpotent 
on  $\mathfrak{s}$. Therefore the eigenvalues of ${\rm ad}\;a$ and ${\rm ad}\;b$ are all $0$. 
Since $[\mathfrak{s}^{[\alpha]}, \mathfrak{s}^{[1-\alpha]}]\subset \mathfrak{s}^{[0]}=0$,
${\rm ad}\;a$ and ${\rm ad}\;b$ commute and hence can be diagonalized simultaneously. 
In particular, the trace of $({\rm ad}\;a)({\rm ad}\;b)$ is $0$. Contradiction. Thus 
 $\mathfrak{s}^{[\alpha]}= 0$. This proves $\mathfrak{s}=0$.  

We have proved that $C_{\g}(\tilde{\mathfrak{h}}^{[0]})=\tilde{\mathfrak{h}}$. Now choose
$a\in \tilde{\mathfrak{h}}^{[0]}$ such that the centralizer $C_{\g}(a)$ of $a$ in $\g$
 is minimal among the collection of all centralizer $C_{\g}(b)$ of $b$ in $\g$
for $b\in \tilde{\mathfrak{h}}^{[0]}$. Note that since elements of $\tilde{\mathfrak{h}}^{[0]}$ are all 
semisimple or ad-diagonalizable, $C_{\g}(b)$ for $b\in \tilde{\mathfrak{h}}^{[0]}$ is equal to the space of
all elements of $\g$ on which ${\rm ad}\;b$ acts  nilpotently. Since $\tilde{\mathfrak{h}}^{[0]}\subset C_{\g}(a)$,
by Lemma A in Subsection 15.2 in \cite{Hum}, we have $C_{\g}(a)\subset C_{\g}(b)$
for $b\in \tilde{\mathfrak{h}}^{[0]}$. But $C_{\g}(\tilde{\mathfrak{h}}^{[0]})=\cap_{b\in \tilde{\mathfrak{h}}^{[0]}}
C_{\g}(b)$. So $C_{\g}(a)\subset C_{\g}(\tilde{\mathfrak{h}}^{[0]})=\tilde{\mathfrak{h}}$. 
But $\tilde{\mathfrak{h}}\subset C_{\g}(a)$. So we must have $\tilde{\mathfrak{h}}= C_{\g}(a)$,
that is, $a\in \tilde{\mathfrak{h}}^{[0]}$ is a regular semisimple element. As the centralizer of 
a fixed point of $g$, $\tilde{\mathfrak{h}}$ is a Cartan subalgebra of $\g$ invariant under $\sigma$. 
In particular, we have a root system $\widetilde{\Phi}$ obtained from $\tilde{\mathfrak{h}}$. 

The regular semisimple element $a$ cannot be orthogonal to any root $\xi$. Otherwise 
$[a, e_{\xi}]=(\xi, a)e_{\xi}=0$ for $e_{\xi}\in \g_{\xi}$ so that $e_{\xi}\in C_{\g}(a)=\tilde{\mathfrak{h}}$,
which is impossible. Thus if we let $\widetilde{\Phi}^{+}=\{\xi\in \widetilde{\Phi}\;|\;(\xi, a)>0\}$, then 
$\widetilde{\Phi}=\widetilde{\Phi}^{+}-\widetilde{\Phi}^{+}$. By Theorem$'$ in Subsection 10.1 in \cite{Hum}, 
the set $\widetilde{\Delta}$ of all indecomposable roots in $\widetilde{\Phi}^{+}$ is a set 
of simple roots of $\widetilde{\Phi}$
and $\widetilde{\Phi}^{+}$ is the set of positive roots.  Since $a$ is fixed by $\sigma$, 
$\sigma$ induces an automorphism of $\widetilde{\Phi}^{+}$. Choose $e_{\xi}\in \g_{\xi}\setminus \{0\}$ for 
$\xi\in \widetilde{\Delta}$. 
Let $\tilde{\mu}$ be the diagram automorphism of $\g$
corresponding to this automorphism of $\widetilde{\Phi}^{+}$ and $e_{\xi}$ for $\xi\in \widetilde{\Delta}$. 
Then $\sigma \tilde{\mu}^{-1}$ fix every element of $\tilde{\mathfrak{h}}$. 
In particular, $\sigma \tilde{\mu}^{-1}$ commutes with ${\rm ad}\;\tilde{a}$ for all $\tilde{a}\in \tilde{\mathfrak{h}}$
and thus can be diagonalized simultaneously together with ${\rm ad}\;\tilde{a}$. 
The root space decomposition $\g=\tilde{\mathfrak{h}}\oplus \coprod_{\xi\in \widetilde{\Phi}}
\tilde{\g}_{\xi}$ gives a diagonalization of ${\rm ad}\;\tilde{a}$. Also $\tilde{\g}_{\xi}$ for $\xi\in \widetilde{\Phi}$
are all one dimensional and $\sigma \tilde{\mu}^{-1}$  acts as the identity on $\tilde{\mathfrak{h}}$,
we see the this root space decomposition also give a diagonalization of $\sigma \tilde{\mu}^{-1}$. 
So on each root space $\tilde{\g}_{\xi}$, it must act as a scalar multiplication by 
$\lambda_{\xi}\in \C^{\times}$. Let $l_{0}(\lambda_{\xi})=\log |\lambda_{\xi}|
+i\arg \lambda_{\xi}$, where $0\le \arg \lambda_{\xi}<2\pi$. 
Let $\tilde{h}\in \tilde{\mathfrak{h}}$
be defined by $(\xi, \tilde{h})=\frac{1}{2\pi i}l_{0}(\lambda_{\xi})$ for $\xi\in \widetilde{\Delta}$.
Then $\lambda_{\xi}e_{\xi}=e^{2\pi i(\xi, \tilde{h})}e_{\xi}=e^{2\pi i({\rm ad}\;\tilde{h})}e_{\xi}$.
 Thus we obtain 
$\sigma\tilde{\mu}^{-1}=e^{2 \pi i({\rm ad}\;\tilde{h})}$, or equivalently, 
$\sigma= e^{2 \pi i({\rm ad}\;\tilde{h})}\tilde{\mu}$. Since 
$\sigma\tilde{\mu}^{-1}$ fix every element of $\tilde{\mathfrak{h}}$, $\tilde{h}$ must be in $\tilde{\mathfrak{h}}$.

Since any two Cartan subalgebras are conjugate to each other, there exists an automorphism $\nu$ of $\g$ 
such that $\nu(\mathfrak{h})=\tilde{\mathfrak{h}}$ and
$\nu(\Delta)=\widetilde{\Delta}$. Let $\mu =\nu \tilde{\mu} \nu^{-1}$ and 
$\breve{h}=\nu^{-1}(\tilde{h})\in \mathfrak{h}$. Then it is clear that $\mu$ is a diagram automorphism of 
$\g$ preserving $\mathfrak{h}$ and $\Delta$ and we have
$$\sigma=  e^{2 \pi i({\rm ad}\; \nu(\breve{h}))} \nu\mu \nu^{-1}
=\nu e^{2 \pi i({\rm ad}\;\breve{h})} \mu  \nu^{-1}.$$
But $\breve{h}$ might not be fixed by $\mu$. We need to find another automorphism such that after 
the conjugation by this automorphism, we have $h\in \mathfrak{h}^{[0]}$ (that is, fixed by $\mu$). 
This argument was in fact given by the proof of Lemma 8.3 in \cite{EMS}: Let $r$ be the order of $\mu$ (in fact 
$r=1, 2$ or $3$), $h=\frac{1}{m}\sum_{k=1}^{r-1}\mu^{k}\breve{h}$ and 
$\eta=e^{\frac{2\pi i}{r}\sum_{k=0}^{r-1}k({\rm ad}\;\mu^{k}\breve{h})}$. 
Then $\mu h=h$ and $h\in \mathfrak{h}$ since $\mathfrak{h}$ is invariant under $\mu$. Moreover,
\begin{align*}
\eta e^{2 \pi i({\rm ad}\;\breve{h})} \mu \eta^{-1}&=
e^{\frac{2\pi i}{r}\sum_{k=1}^{r-1}k({\rm ad}\;\mu^{k}\breve{h})}
e^{2 \pi i({\rm ad}\;\breve{h})} \mu 
e^{-\frac{2\pi i}{r}\sum_{k=1}^{r-1}k({\rm ad}\;\mu^{k}\breve{h})}\nn
&=e^{\frac{2\pi i}{r}\sum_{k=1}^{r-1}k({\rm ad}\;\mu^{k}\breve{h})}
e^{2 \pi i({\rm ad}\;\breve{h})} 
e^{-\frac{2\pi i}{r}\sum_{k=1}^{r-1}k({\rm ad}\;\mu^{k+1}\breve{h})}\mu\nn
&=e^{2\pi i\frac{1}{r}\sum_{k=1}^{r-1}({\rm ad}\; \mu^{k}\breve{h})}\mu\nn
&=e^{2\pi i({\rm ad}\; h)}\mu.
\end{align*}

Let $\tau_{\sigma}=\eta\nu^{-1}$. Then 
$$\sigma=\nu e^{2 \pi i({\rm ad}\;\breve{h})} \mu  \nu^{-1}
=\nu \eta^{-1}e^{2\pi i({\rm ad}\; h)}\mu \eta\nu^{-1}
=\tau_{\sigma}e^{2\pi i({\rm ad}\; h)}\mu\tau_{\sigma}^{-1} .$$
\epfv

For $\mathcal{N}_{g}$, we have the following result:

\begin{prop}\label{N-g}
Assume that $\mathfrak{g}$ is semisimple. Then we have the following:

\begin{enumerate}

\item There exists $a_{\mathcal{N}_{g}}\in \g^{[0]}$ such that 
$\mathcal{N}_{g}b=[a_{\mathcal{N}_{g}}, b]$ for $b\in \g$, that is,
$\mathcal{N}_{g}={\rm ad}\;a_{\mathcal{N}_{g}}$. 

\item On $\hat{\g}$, $\mathcal{N}_{g}(b\otimes t^{m})=[a_{\mathcal{N}_{g}}\otimes t^{0}, b\otimes t^{m}]
=[a_{\mathcal{N}_{g}}, b]\otimes t^{m}$
for $b\in \g$, $m\in \Z$
and $\mathcal{N}_{g}\mathbf{k}=[a_{\mathcal{N}_{g}}\otimes t^{0}, \mathbf{k}]=0$.

\item On $\hat{\g}^{[g]}$, $\mathcal{N}_{g}(b\otimes t^{m})=[a_{\mathcal{N}_{g}}\otimes t^{0}, b\otimes t^{m}]
=[a_{\mathcal{N}_{g}}, b]\otimes t^{m}$
for $b\in \g^{[\beta]}$ and $m\in \beta+\Z$
and $\mathcal{N}_{g}\mathbf{k}=[a_{\mathcal{N}_{g}}\otimes t^{0}, \mathbf{k}]=0$.

\item On $M(\ell, 0)$ or $L(\ell, 0)$, 
$\mathcal{N}_{g}=a_{\mathcal{N}_{g}}(0)$.

\end{enumerate}
\end{prop}
\pf
By Corollary \ref{N-g-der}, $\mathcal{N}_{g}$ is a derivation of $\g$. Since $\g$ is semisimple,
we know that every derivation of $\g$ is inner. So there exists 
$a_{\mathcal{N}_{g}}\in \g$ such that 
$\mathcal{N}_{g}b=[a_{\mathcal{N}_{g}}, b]$ for $b\in \g$. 
Since $\mathcal{S}_{g}$ commutes with  $\mathcal{N}_{g}$, for $b\in \g^{[\beta]}$,
$\mathcal{S}_{g}[a_{\mathcal{N}_{g}}, b]=\mathcal{S}_{g}\mathcal{N}_{g}b
=\mathcal{N}_{g}\mathcal{S}_{g}b=\beta \mathcal{N}_{g}b=\beta [a_{\mathcal{N}_{g}}, b]$.
So $[a_{\mathcal{N}_{g}}, b]\in \g^{[\beta]}$. Thus $a_{\mathcal{N}_{g}}\in \g^{[0]}$. This finishes
the proof of Conclusion 1.

Conclusions 2, 3 and 4 follow immediately from the definitions of the actions of $\mathcal{N}_{g}$
on $\hat{\g}$, $\hat{\g}^{[g]}$, $M(\ell, 0)$ and $L(\ell, 0)$. 
\epfv

In the rest of this section,  we assume that $\g$ is simple with a Cartan subalgebra $\mathfrak{h}$ and 
a set $\Delta$ of simple roots which gives a root system $\Phi$. 
For $a\in \g$, $a=\sum_{\alpha\in P_{g}}a^{\alpha}$,
where $a^{\alpha}\in \g^{[\alpha]}$. Given a
$\hat{\g}$-module $W$, we have introduced $a^{\alpha}(x)=\sum_{n\in \alpha+\Z}a(n)x^{-n-1}$
for $\alpha\in P_{\g}$ above. 

We shall need the following result later:

\begin{prop}\label{comm=0}
Assume that $\g$ is simple with a given Cartan subalgebra $\mathfrak{h}$ and 
a given set $\Delta$ of simple roots. Let $g$, $\mu$, $h$ and $\tau_{\sigma}$ be the same as 
in Proposition \ref{autom-decomp}. Let $\theta$ be the highest root of $\g$ 
and $W$ a $\hat{\g}^{[g]}$-module. Then there exists $r_{\theta}\in \Z$ such that 
$e_{\theta}\in \g^{[(\theta, h) +r_{\theta}]}$ and 
\begin{equation}\label{comm=0-0}
[(\tau_{\sigma}e_{\theta})(x_{1}), (\tau_{\sigma}e_{\theta})(x_{2})]=0.
\end{equation}
\end{prop}
\pf
Since $({\rm ad}\;h)e_{\theta}=[h, \theta]=(\theta, h)e_{\theta}$, 
$e^{2 \pi i({\rm ad}\;h)}e_{\theta}=e^{2\pi i(\theta, h)}e_{\theta}$. Then 
$$e^{2 \pi i({\rm ad}\;\tau_{\sigma}h)}\tau_{\sigma}e_{\theta}
=\tau_{\sigma}e^{2 \pi i({\rm ad}\;h)}e_{\theta}
=e^{2\pi i(\theta, h)}\tau_{\sigma}e_{\theta}.$$
Since $\theta$ is the highest root and $\mu$ is an automorphism of $\Phi_{+}$, 
$\theta$ is fixed under $\mu$ by the definition and by the  uniqueness of the highest root. 
Thus $e_{\theta}$ is also fixed under $\mu$. So $\tau_{\sigma} \mu \tau_{\sigma}^{-1}$ fixes
$\tau_{\sigma}e_{\theta}$. Then by Proposition \ref{autom-decomp}, 
\begin{equation}\label{comm=0-1}
\sigma\tau_{\sigma} e_{\theta}=\tau_{\sigma}e^{2 \pi i({\rm ad}\;h)}\mu e_{\theta}
=e^{2\pi i(\theta, h)}\tau_{\sigma}e_{\theta}.
\end{equation}
%Thus 
%\begin{equation}\label{comm=0-2}
%(\sigma\tau_{\sigma}^{-1} e_{\theta}, \tau_{\sigma}^{-1} e_{\theta})
%=e^{2\pi i(\theta, h)}(\tau_{\sigma}^{-1} e_{\theta}, \tau_{\sigma}^{-1} e_{\theta})=0.
%\end{equation}
From $\sigma=e^{2\pi i\mathcal{S}_{g}}$ and (\ref{comm=0-1}), 
there exists $r_{\theta}\in \Z$ such that $(\theta, h) +r_{\theta}\in P_{\g}$ and
\begin{equation}\label{comm=0-2}
\mathcal{S}_{g}\tau_{\sigma} e_{\theta}=((\theta, h) +r_{\theta})
\tau_{\sigma} e_{\theta}. 
\end{equation}
Thus $e_{\theta}\in \g^{[(\theta, h) +r_{\theta}]}$. 

To prove (\ref{comm=0-0}), first we have
\begin{align}\label{main-2}
[\tau_{\sigma}e_{\theta}, \tau_{\sigma}e_{\theta}]=0.
\end{align}
Since $\theta+\theta\ne 0$, $(e_{\theta}, e_{\theta})=0$. Then by the invariance of 
the bilinear form $(\cdot, \cdot)$, 
\begin{equation}\label{main-3}
(\tau_{\sigma}e_{\theta}, \tau_{\sigma}e_{\theta})=0.
\end{equation} 
Using the invariance of the bilinear form $(\cdot, \cdot)$ and $\mathcal{N}_{g}={\rm ad}\;a_{\mathcal{N}_{g}}$
(Part  1 of Proposition \ref{N-g}), 
we obtain 
\begin{align}\label{main-4}
(\mathcal{N}_{g}\tau_{\sigma}e_{\theta}, \tau_{\sigma}e_{\theta})=([a_{\mathcal{N}_{g}}, \tau_{\sigma}e_{\theta}], \tau_{\sigma}e_{\theta})
=(a_{\mathcal{N}_{g}}, [\tau_{\sigma}e_{\theta}, \tau_{\sigma}e_{\theta}])=0.
\end{align}

Let $\gamma=(\theta, h) +r_{\theta}$. Then  $e_{\theta}\in \g^{[\gamma]}$.
Using (\ref{main-2}), (\ref{main-3}) and (\ref{main-4}), we have
\begin{align}\label{comm=0-3}
[(\tau&_{\sigma}e_{\theta})(x_{1}), (\tau_{\sigma}e_{\theta})(x_{2})]\nn
&=\sum_{m\in \gamma+\Z}\sum_{n\in \gamma+\Z}
[(\tau_{\sigma}e_{\theta})(m), (\tau_{\sigma}e_{\theta})(n)]\nn
&=\sum_{m\in \gamma+\Z}\sum_{n\in \gamma+\Z}
\Bigl([(\tau_{\sigma}e_{\theta}), (\tau_{\sigma}e_{\theta})](m+n)\nn
&\quad\quad\quad\quad\quad\quad\quad\quad\quad
+m((\tau_{\sigma}e_{\theta}), (\tau_{\sigma}e_{\theta}))
\delta_{m+n, 0}\ell+(\mathcal{N}_{g}(\tau_{\sigma}e_{\theta}), (\tau_{\sigma}e_{\theta}))
\delta_{m+n, 0}\ell\Bigr)\nn
&=0.
\end{align}
\epfv

We also need the following general lemma:

\begin{lemma}\label{prod-same-x}
Let $V$ be a grading-restricted vertex algebra (or a vertex operator algebra), $g$ an automorphism of $V$ and 
$W$ a lower-bounded generalized $g$-twisted $V$-module. Assume that for some $u, v\in V$,
$$(Y_{W}^{g})_{0}(u, x_1)
(Y_{W}^{g})_{0}(v, x_2)=(Y_{W}^{g})_{0}(v, x_2)
(Y_{W}^{g})_{0}(u, x_1),$$
where $(Y_{W}^{g})_{0}(v, x)$ for $v\in V$ is the constant term of $Y_{W}^{g}(v, x)$
viewed as a power series in $\log x$. 
Then $(Y_{W}^{g})_{0}(u, x)
(Y_{W}^{g})_{0}(v, x)$ is well defined and 
\begin{equation}\label{L-mod-ss-6}
(Y_{W}^{g})_{0}((Y_{V})_{-1}(u)v, x))=
(Y_{W}^{g})_{0}(u, x)
(Y_{W}^{g})_{0}(v, x).
\end{equation}
\end{lemma}
\pf
For $w\in W$, 
$$Y_{W}^{g}(u, x_1)
Y_{W}^{g}(v, x_2)w=Y_{W}^{g}(v, x_2)Y_{W}^{g}(u, x_1)w$$ 
has only finitely many 
negative power terms in both $x_{1}$ and $x_{2}$. In particular, we can let 
$x_{1}=x_{2}=x$ to obtain a well defined formal series $(Y_{W}^{g})_{0}(u, x)
(Y_{W}^{g})_{0}(v, x)$.

To prove (\ref{L-mod-ss-6}), we use the the Jacobi identity (\ref{jacobi}).
Using $Y_{W}^{g}(u, x)=(Y_{W}^{g})_{0}(x^{-\mathcal{N}_{g}}u, x)$ 
((2.10) in \cite{HY}) for $u\in V$
and $x_{2}^{\mathcal{N}_{g}}Y_{V}(u, x_{0})=Y_{V}(x_{2}^{\mathcal{N}_{g}}u, x)x_{2}^{\mathcal{N}_{g}}$
((2.5) in \cite{H-twist-vo}), and
replacing $u$ and $v$ in (\ref{jacobi}) by $x_{1}^{\mathcal{N}_{g}}u$ and $x_{2}^{\mathcal{N}_{g}}v$,
respectively, we see that that (\ref{jacobi}) becomes
\begin{align}\label{jacobi-1}
x_0^{-1}&\delta\left(\frac{x_1 - x_2}{x_0}\right)
(Y_{W}^{g})_{0}(u, x_1)
(Y_{W}^{g})_{0}(v, x_2)- x_0^{-1}\delta\left(\frac{- x_2 + x_1}{x_0}\right)
(Y_{W}^{g})_{0}(v, x_2)
(Y_{W}^{g})_{0}(u, x_1)\nn
&= x_1^{-1}\delta\left(\frac{x_2+x_0}{x_1}\right)
(Y_{W}^{g})_{0}\left(Y_{V}\left(\left(\frac{x_{2}}{x_{1}}\right)^{\mathcal{S}_{g}}
\left(1+\frac{x_0}{x_2}\right)^{\mathcal{L}_{g}}
u, x_0\right)v, x_2\right).
\end{align}
By the assumption, the left-hand side of (\ref{jacobi-1}) is equal to 
\begin{align*}
&\left(x_0^{-1}\delta\left(\frac{x_1 - x_2}{x_0}\right)-
x_0^{-1}\delta\left(\frac{- x_2 + x_1}{x_0}\right)\right)(Y_{W}^{g})_{0}(u, x_1)
(Y_{W}^{g})_{0}(v, x_2)\nn
&\quad\quad\quad\quad\quad  =x_1^{-1}\delta\left(\frac{x_2+x_0}{x_1}\right)(Y_{W}^{g})_{0}(u, x_1)
(Y_{W}^{g})_{0}(v, x_2).
\end{align*}
Thus from (\ref{jacobi-1}), we obtain
\begin{align}\label{L-mod-ss-2}
 x_1^{-1}&\delta\left(\frac{x_2+x_0}{x_1}\right)
(Y_{W}^{g})_{0}\left(Y_{V}\left(\left(\frac{x_{2}}{x_{1}}\right)^{\mathcal{S}_{g}}
\left(1+\frac{x_0}{x_2}\right)^{\mathcal{L}_{g}}
u, x_0\right)v, x_2\right)\nn
&\quad\quad\quad\quad=x_1^{-1}\delta\left(\frac{x_2+x_0}{x_1}\right)Y_{W}^{g}(u, x_1)
Y_{W}^{g}(v, x_2).
\end{align}
Replacing $u$ in  (\ref{L-mod-ss-2}) by $\bigl(1+\frac{x_0}{x_2}\bigr)^{-\mathcal{L}_{g}}
\bigl(\frac{x_{2}}{x_{1}}\bigr)^{-\mathcal{S}_{g}}u$, we obtain 
\begin{align}\label{L-mod-ss-3}
x_1^{-1}&\delta\left(\frac{x_2+x_0}{x_1}\right)
(Y_{W}^{g})_{0}\left(Y_{V}\left(u, x_0\right)v, x_2\right)\nn
&=x_1^{-1}\delta\left(\frac{x_2+x_0}{x_1}\right)
(Y_{W}^{g})_{0}\left(\left(1 + \frac{x_0}{x_2}\right)^{-\mathcal{L}_{g}}
\left(\frac{x_{2}}{x_{1}}\right)^{-\mathcal{S}_{g}}u, x_1\right)
(Y_{W}^{g})_{0}(v, x_2).
\end{align}
Since $V=\coprod_{\alpha\in P_{V}}V^{[\alpha]}$,
we have $u=\sum_{\alpha\in P_{V}}u^{\alpha}$, where $u^{\alpha}\in V^{[\alpha]}$ for $\alpha\in P_{V}$. 
Also note that $(Y_{W}^{g})_{0}(u^{\alpha}, x)\in x^{-\alpha}({\rm End}\;W)[[x, x^{-1}]]$. 
Then we have 
\begin{align}\label{L-mod-ss-4}
x_1^{-1}&\delta\left(\frac{x_2+x_0}{x_1}\right)
(Y_{W}^{g})_{0}\left(\left(1 + \frac{x_0}{x_2}\right)^{-\mathcal{L}_{g}}
\left(\frac{x_{2}}{x_{1}}\right)^{-\mathcal{S}_{g}}u, x_1\right)\nn
&=\sum_{\alpha\in P_{V}}x_1^{-1}\delta\left(\frac{x_2+x_0}{x_1}\right)
x_{1}^{\alpha}(Y_{W}^{g})_{0}\left(\left(1 + \frac{x_0}{x_2}\right)^{-\mathcal{L}_{g}}
x_{2}^{-\mathcal{S}_{g}}u^{\alpha}, x_1\right)\nn
&=\sum_{\alpha\in P_{V}}x_1^{-1}\delta\left(\frac{x_2+x_0}{x_1}\right)
(x_2+x_0)^{\alpha}(Y_{W}^{g})_{0}\left(\left(1 + \frac{x_0}{x_2}\right)^{-\mathcal{L}_{g}}
x_{2}^{-\mathcal{S}_{g}}u^{\alpha}, x_2+x_0\right)\nn
&=x_1^{-1}\delta\left(\frac{x_2+x_0}{x_1}\right)
(Y_{W}^{g})_{0}\left(\left(1 + \frac{x_0}{x_2}\right)^{-\mathcal{L}_{g}}
\left(\frac{x_{2}}{x_{2}+x_{0}}\right)^{-\mathcal{S}_{g}}u, x_2+x_0\right)\nn
&=x_1^{-1}\delta\left(\frac{x_2+x_0}{x_1}\right)
(Y_{W}^{g})_{0}\left(\left(1 + \frac{x_0}{x_2}\right)^{-\mathcal{N}_{g}}u, x_2+x_0\right).
\end{align}
Using (\ref{L-mod-ss-4}), we see that the right-hand side of (\ref{L-mod-ss-3})
is equal to 
\begin{align}\label{L-mod-ss-5}
x_1^{-1}\delta\left(\frac{x_2+x_0}{x_1}\right)
(Y_{W}^{g})_{0}\left(\left(1 + \frac{x_0}{x_2}\right)^{-\mathcal{N}_{g}}u, x_2+x_{0}\right)
(Y_{W}^{g})_{0}(v, x_2).
\end{align}
Then $\res_{x_{1}}$ of the left-hand side of (\ref{L-mod-ss-3}) and (\ref{L-mod-ss-5}) are also equal, that is,
\begin{equation}\label{L-mod-ss-5.5}
(Y_{W}^{g})_{0}\left(Y_{V}\left(u, x_0\right)v, x_2\right)=
(Y_{W}^{g})_{0}\left(\left(1 + \frac{x_0}{x_2}\right)^{-\mathcal{N}_{g}}u, x_2+x_{0}\right)
Y_{W}^{g}(v, x_2).
\end{equation}
Taking the constant terms in $x_{0}$ in both sides of (\ref{L-mod-ss-5.5}) and then replacing $x_{2}$ by $x$, 
we obtain (\ref{L-mod-ss-6}).
\epfv

Applying Proposition \ref{comm=0} and 
Lemma \ref{prod-same-x} to the lower-bounded (grading-restricted) 
generalized $g$-twisted $M(\ell, 0)$-module
$\overarc{M}_{\ell}^{[g]}$ ($\widetilde{M}_{\ell}^{[g]}$) and using Theorems \ref{ulbtm-voa}
and \ref{first-main}, we have the following consequence:

\begin{cor}
On the lower-bounded $\hat{\g}^{[g]}$-module $\overarc{M}_{\ell}^{[g]}$ 
(the grading-restricted $\hat{\g}^{[g]}$-module $\widetilde{M}_{\ell}^{[g]}$), 
$(\tau_{\sigma}e_{\theta})(x)^{m}$ for $m\in \N$ are well defined
and 
\begin{equation}\label{main-7}
Y_{\overarc{M}_{\ell}^{[g]}}^{g}((\tau_{\sigma}e_{\theta})(-1)^{m}\one, x)
=(\tau_{\sigma}e_{\theta})(x)^{m}
\end{equation}
\begin{equation}\label{main-7.5}
\left(Y_{\widetilde{M}_{\ell}^{[g]}}^{g}((\tau_{\sigma}e_{\theta})(-1)^{m}\one, x)
=(\tau_{\sigma}e_{\theta})(x)^{m}\right).
\end{equation}
In particular,  on
$U(\hat{\g}^{[g]})\otimes_{G(\Omega^{[g]}, \hat{\g}_{+}^{[g]}, \mathbf{k})}M$ and
$U(\hat{\g}^{[g]})
\otimes_{\hat{\g}^{[g]}_{+}\oplus \hat{\g}^{[g]}_{0}}M$,
$(\tau_{\sigma}e_{\theta})(x)^{m}$ for $m\in \N$ are well defined. 
\end{cor}
\pf
By the definition of $\hat{\g}^{[g]}$-module structure on $\overarc{M}_{\ell}^{[g]}$
(see  Proposition \ref{twisted-affine-mod}), we have 
$$(Y_{\overarc{M}_{\ell}^{[g]}}^{g})_{0}((\tau_{\sigma}e_{\theta})(-1)\one, x)
=(\tau_{\sigma}e_{\theta})(x).$$
Then from (\ref{comm=0-0}), we have
\begin{align}\label{L-mod-ss-1}
(Y_{\overarc{M}_{\ell}^{[g]}}^{g})_{0}&((\tau_{\sigma}e_{\theta})(-1)\one, x_{1})
(Y_{\overarc{M}_{\ell}^{[g]}}^{g})_{0}((\tau_{\sigma}e_{\theta})(-1)\one, x_{2})\nn
& =
(Y_{\overarc{M}_{\ell}^{[g]}}^{g})_{0}((\tau_{\sigma}e_{\theta})(-1)\one, x_{2})
(Y_{\overarc{M}_{\ell}^{[g]}}^{g})_{0}((\tau_{\sigma}e_{\theta})(-1)\one, x_{1}).
\end{align}
By (\ref{L-mod-ss-1}),
we can use Lemma \ref{prod-same-x} for $u=v=(\tau_{\sigma}e_{\theta})(-1)$. 
So we have 
$$Y_{\overarc{M}_{\ell}^{[g]}}^{g}((\tau_{\sigma}e_{\theta})(-1)^{2}\one, x)
=(\tau_{\sigma}e_{\theta})(x)^{2},$$
where the right-hand side is well defined.
From (\ref{comm=0-0}), we obtain
$$[(\tau_{\sigma}e_{\theta})(x), (\tau_{\sigma}e_{\theta})(x)^{m}]
=\sum_{i=1}^{m}(\tau_{\sigma}e_{\theta})(x)^{i-1}
[(\tau_{\sigma}e_{\theta})(x), (\tau_{\sigma}e_{\theta})(x)]
(\tau_{\sigma}e_{\theta})(x)^{m-i}=0.$$
Using induction and Lemma \ref{prod-same-x}, we see that for 
$m\in \N$, $(\tau_{\sigma}e_{\theta})(x)^{m}$
is well defined and (\ref{main-7}) holds. The proof for $\widetilde{M}_{L(\ell, 0)}^{[g]}$ is completely the same. 

By  Theorems \ref{ulbtm-voa}
and \ref{first-main}, we see that $(\tau_{\sigma}e_{\theta})(x)^{m}$ for $m\in \N$ are also well defined
on $U(\hat{\g}^{[g]})\otimes_{G(\Omega^{[g]}, \hat{\g}_{+}^{[g]}, \mathbf{k})}M$ and
$U(\hat{\g}^{[g]})
\otimes_{\hat{\g}^{[g]}_{+}\oplus \hat{\g}^{[g]}_{0}}M$.
\epfv

Let $V$ be a grading-restricted vertex algebra (or a vertex operator algebra), $g$ and $h$ automorphisms 
of $V$ and $(W, Y^{g}_{W})$ a  lower-bounded (grading-restricted) generalized $g$-twisted
$V$-module. Recall the  lower-bounded (grading-restricted)
$hgh^{-1}$-twisted
$V$-module $(W, \phi_{h}(Y^{g}))$ (see Proposition 3.2 in \cite{H-twisted-int}),
where 
\begin{align*}
\phi_{h}(Y^{g}_{W}): V\times W&\to W\{x\}[{\rm log} x]\nn
v \otimes w & \mapsto  \phi_{h}(Y^{g}_{W})(v, x)w
\end{align*}
is the linear map defined by
$ \phi_{h}(Y^{g}_{W})(v, x)w=Y^{g}_{W}(h^{-1}v, x)w.$
The $hgh^{-1}$-twisted
$V$-module $(W, \phi_{h}(Y^{g}_{W}))$ is also  denoted by $\phi_{h}(W)$.
%Similarly, for a $\hat{\g}^{[g]}$-module $W$ and an automorphism $h$ of $\g$, we have a 
%$\hat{\g}^{[hgh^{-1}]}$-module. 

Let $M$, as in Subsection \ref{4.2}, be a vector space with actions of 
$g$,  $\mathcal{L}_{g}$,
$\mathcal{S}_{g}$, $\mathcal{N}_{g}$, $L_{M}(0)$, $L_{M}(0)_{S}$ and $L_{M}(0)_{N}$ such that
$M=\coprod_{h\in Q_{M}}M_{[h]}$, where $Q_{M}$ is the set of all eigenvalues of $L_{M}(0)$ and 
$M_{[h]}$ is the generalized eigenspace of $L_{M}(0)$ with eigenvalue $h\in Q_{M}$. 
Then just as in the construction of $\overarc{M}_{\ell}^{[g]}$ 
in Subsection \ref{4.2} for the vertex operator algebra $M(\ell, 0)$, for each $h\in Q_{M}$, 
we have a universal lower-bounded generalized $g$-twisted $L(\ell, 0)$-module generated by $M_{h}$ (see 
Theorem \ref{universal-voa} for its universal property). 
We shall denote them by $\overarc{M}_{L(\ell, 0), h}^{[g]}$ for $h\in Q_{M}$.
Let 
$\overarc{M}_{L(\ell, 0)}^{[g]}=\coprod_{h\in Q_{M}}\overarc{M}_{L(\ell, 0), h}^{[g]}.$
This is the universal lower-bounded generalized $g$-twisted $L(\ell, 0)$-module generated by 
a subspace $M$ annihilated by $\hat{\g}_{+}$. 
If $M$ is in addition a finite-dimensional $\hat{\g}_{\I}$-module such that the action of 
$g$ is compatible, then also as in  Subsection \ref{4.2}, we have a quotient of $\overarc{M}_{L(\ell, 0), h}^{[g]}$. 
We shall denote this quotient by $\widetilde{M}_{L(\ell, 0)}^{[g]}$.

From the definitions of $\overarc{M}_{L(\ell, 0)}^{[g]}$ and $\widetilde{M}_{L(\ell, 0)}^{[g]}$,
we have the following universal properties for them:

\begin{thm}\label{univ-property-L}
Let $(W, Y^{g}_{W})$ be a lower-bounded generalized $g$-twisted $L(\ell, 0)$-module
(when $L(\ell, 0)$ is viewed as a vertex operator algebra).
Let $M^{0}$ a subspace (finite-dimensional $\hat{\g}^{[g]}_{\I}$-submodule) of $W$ invariant 
under the actions of 
$g$, $\mathcal{S}_{g}$, $\mathcal{N}_{g}$, $L_{W}(0)$, $L_{W}(0)_{S}$ and
$L_{W}(0)_{N}$ and annihilated by $\hat{\g}^{[g]}_{+}$. 
Assume that there is a linear ($\hat{\g}^{[g]}_{\I}$-module map) $f: M\to M^{0}$ 
commuting with the actions of 
$g$, $\mathcal{S}_{g}$, $\mathcal{N}_{g}$, $L_{W}(0)|_{M^{0}}$ and $L_{M}(0)$, 
$L_{W}(0)_{S}|_{M^{0}}$ and $L_{M}(0)_{S}$ and
$L_{W}(0)_{N}|_{M^{0}}$ and $L_{M}(0)_{N}$. Then there exists a unique module 
map $\overarc{f}: \overarc{M}^{[g]}_{L(\ell, 0)}\to W$ $(\tilde{f}: \widetilde{M}^{[g]}_{L(\ell, 0)}\to W)$
such that $\overarc{f}|_{M}=f$ $(\tilde{f}|_{M}=f)$. 
If $f$ is surjective and $(W, Y^{g}_{W})$ is generated by the 
coefficients of $(Y^{g})_{WL(\ell, 0)}^{W}(w, x)\one$ for $w\in M^{0}$, 
where $(Y^{g})_{WL(\ell, 0)}^{W}$ 
is the twist vertex operator map obtained from $Y_{W}^{g}$, then 
$\overarc{f}$ $(\tilde{f})$ is surjective. 
\end{thm}
\pf
For $h\in Q_{M}$, by the universal property of 
$\overarc{M}_{L(\ell, 0), h}^{[g]}$ (see the construction of $\overarc{M}_{L(\ell, 0), h}^{[g]}$ and 
Theorem \ref{universal-voa}, ), there is a unique $L(\ell, 0)$-module map $\overarc{f}_{h}$
from $\overarc{M}_{L(\ell, 0), h}^{[g]}$ to the submodule of $W$ generated by $f(M^{0})$.
Let $\overarc{f}: \overarc{M}^{[g]}_{L(\ell, 0)}\to W$ be defined by $\overarc{f}(w)=\overarc{f}_{h}(w)$
for $w\in \overarc{M}_{L(\ell, 0), h}^{[g]}$. Then $\overarc{f}$ is clearly a module 
map. The uniqueness of $\overarc{f}$ follows from the uniqueness of $\overarc{f}_{h}$ for $h\in Q_{M}$. 
It is also clear that the second conclusion holds. 

In the case that $M^{0}$ is a finite-dimensional $\hat{\g}^{[g]}_{\I}$-submodule of $W$, 
the proof is the same as that of Theorem \ref{univ-property} 
 except that we should 
use $\overarc{M}^{[g]}_{L(\ell, 0)}$ instead of $\widehat{M}^{[g]}_{\ell}$. 
\epfv

Let $\overarc{I}_{L(\ell, 0)}^{[g]}$ ($\widetilde{I}_{L(\ell, 0)}^{[g]}$)  be the submodules of 
$U(\hat{\g}^{[g]})\otimes_{G(\Omega^{[g]}, \hat{\g}_{+}^{[g]}, \mathbf{k})}M$ 
($U(\hat{\g}^{[g]})\otimes_{\hat{\g}^{[g]}_{+}\oplus \hat{\g}^{[g]}_{0}}M$)
generated by the coefficients of $(\tau_{\sigma}e_{\theta})(x)^{\ell+1}w$ for 
$w\in U(\hat{\g}^{[g]})\otimes_{G(\Omega^{[g]}, \hat{\g}_{+}^{[g]}, \mathbf{k})}M$
($w\in U(\hat{\g}^{[g]})\otimes_{\hat{\g}^{[g]}_{+}\oplus \hat{\g}^{[g]}_{0}}M$). 
Then on $(U(\hat{\g}^{[g]})\otimes_{G(\Omega^{[g]}, \hat{\g}_{+}^{[g]}, \mathbf{k})}M)/\overarc{I}_{L(\ell, 0)}^{[g]}$,
$(U(\hat{\g}^{[g]})\otimes_{\hat{\g}^{[g]}_{+}\oplus \hat{\g}^{[g]}_{0}}M)/\widetilde{I}_{L(\ell, 0)}^{[g]}$
and their quotients, 
\begin{equation}\label{main-0}
(\tau_{\sigma}e_{\theta})(x)^{\ell+1}=0.
\end{equation}

\begin{thm}\label{main}
Assume that $\g$ is simple and $\ell\in \Z_{+}$.
The universal lower-bounded (grading-restricted) 
generalized $g$-twisted $L(\ell, 0)$-module $\overarc{M}_{L(\ell, 0)}^{[g]}$
($\widetilde{M}_{L(\ell, 0)}^{[g]}$)
is equivalent as a lower-bounded (grading-restricted) $\hat{\g}^{[g]}$-module to 
$(U(\hat{\g}^{[g]})\otimes_{G(\Omega^{[g]}, \hat{\g}_{+}^{[g]}, \mathbf{k})}M)/\overarc{I}_{L(\ell, 0)}^{[g]}$
$((U(\hat{\g}^{[g]})\otimes_{\hat{\g}^{[g]}_{+}\oplus \hat{\g}^{[g]}_{0}}M)/\widetilde{I}_{L(\ell, 0)}^{[g]})$. 
In particular, the lower-bounded generalized $g$-twisted $L(\ell, 0)$-module 
$\widetilde{M}_{L(\ell, 0)}^{[g]}$ is in fact grading restricted.
\end{thm}
\pf
Note that the automorphism $\tau_{\sigma}$ of $\g$ induces
automorphisms, denoted still by $\tau_{\sigma}$, of the vertex operator algebras $M(\ell, 0)$ and $L(\ell, 0)$. 
Then $\tau_{\sigma}(\overarc{M}_{L(\ell, 0)}^{[g]})$ is a lower-bounded generalized
$\tau_{\sigma}^{-1}g\tau_{\sigma}$-twisted $L(\ell, 0)$-module. 
By Proposition \ref{L-mod-t-aff-mod}, $\tau_{\sigma}(\overarc{M}_{L(\ell, 0)}^{[g]})$ is 
also a lower-bounded generalized $\tau_{\sigma}^{-1}g\tau_{\sigma}$-twisted $M(\ell, 0)$-module. 
By Proposition \ref{L-mod-t-aff-mod-1}, We have 
$$\phi_{\tau_{\sigma}}\left(Y_{\overarc{M}_{L(\ell, 0)}^{[g]}}^{g}\right)(e_{\theta}(-1)^{\ell+1}\one, x)=0,$$
or equivalently,
$$Y^{g}_{\overarc{M}_{L(\ell, 0)}^{[g]}}(\tau_{\sigma}e_{\theta}(-1)^{\ell+1}\one, x)=0.$$
Since $\tau_{\sigma}$ is an automorphism of $L(\ell, 0)$ induced from the automorphism $\tau_{\sigma}$
of $\g$, 
we have $\tau_{\sigma}(e_{\theta}(-1)^{\ell+1})\one=(\tau_{\sigma}e_{\theta})(-1)^{\ell+1}\one$.
Hence we have 
\begin{equation}\label{main-6}
Y^{g}_{\overarc{M}_{L(\ell, 0)}^{[g]}}((\tau_{\sigma}e_{\theta})(-1)^{\ell+1}\one, x)=0.
\end{equation}
From (\ref{main-6}) and (\ref{main-7}), we see that (\ref{main-0}) holds
on the $\hat{\g}^{[g]}$-module $\overarc{M}_{L(\ell, 0)}^{[g]}$. 
By Proposition \ref{L-mod-t-aff-mod}, $\overarc{M}_{L(\ell, 0)}^{[g]}$ is also a 
lower-bounded generalized $g$-twisted $M(\ell, 0)$-module generated by $M$.
Then by Corollary \ref{univ-cor-lb}, $\overarc{M}_{L(\ell, 0)}^{[g]}$ is equivalent to a quotient of 
$U(\hat{\g}^{[g]})\otimes_{G(\Omega^{[g]}, \hat{\g}_{+}^{[g]}, \mathbf{k})}M$.
Since on $\overarc{M}_{L(\ell, 0)}^{[g]}$ (\ref{main-0}) holds, we see that 
$\overarc{M}_{L(\ell, 0)}^{[g]}$ is equivalent to a quotient of 
$(U(\hat{\g}^{[g]})\otimes_{G(\Omega^{[g]}, \hat{\g}_{+}^{[g]}, \mathbf{k})}M)/\overarc{I}_{L(\ell, 0)}^{[g]}$.

On the other hand, by Theorem \ref{ulbtm-voa}, 
$(U(\hat{\g}^{[g]})\otimes_{G(\Omega^{[g]}, \hat{\g}_{+}^{[g]}, \mathbf{k})}M)/\overarc{I}_{L(\ell, 0)}^{[g]}$
is equivalent as a $\hat{\g}$-module to the lower-bounded generalized $g$-twisted $M(\ell, 0)$-module 
$\overarc{M}_{\ell}^{[g]}$. Then by Corollary \ref{univ-cor-lb}, 
$W=(U(\hat{\g}^{[g]})\otimes_{G(\Omega^{[g]}, \hat{\g}_{+}^{[g]}, \mathbf{k})}M)/\overarc{I}_{L(\ell, 0)}^{[g]}$ 
is also a lower-bounded generalized $g$-twisted $M(\ell, 0)$-module. But by definition,
(\ref{main-0}) holds on $W$. So we have
\begin{align*}
\phi_{\tau_{\sigma}}(Y_{W}^{g})(e_{\theta}(-1)^{\ell+1}\one, x)
&=Y^{g}_{W}(\tau_{\sigma}e_{\theta}(-1)^{\ell+1}\one, x)\nn
&=Y_{W}^{g}((\tau_{\sigma}e_{\theta})(-1)^{\ell+1}\one, x)\nn
&=(\tau_{\sigma}e_{\theta})(x)^{\ell+1}\nn
&=0.
\end{align*}
By Proposition \ref{L-mod-t-aff-mod-1}, $\phi_{\tau_{\sigma}}(W)$ is a lower-bounded generalized
$\tau_{\sigma}^{-1}g\tau_{\sigma}$-twisted $L(\ell, 0)$-module and thus 
$W$ is a lower-bounded generalized $g$-twisted 
$L(\ell, 0)$-module. Also $M$ can be viewed as a subspace of $W$ invariant 
under the actions of 
$g$, $\mathcal{S}_{g}$, $\mathcal{N}_{g}$, $L_{W}(0)$, $L_{W}(0)_{S}$ and
$L_{W}(0)_{N}$ and with $\hat{\g}^{[g]}_{+}$ acting on $M$ as $0$ and 
we have the identity map  from $M$ to itself. Thus by Theorem \ref{univ-property-L},
there exists a unique surjective $L(\ell, 0)$-module map from $\overarc{M}_{L(\ell, 0)}^{[g]}$
to $W$. In particular, this surjective $L(\ell, 0)$-module map is  a surjective $\hat{\g}$-module map.
Thus we have a surjective $\hat{\g}$-module map from $\overarc{M}_{L(\ell, 0)}^{[g]}$ to 
$W$. Since we have proved that $\overarc{M}_{L(\ell, 0)}^{[g]}$ is a quotient of $W$, 
the existence of such a surjective $\hat{\g}$-module map means that 
$\overarc{M}_{L(\ell, 0)}^{[g]}$ is equivalent to 
$W=(U(\hat{\g}^{[g]})\otimes_{G(\Omega^{[g]}, \hat{\g}_{+}^{[g]}, \mathbf{k})}M)/\overarc{I}_{L(\ell, 0)}^{[g]}$.

The proof for $\widetilde{M}_{L(\ell, 0)}^{[g]}$ is completely the same except that we use the results
in Subsections \ref{4.2} and \ref{4.3} on $\widetilde{M}_{\ell}^{[g]}$ instead of $\overarc{M}_{\ell}^{[g]}$.
Since $(U(\hat{\g}^{[g]})\otimes_{\hat{\g}^{[g]}_{+}\oplus \hat{\g}^{[g]}_{0}}M)/\widetilde{I}_{L(\ell, 0)}^{[g]}$ 
is grading-restricted, we see that $\widetilde{M}_{L(\ell, 0)}^{[g]}$ is grading restricted.
\epfv

We also have the following immediate consequence whose proof is also omitted:

\begin{cor}%\label{univ-cor-gr-L}
Let $W$ be a lower-bounded generalized $g$-twisted $L(\ell, 0)$-module (when $L(\ell, 0)$ is viewed
as a vertex operator algebra) generated by a subspace (finite-dimensional 
$\hat{\g}_{\I}$-submodule) $M$ invariant under $g$, $\mathcal{L}_{g}$, $\mathcal{S}_{g}$,
$\mathcal{N}_{g}$, $L_{W}(0)$, $L_{W}(0)_{S}$ and $L_{W}(0)_{N}$ and annihilated by
$\hat{\g}^{[g]}_{+}$. 
Then $W$ as a lower-bounded $\hat{\g}^{[g]}$-module is equivalent to a quotient of 
$(U(\hat{\g}^{[g]})\otimes_{G(\Omega^{[g]}, \hat{\g}_{+}^{[g]}, \mathbf{k})}M)/\overarc{I}_{L(\ell, 0)}^{[g]}$
$((U(\hat{\g}^{[g]})\otimes_{\hat{\g}^{[g]}_{+}\oplus \hat{\g}^{[g]}_{0}}M)/\widetilde{I}_{L(\ell, 0)}^{[g]}$
and, in particular, 
$W$ is grading restricted).
Conversely, let  $M$ be a vector space (finite-dimensional 
$\hat{\g}_{I_{\I}}$-module) with (compatible) actions of  $g$, $\mathcal{L}_{g}$, $\mathcal{S}_{g}$,
$\mathcal{N}_{g}$, $L_{W}(0)$, $L_{W}(0)_{S}$ and $L_{W}(0)_{N}$ satisfying the conditions
discussed in Section 4. Then 
a quotient module of 
$(U(\hat{\g}^{[g]})\otimes_{G(\Omega^{[g]}, \hat{\g}_{+}^{[g]}, \mathbf{k})}M)/\overarc{I}_{L(\ell, 0)}^{[g]}$
$((U(\hat{\g}^{[g]})\otimes_{\hat{\g}^{[g]}_{+}\oplus \hat{\g}^{[g]}_{0}}M)/\widetilde{I}_{L(\ell, 0)}^{[g]})$
has a natural structure 
of a grading-restricted generalized $g$-twisted $L(\ell, 0)$-module (when $L(\ell, 0)$ is viewed
as a vertex operator algebra). 
\end{cor}

\noindent {\small \sc Department of Mathematics, Rutgers University,
110 Frelinghuysen Rd., Piscataway, NJ 08854-8019}

\noindent {\em E-mail address}: yzhuang@math.rutgers.edu

\end{document}